\documentclass{amsart}

\usepackage[T1]{fontenc}

\usepackage{cite}
\usepackage{tikz}
\usepackage{tikz-cd}
\usepackage{tikz-layers}
\usetikzlibrary{shapes,arrows}
\usetikzlibrary{calc,
	arrows,decorations.pathmorphing,
	backgrounds,fit,positioning,shapes.symbols,chains}
\usetikzlibrary{shapes.geometric, graphs, knots}
\usetikzlibrary{decorations.pathmorphing}
\usetikzlibrary{decorations.markings}
\usetikzlibrary{decorations.pathreplacing,angles,quotes,decorations.markings, arrows.meta,knots}
\usetikzlibrary{positioning,shadows,arrows,trees,fit}
\usetikzlibrary{mindmap}
\usepackage{amssymb,amsmath,amsfonts,latexsym, enumerate}
\usepackage{bm}
\usepackage{yhmath}

\usepackage{longtable, multirow}
\usepackage{multicol}
\usepackage[all,cmtip]{xypic} 
\usepackage{xypic} 
\xyoption{curve}
\usepackage{array, ragged2e}
\newcolumntype{P}[1]{>{\raggedright\arraybackslash}p{#1}}

\usepackage{colortbl}

\usepackage{standalone}
\usepackage{color}
\usepackage{float}
\usepackage{hyperref}
\usepackage[noabbrev, capitalize]{cleveref}
\usepackage{xfrac}

\allowdisplaybreaks

\numberwithin{equation}{section}

\newtheorem{thm}[equation]{Theorem}
\newtheorem{lem}[equation]{Lemma}
\newtheorem{prop}[equation]{Proposition}
\newtheorem{cor}[equation]{Corollary}

\makeatletter
\newcommand*{\centerfloat}{%
	\parindent \z@
	\leftskip \z@ \@plus 1fil \@minus \textwidth
	\rightskip\leftskip
	\parfillskip \z@skip}
\makeatother
\restylefloat*{figure}

\newtheorem{defn}[equation]{Definition}

\theoremstyle{remark}
\newtheorem{rmk}[equation]{Remark}

\newtheorem{ex}[equation]{Example}

\renewcommand{\emptyset}{\font\cmsy = cmsy10 at 10pt
	\hbox{\cmsy \char 59}
}

\definecolor{burntsienna}{rgb}{0.91, 0.45, 0.32}
\definecolor{sapgreen}{rgb}{0.31, 0.49, 0.16}
\definecolor{teal}{rgb}{0.0, 0.5, 0.5}

\definecolor{blue}{rgb}{0.38, 0.51, 0.71} 
\definecolor{darkblue}{RGB}{17, 42, 60} 
\definecolor{red}{RGB}{175, 49, 39} 

\definecolor{orange}{RGB}{217, 156, 55} 
\definecolor{green}{RGB}{144, 169, 84} 
\definecolor{palegreen}{RGB}{197, 184, 104} 

\definecolor{yellow}{RGB}{250, 199, 100} 
\definecolor{brokenwhite}{RGB}{218, 192, 166} 
\definecolor{brokengrey}{rgb}{0.77, 0.76, 0.82} 

\mathchardef\mathdash="2D
\newcommand{\ov }{/}

\newcommand{\defeq}{\overset{\textup{def}}{=}}

\newcommand\R{\mathbb{R}}
\newcommand\N{\mathbb{N}}

\newcommand{\bbb}{\bm{b}}
\newcommand{\ccc}{\bm{c}} 
\newcommand{\ddd}{\bm{d}}

\newcommand{\nul}{\mathbf{0}}

\newcommand{\two}{\mathbf{2}}
\newcommand{\nn}{\mathbf{n}}
\newcommand\I{\mathsf{I}}

\tikzset{
	dot/.style = {circle, fill, minimum size=#1,
		inner sep=0pt, outer sep=0pt},
	dot/.default = 5pt 
}

\tikzset{->-/.style={decoration={
			markings,
			mark=at position #1 with {\arrow{>}}},postaction={decorate}}}
\tikzset{-<-/.style={decoration={
			markings,
			mark=at position #1 with {\arrow{<}}},postaction={decorate}}}

\newcommand\Set{\mathsf{Set}}
\newcommand\Cat{\mathsf {Cat}}

\newcommand\V{\mathsf{X}} 
\newcommand{\CCat}{ \mathsf{C}}
\newcommand{\PCat}{ \mathsf{P}}

\newcommand{\vect}[1][\Bbbk]{ \mathsf{Vect}_{\scriptstyle{#1}}}
\newcommand{\vectf}[1][\Bbbk]{ \mathsf{Vect^f}_{\scriptstyle{#1}}}
\newcommand{\modR}[1][R ]{ {#1}\mathdash\mathsf{Mod}}
\newcommand{\listm} {\small{\mathsf{list}}}
\newcommand\CCC{\mathfrak{C}}
\newcommand\DDD[1][D]{\mathfrak{#1}}
\newcommand{\coev}[1]{{\lfloor{#1}\rfloor}}
\newcommand{\ev}[1]{{\lceil{#1}\rceil}}

\newcommand{\Br}[1][\delta]{\mathsf{Br}_{#1}}

\newcommand{\CBD}[1][(\CCC, \omega)]{\mathsf{BD}^{#1}}
\newcommand{\BD}{\mathsf{BD}}
\newcommand{\BDop}{\mathring{\mathsf{BD}}} 
\newcommand{\CBDd}[1][(\CCC, \omega)]{\mathsf{dBD}^{#1}}

\newcommand{\BDd}{\mathsf{dBD}}

\newcommand{\BDu}{\mathsf{uBD}}
\newcommand{\Sf}[1][f]{\mathfrak S(#1)}
\newcommand{\Tf}[1][f]{\mathfrak T(#1)}
\newcommand\dipal{\{\wideparen{\uparrow,\ \downarrow}\}}
\newcommand\DiBD{\mathsf{OBD}}
\newcommand\DiCBD[1][\DDD]{\mathsf{OBD}^{#1}}

\newcommand\DiBr[1][\delta]{\mathsf{OBr}_{#1}}
\newcommand{\pal}{\mathsf{Pal}}

\newcommand{\kclf}[1][f]{\mathfrak{k}_{#1}}
\newcommand{\kcl}{\mathfrak{k}}

\newcommand{\colop}[1][\lambda]{{{#1}_{\partial}}}
\newcommand{\colcl}[1][\lambda]{{\widetilde{#1}}}

\newcommand{\bfom}{\bm{\overleftarrow \omega}}

\newcommand{\man}{\mathcal M}

\newcommand{\op}{\mathcal O}

\newcommand{\WD}{\mathrm{WD}}
\newcommand{\oWD}{\mathrm{OWD}}

\newcommand{\oCWD}[1][\DDD]{\oWD^{#1}}
\newcommand{\CWD}[1][(\CCC, \omega)]{\WD^{#1}}

\newcommand\WP{\mathsf {WP}}
\newcommand{\CA}{{\mathsf{CA}}}
\newcommand{\CCA}[1][(\CCC,\omega)]{\CA^{#1}}

\newcommand{\alg}[1][A]{\mathcal {#1}}

\newcommand{\xunder}{
	\raisebox{-8pt}{\begin{tikzpicture}[scale = .35]
		
			\begin{knot}[
				clip width=5
				]
				\strand[thick] (0,0) .. controls +(1,0) and +(-1,0) ..
				(2,1);
				\strand [thick] (0,1) .. controls +(1,0) and +(-1,0) ..
				(2,0);
			\end{knot}
\end{tikzpicture}}}
\newcommand{\xover}{
	\raisebox{-8pt}{\begin{tikzpicture}[scale = .35]
			
			\begin{knot}[
				clip width=5
				]
				
				\strand [thick] (0,1) .. controls +(1,0) and +(-1,0) .. 
				(2,0);
				\strand [thick] (0,0) .. controls +(1,0) and +(-1,0) ..
				(2,1);
			\end{knot}
\end{tikzpicture}}}

\newcommand{\oxunder}{
	\raisebox{-8pt}{\begin{tikzpicture}[scale = .35, pics/arrow/.style={code={%
					\draw[line width=0pt,{Computer Modern Rightarrow[line
						width=0.8pt,width=1.5ex,length=1ex]}-] (-0.5ex,0) -- (0.5ex,0);
			}}]
			
			\begin{knot}[
				clip width=5
				]
				\strand[thick] (0,0) .. controls +(1,0) and +(-1,0) ..
				pic[pos=0,sloped]{arrow}	(2,1);
				\strand [thick] (0,1) .. controls +(1,0) and +(-1,0) ..
				pic[pos=0,sloped]{arrow}	(2,0);
			\end{knot}
\end{tikzpicture}}}
\newcommand{\oxover}{
	\raisebox{-8pt}{\begin{tikzpicture}[scale = .35, pics/arrow/.style={code={%
					\draw[line width=0pt,{Computer Modern Rightarrow[line
						width=0.8pt,width=1.5ex,length=1ex]}-] (-0.5ex,0) -- (0.5ex,0);
			}}]
			\begin{knot}[
				clip width=5
				]
				
				\strand [thick] (0,1) .. controls +(1,0) and +(-1,0) .. pic[pos=0,sloped]{arrow}
				(2,0);
				\strand [thick] (0,0) .. controls +(1,0) and +(-1,0) ..pic[pos=0,sloped]{arrow}
				(2,1);
			\end{knot}
\end{tikzpicture}}}

\newcommand\bigWP{\bm{\mathsf {WP}}}
\newcommand{\bigCA}{\bm{\mathsf{CA}}}
\newcommand{\nubigCA}{\bm{\mathsf{CA}}^-}
\newcommand{\bigOCA}{\bm{\mathsf{OCA}}}
\newcommand\VMO[1][\mathsf X]{{#1}\mathdash\bm{\mathsf{MO}}}
\newcommand\VnuMO[1][\mathsf X]{{#1}\mathdash\bm{\mathsf{MO}}^-}

\newcommand\VCMO[1][\mathsf X]{{#1}\mathdash\mathsf{MO}^{(\CCC, \omega)}}
\newcommand\nuVCMO[1][\mathsf X]{{#1}\mathdash{\mathsf{MO}^{(\CCC, \omega)}}^-}

\newcommand{\inica}{\mathcal F}
\newcommand{\kinica}[1][]{\mathcal U_{#1}}
\newcommand{\Cinica}[1][(\CCC, \omega)]{\mathcal F^{#1}}
\newcommand{\Ckinica}[1][(\CCC, \omega)]{\mathcal U^{#1}}
\newcommand{\uniA}[1][\alg]{\alpha_{#1}}
\newcommand{\vca}{\mathcal V}
\newcommand{\vGca}[1][\theta]{\mathcal V_{#1}}
\newcommand{\dvGca}[1][\theta]{\tilde {\mathcal V}_{#1}}
\newcommand{\uniG}[1][\theta]{\alpha_{#1}}

\newcommand\E{\mathsf{E}}

\newcommand{\Fgraph}{ \xymatrix{
		E \ar@(lu,ld)[]_\tau&& H \ar[ll]_s \ar[rr]^t&& V}}
\newcommand{\Fgraphdash}{\xymatrix{
		E \ar@(lu,ld)[]_{\tau } &&H \ar[ll]_{s } \ar[rr]^{t }& &V }}

\newcommand{\Fgraphvar}[6]{\xymatrix{
		*[r] {#1}\ar@(ul,dl)[]_{#6} && {#2} \ar[ll]_-{#4} \ar[rr]^-{#5}&& {#3}}}

\newcommand\X{\mathcal X}

\newcommand{\TT}{\mathbb T}
\newcommand{\DD}{\mathbb D}

\newcommand{\LL}{\mathbb L}

\title[Perspectives on circuit algebras]{Functorial, operadic and modular operadic combinatorics of circuit algebras}
\author{Sophie Raynor}

\date{\today}

\thanks{The author acknowledges the support of Australian
	Research Council grants DP160101519 and FT160100393.}

\begin{document}

\maketitle

\begin{abstract}
	Circuit algebras are a symmetric analogue of Jones's planar algebras introduced to study finite-type invariants of virtual knotted objects. Circuit algebra structures appear, in different forms, across mathematics. This paper provides a dictionary for translating between their diverse incarnations and describing their wider context. A formal definition of a broad class of circuit algebras is established and three equivalent descriptions of circuit algebras are provided: in terms of operads of wiring diagrams, modular operads and categories of Brauer diagrams. As an application, circuit algebra characterisations of algebras over the orthogonal and symplectic groups are given. 
	
\end{abstract}

\section{Introduction}

	Circuit algebras are a symmetric version of Jones's planar algebras\cite{Jon99}. Their basic data consists of a graded monoid equipped with a contraction (or trace) operation and 
		a levelwise symmetric action. They were introduced by Bar-Natan and Dansco \cite{BND17} as a framework for relating local and global features of virtual tangles in the study of finite-type invariants (see also \cite{DF18, Hal16, Tub14}). Recently, Dancso, Halacheva and Robertson have shown \cite{DHR20} that {oriented} circuit algebras are equivalent to wheeled props \cite{MMS09, Mer10}, and used this to describe the
		graded Kashiwara-Vergne and Grothendieck-Teichm\"uller groups $\mathrm{KRV}$ and $\mathrm{GRT}$ as automorphism groups of circuit algebras \cite{DHR21}.
		
		Though the term ``circuit algebra'' is not commonly used outside quantum topology, circuit algebra structures appear in different guises widely across mathematics. 
 This paper defines a broad class of circuit algebras -- including wheeled props -- and explains how they may be equivalently characterised 
as algebras over an operad, 
as monoidal functors, and as modular operads with an extra operation. 

\begin{thm}
	[\cref{thm. lax functor ca} \& \cref{prop CA MO}] \label{thm. def CA} 
	A circuit algebra is, equivalently
	\begin{enumerate}
		\item an algebra over an operad of wiring diagrams,
		\item a symmetric lax monoidal functor from a category of 
		Brauer diagrams, 
		\item a modular operad equipped with an additional graded product.
	\end{enumerate}
	
\end{thm}

To my knowledge, this is the first time that these perspectives (though not new) have been explicitly stated and compared, together in one work and in such generality. Each description relates to structures that arise in different areas of mathematics, so \cref{thm. def CA} provides a dictionary for translating results between these contexts. Moreover, the categorical and operadic structures underlying each version may be generalised (and specialised) in distinct ways, thereby precisely locating circuit algebras within a diverse zoo of related concepts (see \cref{IntroTab} for a partial overview). 

As an application of this combined approach, and building on \cite{DM23}, the following theorem, providing a circuit algebra characterisation of algebras over the orthogonal and symplectic groups, is proved: 

\begin{thm}[\cref{thm. CA inv}]\label{thm. intro CA inv}
	The category of algebras over the $d$-dimensional orthogonal (respectively symplectic) group is equivalent to a subcategory of circuit algebras that satisfy two simple relations. 
\end{thm}

\medspace

In their original low-dimensional topology and quantum algebra context (first in \cite{BND17}, then e.g.,~\cite{DF18, Hal16, Tub14}), circuit algebras are defined as algebras over operads of wiring diagrams (see e.g.,~\cite{DHR20, DHR21}).

Different flavours of circuit algebras -- including nonoriented \cite{MW21}, oriented (wheeled props), and mixed -- are described by different ``coloured'' operads of wiring diagrams (see \cref{ssec. colour BD} and \cref{defn. operad of wiring diagrams}). This also gives an alternative proof that oriented circuit algebras are wheeled props (see \cref{ex. wheeled props}). For any given colouring, several important generalisations of circuit algebras arise as algebras over suboperads of wiring diagrams.
 The columns of Table \ref{IntroTab} are indexed by structures -- including (coloured) planar algebras \cite{Jon99} and modular operads \cite{GK98, HRY19a, Ray20, RayCA2} -- obtained this way. Rows~(2)-(3) describe the operads and categories governing these in the sense of \cref{thm. def CA}. 
 
 	\begin{table}
 	[htb!]
 	\centerfloat{
 		\begin{tabular}{|P{.01\linewidth} |P{.18\linewidth} || P{.24\linewidth} | P{.22\linewidth} | P{.21\linewidth} | P{.15\linewidth}|}
 			\hline
 			
 	\rowcolor{lightgray} 1&\textbf{Structure} &{\sc{Circuit algebras}} (CAs)& \sc{Nonunital (or downward) CAs} &\sc{Modular operads} & \sc{Planar algebras }\\ 
 		&	Special cases: &&&&\\
 		&	Oriented &Wheeled props, \cite{DHR20}  \newline(see e.g.,~\cite{MMS09, DM23})  & Nonunital wheeled props
 		 (see e.g.,~\cite{Sto23})&Wheeled properads\newline (see e.g.,~\cite{HRY15})&\\
 			\hline 
 			\hline
 			
 		2&	\textbf{Governing operad }& Wiring diagrams (WDs)& downward WDs (Koszul \cite{KW24, Sto23}) & connected WDs& planar diagrams
 			\\
 			\hline 
 		3&	\textbf{Classifying category} &Brauer diagrams (BDs)& downward BDs \small{(cospans, Remark \ref{rmk. pushout})}&  \cellcolor{black}& Temperley-Lieb diagrams\\
 			\hline 
 		4&	\textbf{Rep. theory}&(Sections \ref{ssec. representations of BD}~\&~\ref{ssec unital inv})& (Sections \ref{ssec. representations of BD}~\&~\ref{ssec nonunital CA})&(See CA column)&\\
 		&	\small{nonoriented mono.} & $O_d, \ Sp_d$& $O_{\infty}, Sp_{\infty}$ \cite{SS15}&&\small{(quantum $SU_2$)}\\
 		&	\small{oriented mono.}& $GL_d$ \cite{DM23} &$GL_{\infty}$ \cite{SS15}&&\\
 			\hline 
 	5&	Monad in \cite{RayCA2}&$\LL \DD\TT$& $\LL\TT$ (has arities \cite{RayCA2})& $\DD\TT$&\cellcolor{black} \\ \hline
 			
 	\end{tabular}}
 	
 \vskip 1 ex
 	\caption{Comparison of circuit algebras and generalisations obtained from suboperads of wiring diagrams: The subcategories of Brauer diagrams in Row 3, like the operads in Row 2, are dependent on colouring. Modular operads cannot be described as functors from categories of Brauer diagrams. Row 4 indicates (dimension parameter dependent) groups that are related by Schur-Weyl duality to (sub)categories of \textit{monochrome} (non)oriented Brauer diagrams. The planar case is not studied in this work. Row 5 refers to the (colouring-independent) monads described in \cite{RayCA2}.}	\label{IntroTab}
 \end{table}

Statements (1)~\&~(2) of \cref{thm. def CA} are already implicit in the original definition of circuit algebras \cite{BND17}. 
	Their equivalence is a formal consequence of the definition, in \cref{ssec. wd and ca}, of operads of wiring diagrams in terms of a ``category of Brauer diagrams'' (or ``Brauer category'', c.f.~\cref{rmk Brauer diagrams}). 
 Such diagrams have been widely studied since Brauer's 1937 paper \cite{Bra37} extending Schur-Weyl duality to representations of the finite dimensional orthogonal and symplectic groups (see e.g.,~{\cite{Wen88, Koi89}). More recently, categories of Brauer diagrams have been used to simultaneously study 
 	systems of related representations \cite{LZ15, RS20, SS15, SS20I}. So, Statement~(2) of \cref{thm. def CA} implies a link between circuit algebras and classical themes in representation theory. However, the proof of \cref{thm. intro CA inv} in \cref{sec. invariants} does not explicitly use these methods. 
 	Instead, since 
 wheeled props are equivalent to oriented circuit algebras, \cref{thm. intro CA inv} is proved by adapting Derksen and Makam's invariant-theoretic approach to wheeled props \cite{DM23} 
 	from the oriented, to the unoriented case (see \cref{sec. invariants}). 
 
   The final characterisation (3) in \cref{thm. def CA} describes circuit algebras as modular operads equipped with an extra product operation. Modular operads were first introduced in \cite{GK98} to study moduli spaces of higher genus curves. General unital modular operads, as in \cite{HRY19a, HRY19b, Ray20}, may be obtained from \cref{thm. def CA}, (1) by restricting to a suboperad of \emph{connected} wiring diagrams. Unlike the restriction to planar diagrams, which respects the categorical structure (in the sense of \cref{prop. lax functors and algebras}), this is a purely operadic construction and admits no categorical description in terms of Brauer diagrams. 
 
This paper is one a pair that, together, provide a detailed conceptual and technical account of circuit algebra combinatorics. In the companion paper \cite{RayCA2}, I use the modular operadic perspective to build on the results of \cite{Ray20} and construct a monad and graphical calculus and prove an abstract nerve theorem for circuit algebras. 
 Thus, circuit algebras also admit combinatorial characterisations as algebras for a monad on a category of graded symmetric objects, and as ``Segal presheaves'' on  a category of graphs \cite[Section~8]{RayCA2}. 
 
The monad for circuit algebras in \cite{RayCA2} is constructed, using \emph{iterated distributive laws} \cite{Che11}, as a composite $\LL \DD\TT$ of three simpler monads, each governing a different aspect of the circuit algebra structure. This piecewise construction is central to the proof of the nerve theorem \cite[Theorem~8.4] {RayCA2}. It also dovetails with the other perspectives in \cref{thm. def CA}.

For example, algebras for the monad $\LL\TT$ are nonunital circuit algebras  (\cref{IntroTab}), that do not have units for the modular operadic multiplication.  
Their combinatorics (see e.g.,~\cite{KW24, Sto23}) are simpler than the unital case since they avoid the  ``\emph{problem of loops}''  \cite[Section~6]{Ray20}. In the language of Brauer categories (in the sense of \cite{LZ15, SS15}, see \cref{ssec. representations of BD}), this problem of loops refers simply to the dimension parameter associated to the unit trace. Under \cref{thm. def CA}, nonunital circuit algebras correspond precisely with symmetric monoidal functors from subcategories of ``\emph{downward}'' Brauer diagrams, that cannot encode (finite) dimension. 
Sam and Snowden \cite{SS15} have established equivalences between functors from the subcategory of downward monochrome {oriented }Brauer diagrams and representations of the infinite dimensional ({stable}) general linear group $GL_{\infty}$, and between functors from the subcategory of downward (monochrome nonoriented) Brauer diagrams and representations of the infinite dimensional orthogonal and symplectic groups $ O_{\infty}$ and $Sp_{\infty}$ (see \cref{rmk stable rep} \& \cref{ssec nonunital CA}).

Some particularly nice properties of the combinatorics of nonunital circuit algebras are included in Rows 2,3~\&~5 of \cref{IntroTab}.  The modular operadic perspective on nonunital circuit algebras, together with the results of \cite{SS15}, has been exploited in \cite{KRW21} to prove that the Malcev Lie algebras associated to the Torelli groups of surfaces of arbitrary genus are stably Koszul. The relationship is also noted in \cite{Sto23} where nonunital (d.g.)-modular operads are characterised as lax functors from a \textit{``Brauer properad''} obtained by restricting to connected diagrams in the initial circuit algebra.

\medspace

The primary aim of this paper is to provide a precise formal framework for studying a broad class of circuit algebra structures as they arise across mathematics, and thereby extend the toolboxes of 
representation theorists, low-dimensional topologists and operad theorists alike. This presents a plethora of options for generalising circuit algebras and for translating results in new contexts:

A particular motivation for a formal study of circuit algebra structures (here and in \cite{RayCA2}) comes from the work of  Dansco, Halacheva and Robertson 
\cite{DHR21} who have used circuit algebras to obtain results relating the \textit{graded} {Grothendieck-Teichm\"uller} and {Kashiwara-Verne} 
groups $\mathrm{GRT}$ and $\mathrm{KRV}$. In order to extend these results to the ungraded groups $\mathrm{GT}$ and $\mathrm{KV}$, it is necessary to relax the circuit algebra axioms up to homotopy \cite[Introduction,~ Remark~1.1]{DHR21}. Weakening the characterisation in \cite[Theorem~8.4]{RayCA2} of circuit algebras in terms of Segal functors suggests one way to do this. However, there are difficulties adapting the methods, used in \cite{HRY19b} and 
\cite{Ray20} to construct Segal models for homotopy modular operads, to model homotopy circuit algebras \cite[Section~8.4]{RayCA2}. 

Stoeckl's construction \cite{Sto23} of a model for nonunital $(\infty, 1)$-wheeled props, and the proof, in \cite{KW24}, that the operad for monochrome nonunital circuit algebras is Koszul, potentially provide another (operadic) approach to constructing a model. From the categorical perspective, Sharma's model structure for compact closed categories \cite{Sha21} may also shed light on this question. 

Several questions about duality arise from the circuit algebra characterisations in \cref{thm. def CA}.
For example: Can the operadic perspective provide new insights into the Schur-Weyl duality of the classical groups and their quantisations? Given that the operads governing nonunital wheeled props and circuit algebras are Koszul \cite{KW24, Sto23}, it is natural and useful to ask whether this is also true of the operads for unital circuit algebras (i.e., operads of wiring diagrams). How can this be interpreted in terms of the (downward) Brauer diagram categories? Is there a general Tannakian formalism \cite{JS91} for such questions? (I thank Ross Street for helpful discussions on duality.)

Finally, the categorical and graphical structures governing circuit algebras are seeing increasing applications outside pure mathematics. 
They provide a powerful formal framework for organising, understanding and classifying complex networked systems, by studying their local-global-local structure. Potentially, these methods could help define the theoretical limits of emerging technologies, as well as improving transparency (e.g.,~in AI) and informing efficient design of algorithms and software. For example, the ZX-calculus \cite{CK17}, that provides a rigorous graphical formalism for quantum computation (and could, potentially, make quantum computation accessible to a general audience \cite{QiP23}), admits a circuit algebra description. It would be interesting to compare this with circuit algebras that arise in quantisation problems \cite{BND17, DHR21}.

\subsection{Overview}

Categorical preliminaries are given in \cref{sec. categories} to establish notation and terminology for the (symmetric monoidal category) concepts rest of the paper. 

 \cref{sec. BD} provides a detailed discussion of the categories of (coloured) Brauer diagrams, and describes their relation to several known results on the invariant theory of classical groups (c.f.~ \cite{LZ15, SS15}). 
 
Categories of Brauer diagrams are used, in \cref{sec. CA}, to define circuit algebras. 
\cref{ssec operads} provides a quick introduction to operads and their algebras. In \cref{ssec. wd and ca}, operads of wiring diagrams and circuit algebras are introduced and defined using categories of Brauer diagrams from \cref{sec. BD}.

In \cref{sec: definitions}, an axiomatic characterisation of circuit algebras is given and it is shown that they are modular operads that admit an extra graded product. Finally, in \cref{sec. invariants}, \cref{thm. intro CA inv} is proved as an application of the preceding ideas. The method is then extended to give a nonunital circuit algebra characterisation (\cref{thm. Ginfty}) of $O_{\infty}$ and $Sp_{\infty}$.

The companion paper \cite{RayCA2} builds on the modular operadic perspective to obtain a graphical calculus, monad and nerve theorem for circuit algebras. The machinery used, involving a combined application of iterated distributive laws \cite{Che11} and abstract nerve theory \cite{BMW12}, is also explained in detail \cite[Section~2]{RayCA2}.

\subsection*{Acknowledgements}

I thank Marcy Robertson, Zsuzsanna Dancso, and Chandan Singh and Kurt Stoeckl for encouraging my interest and learning in this field. I am grateful to Ole Warnaar for all his support, to Kevin Coulembier for patiently explaining some representation theory and to my students and colleagues at James Cook University, Bindal Country, 
for their curiosity and friendship. I thank the members of Centre of Australian Category Theory, Macquarie University, Dharug Country, where I first began thinking about this work. I am particularly grateful to Ross Street for his friendship and patience discussing duality with me, and Richard Garner for a remark that led to new perspectives. 

\section{Key categorical concepts}\label{sec. categories}

This section provides a brief outline  of the notation and terminology conventions for symmetric monoidal categories 
that will be used in the rest of the paper. For precise definitions and a detailed discussion of symmetric monoidal categories, see e.g.,~\cite[Chapters~2~\&~8]{EGN15}.

\subsection{Symmetric monoidal categories}\label{ssec. SMC}

A \textit{monoidal category} is a category $\V$ together with a bifunctor $\otimes \colon \V \times \V \to \V$ (the \textit{monoidal product}) that is associative up to natural \textit{associator} isomorphism, and for which there is an object $I$ of $\V$ (the \textit{monoidal unit}) that acts as a two-sided identity for $\otimes$ up to natural 
\textit{unitor} isomorphisms. The monoidal product and the associator and unitor isomorphisms are required to satisfy axioms that mean that certain sensible diagrams commute. If the associator and unitor isomorphisms are the identity, then the monoidal category 
is called \textit{strict monoidal}.

A \textit{braiding} on a monoidal category $(\V, \otimes, I )$ is a collection of  isomorphisms $\sigma_{x, y} \colon x \otimes y \rightarrow y \otimes x $ (defined for all $x, y \in \V$)
that satisfy the braid 
identities 
\begin{equation} \label{eq. braids} (\sigma_{y, z}\otimes id_x)(id_y \otimes\sigma_{x, z})(\sigma_{x, y} \otimes id_z) = (id_z\otimes \sigma_{x, y})(\sigma_{x, z}\otimes id_y)( id_x\otimes\sigma_{y, z})  \ \text{ for all } x, y, z. \end{equation} If $\sigma_{y,x} =  \sigma_{x,y}^{-1}$ for all $x, y$, then the monoidal structure on $\V$ is \textit{symmetric}. 

\begin{rmk}
	
	In this paper, associators, unitors and symmetry (braiding) isomorphisms will be ignored in the notation, and (symmetric) monoidal categories will be denoted simply by $\V$ or $(\V, \otimes, I)$.

\end{rmk}

\begin{ex}\label{ex. slice}
	
	For any category $\V$ and object $ x \in \V$, objects of the \emph{slice category $ x \ov \V$ of $\V$ under $x$}  are pairs $(y, f)$ where $f \in \V(x,y)$. Morphisms $(y,f) \to (y', f')$ are commuting triangles in $\V$ of the form:
	\[ \xymatrix{&x \ar[dl]_f \ar[dr]^{f'} &\\ y \ar[rr]^g && y'.}\] 
	The \emph{slice category $\V \ov x$ of $\V$ over $x$} is defined similarly, with objects $(y,f) \colon f \in \V(y,x)$ and morphisms $g \colon (y,f) \to (y',f')$ given by morphisms $ g \in \V(y', y)$ such that $ f \circ g = f'$.

	If $(\V, \otimes, I)$ is a monoidal category, then in general $\V \ov x$ (respectively $x \ov \V$) does not inherit a monoidal structure from $\V$. However, since $I \otimes I \cong I$ by definition, $\otimes$ defines a monoidal product on $I \ov \V$ (respectively $\V \ov I$) with unit $id_I \in \V(I,I)$. 
\end{ex}

\begin{defn}\label{def perm prop}
	
	Symmetric strict monoidal categories are called \emph{permutative categories}. The notation $\oplus$ and $0$ will often be used to designate the monoidal product and unit of a permutative category. 
	
	An (ordinary) \emph{$\DDD$-coloured prop} is a small permutative category $\PCat$ whose object monoid is free on a set $\DDD$. When $\DDD = \{1\}$ is a singleton, then $\PCat$ is a (monochrome) prop (with object set $\N$) in the original sense of \cite{MacL65}. 
\end{defn}

\begin{ex}\label{ex. Sigma}
	For each $n \in \N$, let $ \nn$ denote the set $\{1, 2, \dots, n\}$ (so $\nul = \emptyset$), and let $ \Sigma_n$ be the group of permutations on $\nn$. Let $\Sigma$ be the \emph{symmetric groupoid} with $\Sigma(n,n) = \Sigma_n$ for all $n$, and $\Sigma (m,n) = \emptyset$ when $m \neq n$. Addition of natural numbers gives $\Sigma$ a (monochrome) prop structure.

	More generally, let $\DDD$ be a set, and let $\listm (\DDD) =\coprod_{n \in \N} \DDD^n$ denote the set of finite ordered sets $\ccc  = (c_1, \dots, c_n)$ of elements of $\DDD$. So $\listm (\DDD)$ underlies the free associative monoid on $\DDD$. For $\ccc = (c_1, \dots, c_m)$ and $\ddd = (d_1, \dots, d_n)$ in $\listm (\DDD)$, their (concatenation) product $\ccc \ddd = \ccc\oplus \ddd$ is given by \[  \ccc \ddd  \defeq (c_{1}, \dots, c_{m}, d_1, \dots, d_n) .\] The empty list is the unit for $\oplus$ and is denoted by $\varnothing$ (or $\varnothing_{\DDD}$).
	
	The symmetric groupoid $\Sigma$ acts on $\listm (\DDD)$ from the right by $\sigma \colon  (\ccc \sigma) \defeq (c_{\sigma 1}, \dots, c_{\sigma m}) \mapsto \ccc $, for all $\ccc= (c_1, \dots, c_m)$ and $\sigma \in \Sigma_m$. 
	
	The $\DDD$-coloured prop so obtained is the \emph{free symmetric groupoid $\Sigma^{\DDD}$ on $\DDD$}.
\end{ex}
\begin{ex}
	\label{ex. Symmetric objects} For any category $\V$, a functor $S \colon \Sigma \to \V$ is equivalently described by a sequence $(S(n))_n$ of objects of $\V$ such that $\Sigma_n$ acts on $ S(n)$ for all $n$. 

 A \emph{$\listm (\DDD)$-graded symmetric object in $\V$} is a functor $ B \colon \Sigma ^{\DDD} \to \V$. Equivalently, it is a collection $(B (\ddd))_{ \ddd \in \listm (\DDD)}$ of $ \V$-objects, and $\V$-isomorphisms 
		$ B (\sigma \ddd) \xrightarrow{\cong}B (\ddd )$,
		defined for all $ \ddd = (d_1, \dots, d_n) \in \DDD^n$ and all $\sigma \in \Sigma_n$.

\end{ex}

 Let $\mathsf V$ be a symmetric monoidal category. In a $\mathsf V$-(enriched) category, the hom sets are instead $\mathsf V$-objects and composition is a $\mathsf V$-morphism such that compatibility axioms are satisfied. Other than (ordinary) $\Set$-enriched categories, this paper will also consider categories \emph{enriched} in the categories $\vect$ of $\Bbbk$-vector spaces (where $\Bbbk$ is an algebraically closed field of characteristic 0), and $\modR$ of modules over a commutative ring $R$.

\begin{ex}
	\label{ex symmtric algebra}

	 Let $R$ be a commutative ring. Then $ R[\Sigma] \defeq \bigoplus_{n \in \N} R[\Sigma_n]$, where $R[\Sigma_n]$ denotes the group algebra (for $n \in \N$), describes the free $\modR$-prop on $\Sigma$.  
\end{ex}

\begin{ex}
	\label{ex. endo prop}

Given a vector space $V$, the ($\vect$-enriched) \emph{endomorphism prop associated to $V$} is denoted by $T(V)$, with $T(V)(m,n) = Hom_{\Bbbk} (V^{\otimes m}, V^{\otimes n})$, 
the space of linear transformations $V^{\otimes m } \to V^{\otimes n}$. 
By convention, $V^{\otimes 0} = \Bbbk$, so $T(V)(0,0) = \Bbbk$, and $T(V)$ embeds canonically in $\vect$ as the full subcategory with objects $V^{\otimes n}$, $n \in \N$. These categories are identified in what follows.

For each $n \in \N$, the symmetric group $\Sigma_n$ acts on $V^{\otimes n}$ by permuting factors. Hence $\Sigma$ acts on $T(V)$ levelwise.
\end{ex}

\begin{defn}\label{def monoidal functor}

	A \emph{(lax) monoidal functor} $(\Theta,\eta_{\Theta} , \theta)\colon (\V_1, \otimes_1, I_1) \to (\V_2, \otimes_2, I_2)$ consists of a functor 
	$\Theta \colon \V_1 \to \V_2$, together with a morphism $\eta = \eta_{\Theta} \colon I_2 \rightarrow \Theta (I_1)$ in $\V_2$ and a natural transformation $\pi = \pi_{\Theta}\colon \Theta(-) \otimes_{2} \Theta(-) \Rightarrow \Theta(- \otimes_1 -)$ such that all the expected structure diagrams commute. A monoidal functor $ (\Theta, \pi, \eta )$ is called \emph{strong} if $\pi$ and $ \eta$ are invertible, and \emph{strict} if they are the identity.

\end{defn}

\begin{ex}
	\label{ex. universal Sigma}
		As in \cref{ex. Sigma}, let $\Sigma$ be the symmetric groupoid. For any symmetric monoidal category $(\V, \otimes, I)$ and any choice of object $x \in \V$, there is a unique symmetric strict monoidal functor $\Sigma \to \V$ with $0\mapsto I$ and $1 \mapsto x $.

\end{ex}

\begin{defn}\label{defn symmetric monoid}

	 A \emph{$\listm (\DDD)$-graded symmetric monoid in $\V$} is a symmetric monoidal functor \[ (B, \pi, \eta) \colon (\Sigma^{\DDD}, \oplus, \varnothing) \to (\V, \otimes, I)\] where  $(\Sigma^{\DDD}, \oplus, \varnothing) $ is the prop defined in \cref{ex. Symmetric objects}. The structure maps $(\pi, \eta)$ describe a commutative and associative (up to symmetry and associators in $\V$) unital monoid structure on the underlying graded symmetric object $(B(\ddd))_{\ddd}$.

\end{defn}

\begin{rmk}
	\label{rmk. tensor categories} 
		Enriched (lax) monoidal $\mathsf V$-functors between monoidal $\mathsf V$-categories are defined as in \cref{def monoidal functor} except that the underlying functor is $\mathsf V$-enriched and the structure maps are $\mathsf V$-morphisms such that the relevant diagrams commute in $\mathsf V$.
		
	Symmetric monoidal categories enriched in a linear category (such as $\vect$) are often called \emph{tensor categories} \cite{EGN15}. In the tensor category literature, \emph{(tensor) functors} between tensor categories are usually assumed to preserve the monoidal product strictly. This contrasts with the approach of this paper where all monoidal functors are assumed to be lax, unless explicitly stated otherwise. 
	
\end{rmk}

\subsection{Categorical duality and trace}\label{ssec: duality}

\begin{defn}
	\label{defn. compact closed} 
	An object $x$ of a symmetric monoidal category $\V$ has a \emph{dual object} $x^* $ in $\V$ if there are morphisms  $\cup_x\colon I\to x \otimes x^*$ and $\cap_x  \colon x^* \otimes x \to I$ that satisfy the \textit{triangle identities} (illustrated in \cref{fig. triangle 1}): 
	\begin{equation} \label{eq. dual}
		(\cap_x \otimes id_x) \circ (id_x \otimes \cup_x)= id_x = (id_x \otimes \cap_{x^*} )\circ (\cup_{x^*} \otimes id_x).
	\end{equation}.

	A \emph{compact closed category} is a symmetric monoidal category such that every object has a dual \cite{MacL65}.
	
\end{defn}

\begin{figure}	[htb!]
	
	\begin{tikzpicture}[scale = .16]
		\begin{pgfonlayer}{background}
			
			\node  (2) at (3.75, 5) {};
			
			\node  (4) at (1.25, 5) {};
			
			\node  (7) at (6.25, 5) {};
			
			\node  (10) at (1.25, 10) {};
			
			\node  (13) at (6.25, 0) {};
			
			
			\node  (15) at (25, 5.25) {};
			\node  (16) at (21, 5.25) {};
			\node  (17) at (21.75, 5.25) {};
			\node  (18) at (22.5, 5.25) {};
			\node  (19) at (24.25, 5.25) {};
			\node  (20) at (23.5, 5.25) {};
			\node  (21) at (18.5, 5.25) {};
			\node  (22) at (19.25, 5.25) {};
			\node  (23) at (20, 5.25) {};
			\node  (24) at (25, 10.25) {};
			\node  (25) at (24.25, 10.25) {};
			\node  (26) at (23.5, 10.25) {};
			\node  (27) at (18.5, 0.25) {};
			\node  (28) at (19.25, 0.25) {};
			\node  (29) at (20, 0.25) {};

			\node  (30) at (12.75, 10) {};
			\node  (32) at (12, 10) {};
			\node  (33) at (12, 0) {};
			\node  (35) at (12.75, 0) {};
			
		\end{pgfonlayer}
		\begin{pgfonlayer}{above}
			\draw (10.center) to (4.center);
			\draw (7.center) to (13.center);
			\draw [bend right=90, looseness=2.25] (4.center) to (2.center);
			\draw [bend left=90, looseness=2.50] (2.center) to (7.center);
			\draw (25.center) to (19.center);
			\draw (22.center) to (28.center);
			\draw [bend left=90, looseness=2.00] (19.center) to (17.center);
			\draw [bend right=90, looseness=2.25] (17.center) to (22.center);
			\draw (32.center) to (33.center);

			\draw 
			(8, 4.8)--(9.5,4.8)
			(8, 5.3)--(9.5,5.3);
			
			\draw 
			(14.5, 4.8)--(16,4.8)
			(14.5, 5.3)--(16,5.3);

		\end{pgfonlayer}
	\end{tikzpicture}

	\caption{ String diagram representation of the triangle identities. } \label{fig. triangle 1} \end{figure}
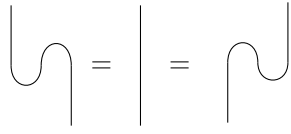

Let $\V$ be a compact closed category. For all morphisms $f \in \V(x,y)$, there is a corresponding \emph{evaluation morphism} $\ev{f} \in \V(y^* \otimes x, I)$ induced by composition with $\cap_y$ (\cref{fig eval} (a)) and  \emph{coevaluation morphism} $\coev{f} \in \V(I, y \otimes x^*)$ induced by composition with $\cup_x$ (\cref{fig eval} (b)):
\begin{equation}\label{eq ev coev dual}
	\ev{f} \defeq \cap_y\circ (id_{y^*} \otimes f)  \text{ and }  \coev{f} \defeq  (f\otimes  id_{x^*})\circ \cup_x,
\end{equation}
and a \emph{dual morphism} (called the \textit{transpose morphism} in e.g.,~\cite{Sha21}) $f^* \in \V(y^*, x^*)$ (\cref{fig eval} (c)):
\begin{equation}\label{eq. BD dual}
	f^* \defeq (\cap_y \otimes id_{x^*})\circ ( id_{y^*} \otimes f \otimes  id_{x^*} ) \circ ( id_{y^*} \otimes  \cup_x). 
\end{equation}

\begin{figure}	[htb!]
	
	\begin{tikzpicture}
		\node at (-1,1) {(a)};		
		\node at (0,0){\begin{tikzpicture}[scale = .2]
				\begin{pgfonlayer}{background}
					
					\node  (0) at (-2, 0) {};
					\node  (1) at (2, 0) {};
					\node  (2) at (-2, -5) {};
					\node  (3) at (2, -5) {};
					
					\draw (0.center) to (2.center);
					\draw (1.center) to (3.center);
					\draw [bend left=90, looseness=5.00] (0.center) to (1.center);
				\end{pgfonlayer}
				
				\node at (-2,-6){\scriptsize{$y^*$}};
					\node at (2,-6){\scriptsize{$x$}};

				\filldraw[white] (1,-1) rectangle ++(2cm,2cm);
				\draw (3, 1)--(3, -1)--(1, -1)--(1, 1)-- (3,1);
				\node at (2.1,-0.1){\scriptsize{$f$}};
		\end{tikzpicture}};
		\node at (2.5,1) {(b)};

		\node at (3.5,0){\begin{tikzpicture}[scale = .2]
				\begin{pgfonlayer}{background}
					\node  (0) at (-2, 0) {};
					\node  (1) at (2, 0) {};
					\node  (2) at (-2, 5) {};
					\node  (3) at (2, 5) {};
					\draw (2.center) to (0.center);
					\draw (3.center) to (1.center);
					\draw [bend right=90, looseness=5.00] (0.center) to (1.center);
				\end{pgfonlayer}
				
				\node at (-2,6){\scriptsize{$y$}};
				\node at (2,6){\scriptsize{$x^*$}};
				
				\filldraw[white] (-3,-1) rectangle ++(2cm,2cm);
				\draw (-3, 1)--(-3, -1)--(-1, -1)--(-1, 1)-- (-3,1);
				\node at (-2,-0.1){\scriptsize{$f$}};

		\end{tikzpicture}};
		\node at (6, 1) {(c)};		
		\node at (8,0){\begin{tikzpicture}[scale = .2]
				\begin{pgfonlayer}{background}
					
					\node (0) at (-2, 0) {};
					\node (1) at (2, 0) {};
						\node(2) at (-6, -5){};
							\node(3) at (2,5){};
					\node (4) at (-6, 0) {};
					\node (5) at (-6, 5) {};
					\node (6) at (2, -5) {};

					\draw [bend left=90, looseness=2.70] (4.center) to (0.center);
					\draw (2.center) to (4.center);
					\draw [bend right=90, looseness=2.70] (0.center) to (1.center);
					\draw (3.center) to (1.center);
				\end{pgfonlayer}

					\node at (-6,-6){\scriptsize{$y^*$}};
				\node at (2,6){\scriptsize{$x^*$}};
				
				\filldraw[white] (-3,-1) rectangle ++(2.5cm,2cm);
				\draw (-3, 1)--(-3, -1)--(-1, -1)--(-1, 1)-- (-3,1);
				\node at (-2,-0.1){\scriptsize{$f$}};

		\end{tikzpicture}};
	\end{tikzpicture}
	\caption{(a) $\ev{f}\colon y^* \otimes x \to I$; (b)$\coev{f}\colon I \to y \otimes  x^*$; (c) $f^*\colon y^*\to x^*$ .} \label{fig eval}
\end{figure}

In particular $\cap_x = \ev{id_x}$, $\cup_x = \coev{id_x}$ and $(id_x)^* = id_{x^*}$ for all objects $x$. And, for composable morphisms $f$ and $g$, $(g \circ f)^* = f^* \circ g^*$ in $\V$.

\begin{ex}
	\label{ex. compact closed}
Let $\Bbbk$ be a field. The monoidal category $(\vectf, \otimes , \Bbbk)$ of finite dimensional $\Bbbk$-vector spaces has a canonical compact closed structure given by $V^* = \vect (V, \Bbbk)$. For each $ V \in \vectf$, its dimension $dim(V)$ over $\Bbbk$ is equal to its \emph{categorical dimension} given by $\cap_V \circ \cup_V \in \Bbbk$.

\end{ex}

A \emph{traced symmetric monoidal category} \cite{JSV96} is a monoidal category $(\V, \otimes, I)$ equipped with a family of \emph{(partial) trace} functions $tr_{x,y}^z \colon \V(x \otimes z, y \times z) \to \V(x,y)$, natural in  objects $x,y,z \in \V$ and satisfying:
\begin{description}
	\item [Vanishing] For all objects $x,y,a,b \in \V$,  $tr_{x,y}^I$ is the identity on $\V(x,y) = \V(x \otimes I, y \otimes I)$, and $ tr_{x,y}^{a\otimes b} = tr_{x,y}^a \circ tr^b_{x \otimes a, y \otimes a} \colon \V(x \otimes a \otimes b, y \otimes a \otimes b) \to \V (x,y)$
	\item [Superposing] For all $ f \in \V(x \otimes a, y \otimes a)$ and $g \in \V(w, z)$, $tr_{w \otimes x, z \otimes y}^a(g  \otimes f) = g \otimes tr_{ x, y}^a( f).$
	\item [Yanking] Let $\sigma_{x,y} \colon x \otimes y \to y \otimes x$ denote the symmetry in $ \V$.  For all $x \in \V$, $tr_{x,x}^x (\sigma_{x,x}) = id_x$.
\end{description}

(In a $\mathsf V$-enriched traced monoidal category, the trace $tr$ is described by $\mathsf V$-morphisms.)

A compact closed category $(\V, \otimes , I, *)$ is traced monoidal with trace defined by 
\[ tr_{x,y}^a (f) \defeq (id_y \otimes \cap_{a^*}) \circ f \circ (id_x \otimes \cup_a) \in \V(x,y)\]  for all $a,x,y \in \V$ and $f\in \V(x\otimes a, y \otimes a)$. As \cref{ex wheeled endo prop} shows, the converse is not true. However, via the \textit{``Int construction''} \cite[Section~4]{JSV96}, any traced symmetric monoidal category $ \V$ embeds fully  faithfully in its compact closed completion $Int(V)$.

Another special class of traced symmetric monoidal category is given by wheeled props. These appear in a variety of contexts involving algebraic structures with trace operations (see e.g.,~\cite{Mer10, MMS09}).

\begin{defn}\label{defn WP}
	A \emph{($\DDD$-coloured) wheeled prop} $(P, \oplus
	_P, 0, tr_P)$ is a ($\DDD$-coloured) prop  $(P, \otimes _P, I)$ equipped with a trace $tr_P$ satisfying the axioms of \cite{JSV96}. 
	
\end{defn}
Most applications consider wheeled props enriched in a linear category such as $\vect$. (\cref{prop: wheeled prop oriented brauer} describes wheeled props enriched in $\mathsf V$ in terms of symmetric monoidal functors to $\mathsf V$.)

\begin{rmk}\label{rmk WP monad}
	Wheeled props are usually defined as algebras for a \textit{graph substitution} monad (see e.g., \cite{Mer10,MMS09,DHR20}). The equivalence of \cref{defn WP} with the graph substitution definition follows from \cite[Theorem~7.9]{RayCA2}.
\end{rmk}

\begin{ex}\label{ex wheeled endo prop}

Let $V$ be a finite dimensional $\Bbbk$-vector space. The endomorphism prop $T(V)$ described in \cref{ex. endo prop} is not compact closed since the dual space $V^* = \vect(V, \Bbbk)$ is not an object of $T(V)$. However, the canonical isomorphism $T(V)(m,n) = \vect (V^{\otimes m}, V^{\otimes n}) \cong (V^*)^{\otimes m} \otimes  V^{\otimes n}$ (for all $m, n \in \N$) induces a trace on $T(V)$ by
	\[ v_1 \otimes \dots \otimes v_n \otimes \alpha_1 \otimes \dots \otimes \alpha_m  \mapsto  \alpha_m(v_n)( v_1 \otimes \dots \otimes v_{n-1}\otimes \alpha_1 \otimes \dots \otimes \alpha_{m-1} ). \] 
		Henceforth, 
		$T(V)$ will be assumed to be a wheeled prop with the canonical trace.

The ($\{V, V^*\}$-coloured) mixed tensor prop $ T^{\dipal}(V) \subset   \vectf$ is closed under duals and thus inherits the compact closed structure from $\vectf$. It is straightforward to check that $T^{\dipal}(V)$ is equivalent -- via shuffle permutations of mixed tensor products $(V^*)^{\otimes m} \otimes V^{\otimes n}$ -- to the compact closed category $Int (T(V))$ obtained by applying the Int construction of \cite{JSV96}.  

\end{ex}

\section{Brauer diagrams}\label{sec. BD}

Circuit algebras are defined in \cref{sec. CA} as algebras over an operad of wiring diagrams. It will follow from \cref{thm. lax functor ca} that they admit an equivalent description as symmetric monoidal functors from categories of (coloured) Brauer diagrams. These diagrams are an important tool in the representation theory of orthogonal, symplectic and general linear groups \cite{Bra37, SS20I}.

The category $\BD$ of monochrome Brauer diagrams is described in \cref{ssec. mono BD}. In \cref{ssec. colour BD}, this definition is generalised to categories of coloured Brauer diagrams, of which oriented Brauer diagrams -- that encode the combinatorics of wheeled props (c.f.,~\cref{prop: wheeled prop oriented brauer}) -- are a special case.

\begin{rmk}
	\label{rmk Brauer diagrams} Several variations of the categories of Brauer diagrams defined in this work have appeared in diverse contexts, usually under the name \emph{``Brauer category''}: For some authors (e.g.,~\cite{Ban16,MW21}), Brauer categories are ordinary categories and coincide with the categories $\BD$ (and $\CBD$) described in this section. However, most works (e.g.,~\cite{LZ15, SS15, RS20,SS20I}) define linear Brauer categories, enriched in the category $\modR$ of $R$-modules for some commutative ring $R$. The definition of these categories is dependent on a choice of parametrising element of the ground ring. 
	
Hence, to distinguish them from linear versions, the categories $\BD$, $\CBD$ described here are called ``categories of Brauer diagrams''.
\end{rmk}

\subsection{Monochrome Brauer diagrams}\label{ssec. mono BD}

The category $\BD$ of (nonoriented monochrome) Brauer diagrams may be pithily defined 
as the free compact closed category generated by a single self-dual object. This section gives a more concrete description of $\BD$, in terms of pairings on finite sets.

\begin{defn}
	
	A \emph{pairing} (\emph{perfect matching}) on a set $X$ is a fixed point free involution $\tau $ on $ X$. 
\end{defn} 
Equivalently, a pairing $\tau$ on $X$ is a partition of $X$ into two-element subsets. In particular, a finite set $X$ admits a pairing if and only if it has even cardinality. The empty set has trivial pairing $\varnothing$ by convention. 

\begin{ex}\label{ex. manifold pairing} If $\man$ is a compact 1-manifold, then its boundary $\partial \man$ has a canonical pairing $\tau^\man$ such that $x  = \tau^\man y$ if $x$ and $y$ are in the same connected component of $\man$ and $x \neq y$.
\end{ex}

\begin{defn}
	\label{def Brauer diag}

	A \emph{(monochrome) Brauer diagram} $f$ between natural numbers $m$ and $n$ is a pair $(\tau_f, \kclf)$ of 
	a pairing $\tau$ on the disjoint union  $ \Sf \amalg \Tf$ -- where $\Sf = \{ s_1 , \dots, s_m\}$ is the \textit{source}, and $\Tf = \{t_1, \dots, t_n\}$ is the \textit{target}, of $f$ -- and a natural number $\kclf$ called the \emph{number of closed components of $f$}. An \emph{open Brauer diagram} is a Brauer diagram  $\tau = (\tau, 0)$ with no closed components. 
	
\end{defn}

Let $\BD(m,n)$ denote the set of Brauer diagrams from $m$ to $n$.

\begin{ex}
	\label{ex symmetric groupoid}
	For all $n$, there is a canonical inclusion $\Sigma_n \hookrightarrow \BD(n,n)$ that takes $\sigma \in \Sigma_n$ to the open Brauer diagram induced by the pairing $s_i \mapsto t_{\sigma i}$ on $\{s_1, \dots, s_n\} \amalg \{t_1, \dots, t_n\}$, $1 \leq i \leq n$.
	
	In particular, the pairing $s_i \mapsto t_i$, $1 \leq i \leq n$ defines the \emph{identity (open) Brauer diagram} $id_n$ on $n$.
	
\end{ex}

Brauer diagrams may be represented graphically as follows: a pairing $\tau$ on the disjoint union $X \amalg Y$ of finite sets $X$ and $Y$ is described by a univalent graph whose vertices are indexed by $X \amalg Y$, with elements of $X$ below those of $Y$, and edges connecting vertices $v_1 $ and $v_2$ if and only if the corresponding elements of $X \amalg Y$ are identified by $\tau$. A Brauer diagram $f = (\tau, \kcl) \colon m \to n$ may be represented by the graph for $\tau$, together with 
$\kcl$ closed circles (called \textit{bubbles} in \cite{RS20}) drawn next to it. 

Given finite sets $X,Y, Z$, and pairings $\tau_{X,Y}$ and $\tau_{Y,Z}$ on $X \amalg Y$ and $Y \amalg Z$ respectively, one may vertically stack the diagrams for $\tau_{X,Y}$ and $\tau_{Y,Z}$ as in  \cref{fig. pairing comp} to obtain a pairing on $X \amalg Z$: 

Namely, 
$\tau_{X,Y}$ and $\tau_{Y,Z}$ generate an equivalence relation on $X \amalg Y \amalg Z$ 
where objects $x$ and $y$ are equivalent if and only if they are related by a sequence of (alternating) applications of $\tau_{X,Y}$ and $\tau_{Y,Z}$ 
(\cref{fig. pairing comp}(b)(i)-(iv)). Each equivalence class contains precisely zero or two elements of $X \amalg Z$. The classes that contain two elements of $X \amalg Z$ 
-- the \textit{open components} of the composition -- describe the desired pairing on $X \amalg Z$. The remaining equivalence classes -- that describe cycles of elements of $Y$ -- are called \textit{closed components formed by the composition} of $\tau_{X,Y}$ and $\tau_{Y,Z}$.

Likewise, Brauer diagrams $f = (\tau_f, \kclf) \in \BD (l,m)$ and $g =  (\tau_g, \kclf[g]) \in \BD (m,n)$  may be composed vertically to obtain a Brauer diagram $ g \circ f = (\tau_{gf}, \kclf[gf]) \in \BD(l,n)$ with 
\begin{itemize}
	\item the pairing $\tau_{gf}$ is the composition pairing $\tau_g \circ \tau_f$ obtained by identifying $ \Tf = \Sf[g]$ according to $t_{f,i} \mapsto s_{g,i}$;
	\item the number $\kclf[gf]$ of closed components in $g \circ f$ satisfies  
	$ \kclf[gf] = \kclf + \kclf[g]\ + \kcl(\tau_{f}, \tau_{g})$ where $\kcl(\tau_{f}, \tau_{g})$ is the number of closed components \emph{formed} by the composition of $\tau_f$ and $\tau_g$.
\end{itemize}

	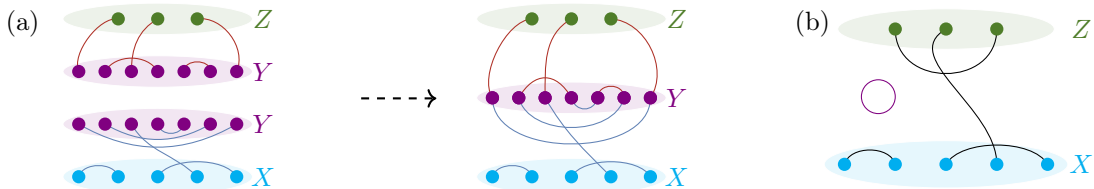
\begin{figure}
		[htb!]
		
		\begin{tikzpicture}
			\node at (-2.5,1){(a)};
			\node at (-.5,0){	\begin{tikzpicture}[scale = .35]
					\begin{pgfonlayer}{above}
						\node at (5,18){\color{sapgreen}$Z$};
						\node at (5,16){\color{violet}$Y$};
						\node at (5,14){\color{violet}$Y$};
						\node at (5,12){\color{cyan}$X$};
						\node [dot, violet]  (12) at (-2, 16) {};
						\node  [dot, violet] (13) at (-1, 16) {};
						\node  [dot, violet] (14) at (0, 16) {};
						\node  [dot, violet] (15) at (1, 16) {};
						\node [dot, violet]  (16) at (2, 16) {};
						\node  [dot, violet] (17) at (3, 16) {};
						\node  [dot, violet] (18) at (4, 16) {};
						\node  [dot, sapgreen] (20) at (-0.5, 18) {};
						\node  [dot, sapgreen](21) at (1, 18) {};
						\node  [dot, sapgreen](22) at (2.5, 18) {};
						\node  [dot, cyan](28) at (-2, 12) {};
						\node   [dot, cyan](29) at (-0.5, 12) {};
						\node   [dot, cyan](30) at (1, 12) {};
						\node   [dot, cyan](31) at (2.5, 12) {};
						\node  [dot, cyan] (32) at (4, 12) {};
						\node [dot, violet]  (33) at (-2, 14) {};
						\node  [dot, violet] (34) at (-1, 14) {};
						\node  [dot, violet] (35) at (0, 14) {};
						\node  [dot, violet] (36) at (1, 14) {};
						\node [dot, violet]  (37) at (2, 14) {};
						\node  [dot, violet] (38) at (3, 14) {};
						\node [dot, violet]  (39) at (4, 14) {};
					\end{pgfonlayer}
					\begin{pgfonlayer}{background}
						\filldraw[draw = white, fill = violet, fill opacity  = .1]	(1,16) ellipse (3.6cm and .6cm); 
						\filldraw[draw = white, fill = violet, fill opacity  = .1]	(1,14) ellipse (3.6cm and .6cm);
						\filldraw[draw = white, fill = sapgreen, fill opacity  = .1]	(1,18) ellipse (3.6cm and .6cm);
						\filldraw[draw = white, fill =cyan , fill opacity  = .1]	(1,12) ellipse (3.6cm and .7cm);
						\draw [red, bend left=300, looseness=0.75] (20.center) to (12.center);
						\draw [red, in=90, out=-150] (21.center) to (14.center);
						\draw [red, bend left=60] (22.center) to (18.center);
						\draw [red, bend left=90, looseness=1.25] (16.center) to (17.center);
						\draw [red, bend left=60] (13.center) to (15.center);
						\draw [blue, bend right=75, looseness=1.25] (36.center) to (37.center);
						\draw [blue, in=135, out=-60] (35.center) to (31.center);
						\draw [blue, bend left=75] (28.center) to (29.center);
						\draw [blue, bend left=60, looseness=0.75] (30.center) to (32.center);
						\draw [blue, bend right] (34.center) to (38.center);
						\draw [blue, bend right] (33.center) to (39.center);
					\end{pgfonlayer}
			\end{tikzpicture}};
			\draw [thick, dashed, -> ] (2,0)--(3,0);
			\node at (5,0){\begin{tikzpicture}[scale = .35]
					\begin{pgfonlayer}{above}
						\node at (5,18){\color{sapgreen}$Z$};
						\node at (5,15){\color{violet}$Y$};
						
						\node at (5,12){\color{cyan}$X$};
						\node [dot, violet]  (12) at (-2, 15) {};
						\node  [dot, violet] (13) at (-1, 15) {};
						\node  [dot, violet] (14) at (0, 15) {};
						\node  [dot, violet] (15) at (1, 15) {};
						\node [dot, violet]  (16) at (2, 15) {};
						\node  [dot, violet] (17) at (3, 15) {};
						\node  [dot, violet] (18) at (4, 15) {};
						\node  [dot, sapgreen] (20) at (-0.5, 18) {};
						\node  [dot, sapgreen](21) at (1, 18) {};
						\node  [dot, sapgreen](22) at (2.5, 18) {};
						\node  [dot, cyan](28) at (-2, 12) {};
						\node   [dot, cyan](29) at (-0.5, 12) {};
						\node   [dot, cyan](30) at (1, 12) {};
						\node   [dot, cyan](31) at (2.5, 12) {};
						\node  [dot, cyan] (32) at (4, 12) {};
						\node [dot, violet]  (33) at (-2, 15) {};
						\node  [dot, violet] (34) at (-1, 15) {};
						\node  [dot, violet] (35) at (0, 15) {};
						\node  [dot, violet] (36) at (1, 15) {};
						\node [dot, violet]  (37) at (2, 15) {};
						\node  [dot, violet] (38) at (3, 15) {};
						\node [dot, violet]  (39) at (4, 15) {};
					\end{pgfonlayer}
					\begin{pgfonlayer}{background}
						\filldraw[draw = white, fill = violet, fill opacity  = .1]	(1,15) ellipse (3.6cm and .6cm); 
						
						\filldraw[draw = white, fill = sapgreen, fill opacity  = .1]	(1,18) ellipse (3.6cm and .6cm);
						\filldraw[draw = white, fill =cyan , fill opacity  = .1]	(1,12) ellipse (3.6cm and .7cm);
						\draw [red, bend left=300, looseness=0.75] (20.center) to (12.center);
						\draw [red, in=90, out=-150] (21.center) to (14.center);
						\draw [red, bend left=60] (22.center) to (18.center);
						\draw [red, bend left=90, looseness=1.5] (16.center) to (17.center);
						\draw [red, bend left=60, looseness=1.5] (13.center) to (15.center);
						\draw [blue, bend right=75, looseness=1.5] (36.center) to (37.center);
						\draw [blue, in=135, out=-60] (35.center) to (31.center);
						\draw [blue, bend left=75] (28.center) to (29.center);
						\draw [blue, bend left=60, looseness=0.75] (30.center) to (32.center);
						\draw [blue, bend right = 75, looseness=1] (34.center) to (38.center);
						\draw [blue, bend right = 90, looseness=1] (33.center) to (39.center);
					\end{pgfonlayer}
				\end{tikzpicture}
			};
			
			\node at (8,1){(b)};
			\node at (10,0){\begin{tikzpicture}[scale = .45]
					\begin{pgfonlayer}{above}
						\node   (12) at (-2, 15) {};
						\node   (13) at (-1, 15) {};
						\node   (14) at (0, 15) {};
						\node   (15) at (1, 15) {};
						\node   (16) at (2, 15) {};
						\node   (17) at (3, 15) {};
						\node   (18) at (4, 15) {};
						\node  [dot, sapgreen] (20) at (-0.5, 16) {};
						\node  [dot, sapgreen](21) at (1, 16) {};
						\node  [dot, sapgreen](22) at (2.5, 16) {};
						\node  [dot, cyan](28) at (-2, 12) {};
						\node   [dot, cyan](29) at (-0.5, 12) {};
						\node   [dot, cyan](30) at (1, 12) {};
						\node   [dot, cyan](31) at (2.5, 12) {};
						\node  [dot, cyan] (32) at (4, 12) {};
						\node   (33) at (-2, 15) {};
						\node   (34) at (-1, 15) {};
						\node   (35) at (0, 15) {};
						\node   (36) at (1, 15) {};
						\node   (37) at (2, 15) {};
						\node  (38) at (3, 15) {};
						\node   (39) at (4, 15) {};
					\end{pgfonlayer}
					\node at (5,16){\color{sapgreen}$Z$};

				\node at (5,12){\color{cyan}$X$};
					\begin{pgfonlayer}{background}

					\filldraw[draw = white, fill = sapgreen, fill opacity  = .1]	(1,16) ellipse (3.2cm and .6cm);
					\filldraw[draw = white, fill =cyan , fill opacity  = .1]	(1,12) ellipse (3.6cm and .7cm);
						\draw [ in=90, out=-150] (21.center) to (31.center);
						
						\draw [ bend right=90, looseness=1.5] (20.center) to (22.center);
					
						\draw [ bend left=75] (28.center) to (29.center);
						\draw [ bend left=60, looseness=0.75] (30.center) to (32.center);
						\draw [violet] (-1, 14) circle (.5cm);
						
					\end{pgfonlayer}
				\end{tikzpicture}
			};
		\end{tikzpicture}
		\caption{(a) Composition of pairings on $X \amalg Y$ and $Y \amalg Z$; 
			(b) the resulting pairing on $X \amalg Z$, together with the single closed component formed in the composition. }\label{fig. pairing comp}
	\end{figure}

	This composition is associative, with two-sided units $(id_n, 0) \in \BD (n,n)$. Hence, we may define:
	
	\begin{defn}\label{def. BD}
		
		The category $\BD$ of \emph{(monochrome unoriented) Brauer diagrams} has objects $n \in \N$, morphism sets
		$\BD(m, n)$ and composition given by vertical composition of Brauer diagrams.

	\end{defn}
	
	The category $\BD$ is a prop with monoidal product (horizontal sum) induced by addition of natural numbers and juxtaposition of Brauer diagrams: for $(\tau_1, \kcl_1)\colon m_1 \to n_1$ and $(\tau_2, \kcl_2) \colon m_2 \to n_2$,
	\[ (\tau_1, \kcl_1) \oplus (\tau_2, \kcl_2) = (\tau_1 \amalg \tau_2, \kcl_1 +\kcl_2) \colon m_1 + m_2 \to n_1 + n_2.\] The monoidal unit is given by the trivial open Brauer diagram $(\varnothing, 0) \colon 0 \to 0$.
		
		Note, in particular, that any Brauer diagram $f = (\tau, \kcl)\colon m \to n$ may be written as a horizontal sum $(\tau, 0) \oplus (\varnothing, \kcl)$ of an open Brauer diagram $(\tau, 0) \colon m \to n$ and a scalar $(\varnothing, \kcl) =\bigoplus_{i = 1}^\kcl (\varnothing, 1)\colon 0 \to 0$. 
		
			Let $id_1 \in \BD(1,1)$, $\cup \in \BD(0,2)$ and $\cap \in \BD(2,0)$ be the morphisms induced by the unique pairing on the two-element set. 
For all $n \in \N$, 
		$ id_n =  \bigoplus_{i = 1}^n id_1\in \BD (n,n)$, and 
		$ \cup_n \defeq \coev{id_n}  \in \BD(0, 2n)$ and $\cap_n \defeq \ev{id_n} \in \BD(2n, 0)$
		satisfy the $n$-fold triangle identities. 
		\begin{equation}\label{eq. n triangle} 
			(\cap_n \oplus id_n) \circ (id_n \oplus \cup_n)= id_n = (id_n\oplus \cap_n )\circ (\cup_n \oplus id_n).
		\end{equation}

	As such, $\BD$ is the free compact closed category generated by one self-dual object. Hence, it has the following universal property:
	
	\begin{lem}
	\label{lem. BD universal}
		For any symmetric monoidal category $ \mathsf C$ and any self-dual object $x \in \mathsf C$, there is a unique symmetric strict monoidal functor $\xi_x \colon \BD \to \mathsf C$ such that $\xi_x(1) = x$.
	\end{lem}

		\begin{rmk}\label{rmk open not subcat}
			
			It is important to note that the subsets $\BDop(m,n) \subset  \BD(m,n)$ of open Brauer diagrams do not describe a subcategory of $ \BD$. Namely, the \emph{unit trace} $tr(id_i) = \bigcirc = \cap \circ \cup $ 
			satisfies $\bigcirc  =(\varnothing, 1)\in \BD(0,0)$ which is not open. 
			
		\end{rmk}

	\begin{ex}
		\label{ex. manifold components} Brauer diagrams may equivalently be defined as tangles in some high (>3)-dimensional space (e.g.,~\cite{Ban16}). In fact, $\BD$ is a skeletal subcategory of the 1-dimensional cobordism category whose morphisms are boundary-preserving isotopy classes of compact 1-manifolds. Hence, monoidal functors from $\BD$ may be referred to as \emph{lax} TQFTs. 
		
		Let $\I$ denote the unit interval $[0,1]$, and let $\man \cong n_o (\I)\amalg n_c (S^1)$ be a compact 1-manifold with canonical pairing $\tau^\man$ on $\partial \man$ as in \cref{ex. manifold pairing}. 
		If $m,n \in \N$ satisfy $m+n = 2 n_o$, and $\phi \colon \{s_1, \dots, s_m\} \amalg \{t_1, \dots, t_n\} \to \partial \man$ is any isomorphism, then 
		$(\phi^{-1} \tau^\man \phi, n_c) \in \BD(m,n)$. 
		Conversely, given a morphism $ f = (\tau, \kcl )\in \BD(m,n)$, there is a unique (up to boundary-preserving isotopy) 
		compact 1-manifold $\man_f\cong \frac{m+n}{2}(\I) \amalg \kcl (S^1)$ and isomorphism $\phi_f \colon 
		\Sf \amalg \Tf \to \partial \man$ such that $\phi^{-1}_f \tau^{\man_f} \phi_f = \tau$. 
		
	\end{ex}

Let $f = (\tau, \kcl)\in \BD(m,n)$. Following \cref{ex. manifold components}, $ \partial f  \defeq \Sf \amalg \Tf$ is called the \textit{boundary of $f = (\tau, \kcl)\in \BD(m,n)$}. 
A \emph{component of $f$ } is an element of the set $\pi_0(f)$ of connected components of a compact manifold $\man_f$ as in \cref{ex. manifold components}. So, $|\pi_0(f)|  = \frac{(m+n)}{2} + \kcl$.

		There is a canonical map $\partial f \to \pi_0 (f)$ so that $f $ is described by a diagram of cospans of finite sets:
		\begin{equation}\label{eq. cospan def BD} \xymatrix{  \Sf \ar@{-->}[rrd] \ar@{>->}[rr]^{S_i \ \mapsto \tau (S_i)}&& \partial f \ar[d] && \ar@{>->}[ll]_-{t_j \ \mapsto \tau(t_j) }\Tf  \ar@{-->}[lld]\\&& \pi_0(f).&&}\end{equation}

		\begin{rmk}\label{rmk. pushout}

			By (\ref{eq. cospan def BD}), for composable morphisms $f \in \BD(k,m)$ and $g \in \BD(m,n)$, we may consider the pushout diagram: 
			\begin{equation}\label{eq. cospan pushout} 
				\xymatrix@C=.6cm {
					\Sf\ar[drr]_-{\tau_f} &&&& \Tf =  \Sf[g]  \ar[dll]^-{\tau_f} \ar[drr]_-{\tau_g}&&&& \Tf[g] \ar[dll]^-{\tau_g}\\
					&& \partial f \ar@{-->}[drr]\ar[d] &&&&  \partial g\ar@{-->}[dll]\ar[d] &&\\
					&& \pi_0 (f)  \ar@{-->}[drr]&& P(gf)\ar[d] && \pi_0 (g)\ar@{-->}[dll]&& \\
					&&&& \pi^P(gf). &&&&
				}
			\end{equation}

			However, $\BD$ is not a cospan category since, in general, $P(gf) \not \cong \partial (gf) = \Sf \amalg \Tf[g]$ and hence composition of morphisms in $\BD$ is not described by compositions (pushouts) of cospans as in (\ref{eq. cospan pushout}). 
			
			For example, in 
			the pushout (\ref{eq. cospan pushout}) for the composition $\cap \circ \cup = \bigcirc$, $P(\cap \circ \cup)$ has two elements, but 
			$\partial \bigcirc = \emptyset$. This is equivalent to the observation that open Brauer diagrams do not describe a subcategory of $\BD$ and is closely related to the \textit{problem of loops} discussed in detail in \cite[Section~6]{Ray20}. 
			
	\end{rmk}
	
	By e.g.,~\cite[Theorem~2.6]{LZ15} or \cite[Proposition~2.15]{Ban16}, the category $\BD$ is generated, under horizontal and vertical composition, by the open morphisms $id_1$, $\cup$, $\cap$ and the unique non-identity permutation $\sigma_{\two} \in \Sigma_2 \subset  \BD(2,2)$, with the obvious identity, symmetry and triangle relations (\cref{fig. triangle 1}). Interesting subcategories of $ \BD$ may be obtained by taking subsets of the generating set.
	
		\begin{defn}\label{def. dw BD}
		The category $\BDd \subset  \BD$ of \emph{downward Brauer diagrams} is the subcategory of open morphisms $  (\tau^{\downarrow}, 0)\in \BD(m,n)$ such that, for all $y \in \Tf$, 
		$\tau^{\downarrow}(y)  \in  \Sf$. 
		
		The category $\BDu \subset  \BD$ of \emph{upward Brauer diagrams} is the opposite category of $\BDd$. 
	\end{defn}

		The category $\BDd$ is generated by $id_1, \sigma_{\two}$, and $\cap$ (and $\BDu$ is generated by $id_1, \sigma_{\two}$, and $\cup$) under horizontal and vertical composition, according to the relations in $\BD$.
	In particular, $\BDd(m,n) $ is empty whenever $n >m$, so $\cup$ is not a morphism in $\BDd$ (and $\cap$ is not a morphism in $\BDu$). 
	Since morphisms in $\BDd$ are open, $\BDd(m,n)$ is finite for all $m,n$. Moreover, composition in $\BDd$ (respectively $\BDu$) may be described by pushouts of cospans as in (\ref{eq. cospan pushout}).
	
In \cite{SS15} and \cref{ssec nonunital CA}, representations of the infinite orthogonal and symplectic groups are described in terms of $\BDd$ and, in \cref{def nonunital ca}, $\BDd$ is used to define nonunital circuit algebras. (See also \cite[Section~5]{RayCA2}.)

\begin{rmk}
	\label{rmk TL subcategory}
	Other interesting subcategories of $\BD$ may be obtained by restricting to different subsets of the generating morphisms. Of course, the intersection of $\BDd$ and $\BDu$ in $\BD$ is the permutation groupoid $\Sigma $  generated by $id_1$ and $ \sigma_{\two}$. The \emph{Temperley-Lieb category} $\mathsf{TL} \subset \BD$ is the subcategory of \emph{planar Brauer diagrams} generated by $id_1, \cup, \cap$, but not the symmetry morphism $\sigma_{\two}$ (see \cite{GW02}).
\end{rmk}

		\subsection{Coloured Brauer diagrams, orientations and wheeled props}\label{ssec. colour BD}
		
Generalisations of categories of Brauer diagrams are obtained by {colouring} the diagram components. By considering involutions on colours, the same constructions also serve to describe (coloured) oriented, and mixed Brauer diagrams. (See also \cite{DCH19,Ray20}.)

	\begin{defn} \label{def. palette}
		
		A pair $(\CCC, \omega )$ of a set $\CCC$ together with an involution $\omega \colon \CCC \to \CCC$ 
		is called an \emph{(involutive) palette}. 
	Elements $c \in \CCC$ are called \emph{colours in $(\CCC, \omega )$}. The set of orbits of $\omega$ in $\CCC$ is denoted by $\widetilde \CCC$.
	
	For any palette $(\CCC, \omega)$, there is an induced\emph{ free monoid palette} $(\listm (\CCC), \bfom)$ with involution \begin{equation}\label{eq list inv}
		\bfom \colon (c_1, \dots, c_n) \mapsto (\omega c_n , \dots, \omega c_1). 
	\end{equation}
	
	Objects of the category $\pal$ are palettes $(\CCC, \omega)$, and morphisms $(\CCC, \omega) \to (\CCC', \omega')$ are given by morphisms $\lambda \in \Set(\CCC, \CCC')$ such that $\lambda \circ \omega = \omega' \circ \lambda$.

\end{defn}

Now let $(\CCC, \omega)$ be any palette and $(X, \tau)$ be the palette described by a pairing $\tau$ on a finite set $ X $. 
\begin{defn}

\label{def. colouring}

A \emph{$(\CCC, \omega)$-colouring of $\tau$} is a morphism $\lambda \colon (X, \tau )\to (\CCC, \omega)$ in $\pal$. 

A \emph{$(\CCC, \omega)$-colouring $\lambda$ of a Brauer diagram $f  = (\tau, \kcl)\in \BD(m, n)$} is given by a pair $\lambda = (\colop, \colcl)$ where $\colop$ is a colouring of $\tau$ and $\colcl$ is a map $\pi_0(f)\to \widetilde {\CCC}$ such that the following diagram of sets commutes: 

\begin{equation}\label{eq. colour} \xymatrix{ \partial f \ar[rr]^-{\lambda_\partial }\ar[d]^{\cong}_{\tau}&& \CCC \ar[d]^{\omega}_{\cong}\\
		\partial f  \ar[rr]^-{\lambda_\partial }\ar@{->}[d]&& \CCC \ar@{->>}[d]\\
		\pi_0(f) \ar[rr]^-{\widetilde \lambda} && \widetilde \CCC.}\end{equation} 

The \emph{type of the colouring $\lambda$} 
is the pair $(\ccc, \ddd)\in (\listm (\CCC))^2$ -- where $\ccc$ is called the \emph{input type}, and $\ddd$ is called the \emph{output type}, of $(f, \lambda)$ -- defined by:
\begin{equation}\label{eq. type of colouring}
	\ddd = (d_1, \dots, d_n)  =  \lambda_{\partial}( \Tf), \text{ and } \ccc = (c_1, \dots, c_m)  =\omega \circ \lambda_{\partial}(\Sf ).
\end{equation}
\end{defn}

\begin{rmk}\label{rmk target involution}

The application of $\omega$ in the definition of the input type $\ccc = \omega \circ \colop (\Sf)$ is necessary to define categorical composition of coloured Brauer diagrams in \cref{def. colour BD}.

\end{rmk}

Given $\ccc = (c_1, \dots, c_m)$ and $ \ddd =(d_1, \dots, d_n) $ in $ \listm (\CCC)$, objects of the set $\CBD(\ccc, \ddd)$ of \emph{$(\CCC, \omega)$-coloured Brauer diagrams from $\ccc$ to $\ddd$} are pairs $(f, \lambda)$ where $f = (\tau, \kcl)$ is a morphism in $\BD(m,n)$, and $\lambda$ is a colouring of $f$ of type $(\ccc, \ddd)$.

Horizontal composition $\oplus$ of coloured Brauer diagrams $(f, \lambda ) \in \CBD(\ccc_1, \ddd_1)$ and $(g, \gamma) \in \CBD (\ccc_2, \ddd_2)$ is given by juxtaposition and concatenation:
\[ (f, \lambda) \oplus (g, \gamma) = (f \oplus g, \lambda \amalg \gamma)\in \CBD (\ccc_1\ccc_2, \ddd_1\ddd_2).\]

To define vertical composition, let $(f, \lambda) \in \CBD (\bbb, \ccc)$ and $(g, \gamma) \in \CBD (\ccc, \ddd)$ with $ f = (\tau_f, \kclf) \in \BD(k,m)$ and $g = (\tau_g, \kclf[g]) \in \BD(m,n)$ be such that $ gf  = (\tau_{gf}, \kclf[gf]) \in \BD(k,n)$. By definition, $\gamma_\partial (y) = \omega \lambda_\partial (y)$ for each $y \in \Tf  = \Sf[g]$. So 
$\lambda$ and $\gamma$ induce a well-defined colouring $\gamma \lambda$ on $g \circ f$. 

\begin{figure}
[htb!]
\begin{tikzpicture}
	\node at (0,0){\begin{tikzpicture}[scale = .6]
			\begin{pgfonlayer}{above}
				\node [dot, red]  (12) at (-2, 15) {};
				\node [left, above] at (-2, 15) {\tiny{$\omega c$}};
				\node  [dot, sapgreen] (13) at (-1, 15) {};
				\node [right, above] at (-1, 15) {\tiny{$d$}};
				\node [right, below] at (-1, 15) {\tiny{$\omega d$}};
				\node  [dot, blue] (14) at (0, 15) {};
				\node  [dot, sapgreen] (15) at (1, 15) {};
				\node [right, above] at (1, 15) {\tiny{$\omega d$}};
				\node [right, below] at (1, 15) {\tiny{$d$}};
				\node [dot, sapgreen]  (16) at (2, 15) {};
				\node [right, above] at (2, 15) {\tiny{$ d$}};
				\node [right, below] at (2, 15) {\tiny{$\omega d$}};
				\node  [dot, sapgreen] (17) at (3, 15) {};
				\node [right, above] at (3, 15) {\tiny{$\omega d$}};
				\node [right, below] at (3, 15) {\tiny{$d$}};
				\node  [dot, red] (18) at (4, 15) {};
				\node [left,above] at (4, 15) {\tiny{$ c$}};
				\node  [dot, red] (20) at (-0.5, 18) {};
				\node [left, below] at (-0.5, 18) {\tiny{$ c$}};
				\node  [dot, blue](21) at (1, 18) {};
				\node  [dot, red](22) at (2.5, 18) {};
				\node [right,below] at (2.5,18) {\tiny{$\omega c$}};
				\node  [dot, violet](28) at (-2, 12) {};
				\node   [dot, violet](29) at (-0.5, 12) {};
				\node   [dot, red](30) at (1, 12) {};
				\node [left, above] at (1,12) {\tiny{$ \omega c$}};
				\node   [dot, blue](31) at (2.5, 12) {};
				\node  [dot, red] (32) at (4, 12) {};
				\node [right, above] at (4,12) {\tiny{$  c$}};
				\node [dot, red]  (33) at (-2, 15) {};
				\node [left,below] at (-2, 15) {\tiny{$ c$}};
				\node  [dot, sapgreen] (34) at (-1, 15) {};
				\node  [dot, blue] (35) at (0, 15) {};
				\node  [dot, sapgreen] (36) at (1, 15) {};
				\node [dot, sapgreen]  (37) at (2, 15) {};
				\node  [dot, sapgreen] (38) at (3, 15) {};
				\node [dot, red]  (39) at (4, 15) {};
				\node [left,below] at (4, 15) {\tiny{$\omega c$}};
				
			\end{pgfonlayer}
			\begin{pgfonlayer}{background}
				\draw [red, bend left=300, looseness=0.75] (20.center) to (12.center);
				\draw [blue, in=90, out=-150] (21.center) to (14.center);
				\draw [red, bend left=60] (22.center) to (18.center);
				
				\draw [sapgreen, bend left=90, looseness=1.5] (16.center) to (17.center);
				\draw [sapgreen, bend left=60, looseness=1.5] (13.center) to (15.center);
				\draw [sapgreen, bend right=75, looseness=1.5] (36.center) to (37.center);
				\draw [blue, in=135, out=-60] (35.center) to (31.center);
				\draw [violet, bend left=75] (28.center) to (29.center);
				\draw [red, bend left=60, looseness=0.75] (30.center) to (32.center);
				fill opacity=0.2,  bend right = 75, looseness=1] (34.center) to (38.center);
				\draw [sapgreen, bend right = 75, looseness=1] (34.center) to (38.center);
				\draw [red, bend right = 90, looseness=1] (33.center) to (39.center);
			\end{pgfonlayer}
		\end{tikzpicture}
	};
	\draw [thick, dashed, -> ] (3.4,0)--(4.5,0);
	\node at (7,0){\begin{tikzpicture}[scale = .5]
			\begin{pgfonlayer}{above}
				
				\node  [dot, red]  (20) at (-0.5, 16) {};
				\node [left, below] at (-0.7, 16) {\tiny{$ c$}};
				\node  [dot, blue](21) at (1, 16) {};
				
				\node  [dot, red] (22) at (2.5, 16) {};
				\node [right, below] at (2.8, 16) {\tiny{$\omega c$}};
				
				\node  [dot, violet](28) at (-2, 12) {};
				
				\node   [dot, violet](29) at (-0.5, 12) {};
				
				\node  [dot, red] (30) at (1, 12) {};
				\node [right, above] at (1, 12) {\tiny{$\omega c$}};
				\node   [dot, blue](31) at (2.5, 12) {};
				
				\node  [dot, red]  (32) at (4, 12) {};
				\node [right, above] at  (4, 12) {\tiny{$c$}};
				
			\end{pgfonlayer}
			\begin{pgfonlayer}{background}
				\draw [blue, in=90, out=-150] (21.center) to (31.center);
				\draw [red, bend right=90, looseness=1.5] (20.center) to (22.center);
				\draw [violet, bend left=75] (28.center) to (29.center);
				\draw [red, bend left=60, looseness=0.75] (30.center) to (32.center);
				\draw [sapgreen] (4, 14) circle (.5cm);
				\node at (4.8,14){\tiny{$[d]$}};
			\end{pgfonlayer}
		\end{tikzpicture}
	};
\end{tikzpicture}
\caption{Composing coloured pairings.} \label{fig. composing colours}
\end{figure}
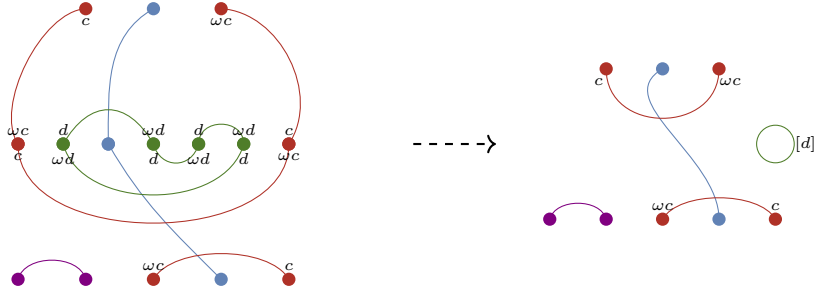

\begin{defn}
\label{def. colour BD}
Objects of the category $\CBD$ of \emph{$(\CCC, \omega)$-coloured Brauer diagrams} are elements of $\listm (\CCC)$. Morphisms in $\CBD(\ccc, \ddd)$ are $(\CCC, \omega)$-coloured Brauer diagrams of type $(\ccc, \ddd)$, with composition of morphisms $(f, \lambda) \in \CBD (\bbb, \ccc)$ and $(g, \gamma) \in \CBD (\ccc, \ddd)$ is given by $(gf, \gamma \lambda) \in \CBD(\bbb,\ddd)$.

\end{defn}

\begin{rmk}
\label{rmk coloured pairings} Let $\ccc = (c_1, \dots, c_m), \ddd  = (d_1, \dots, d_n)$ and  let $ (f, \lambda) \in \CBD(\ccc,\ddd)$ be a morphism with underlying Brauer diagram $f = (\tau , \kclf) \in \BD(m,n)$. 

The pairing $\tau$ induces a pairing on $\{c_1, \dots, c_m\} \amalg \{d_1, \dots, d_n\}$ in the obvious manner and $\tilde \lambda$ describes an unordered $\kclf$-tuple in $\widetilde \CCC$. Hence, a $ (f, \lambda) \in \CBD(\ccc,\ddd)$ may also be denoted simply by $(\tau, \tilde \lambda)$. 
\end{rmk}

The category $\CBD$ is a $\CCC$-coloured prop (see \cref{ssec. SMC}), with monoidal structure $\oplus$ 
induced by concatenation of object lists and disjoint union of coloured Brauer diagrams. It has a compact closed structure given by $\ccc^* =  \bfom(\ccc)$ for all $\ccc$. 

\begin{rmk}
	\label{rmk chromatic Brauer} When $\omega = id_{\CCC}$ is the identity, $\CBD$ is a category of nonoriented $\CCC$-coloured Brauer diagrams, called a \emph{chromatic Brauer category} in \cite{MW21}. Extending \cite{Ban16}, these are used in \cite{MW21} to distinguish exotic smooth spheres.
\end{rmk}

Of particular importance is the palette $\dipal$ given by the unique non-trivial involution $(\uparrow) \leftrightarrow  (\downarrow)$ on the two-element set $\{\uparrow , \downarrow\}$. A $\dipal$-coloured Brauer diagram is called \emph{oriented} and $\DiBD \defeq \BD^{\dipal}$ is the category of \textit{ (monochrome) oriented Brauer diagrams}. Objects of $\DiBD$ are 
finite words in the alphabet $\{\uparrow , \downarrow\}$. Let $ \uparrow ^n$ (respectively $\downarrow ^n$) denote the object of $\DiBD$ given by $ n$ copies of $\uparrow$ (respectively $\downarrow$) in $ \listm \{\uparrow , \downarrow\}$. So objects of $\DiBD$ are concatenations of words of the form $ \uparrow ^m$ and $ \downarrow ^n$. 
Morphisms in $\DiBD$ are represented, as in \cref{fig. composing directions}, by diagrams of oriented intervals and (unoriented) circles.

 More generally, if $\DDD$ is a set, and $\CCC = \DDD \times\dipal$, then the category $\DiCBD \defeq \CBD$ of $\DDD$-coloured oriented Brauer diagrams is the free compact closed prop generated by elements of the set $\DDD$ and their formal duals. 
 	For $\ddd=  (d_1, \dots, d_n)\in \listm (\DDD)$, let $ \uparrow ^{\ddd}$ (respectively $ \downarrow ^{\ddd}$) denote $((d_1, \uparrow), \dots , (d_n, \uparrow )) \in \listm (\DDD \times \dipal)$.  If $ \downarrow ^{\ddd}$ is defined similarly, then objects of $\DiCBD$ are concatenations of words of the form $\uparrow ^{\ddd}$ and $  \downarrow ^{\ccc} $. Note that $\bfom(\uparrow^{\ddd}) = \downarrow^{\ddd^\dagger}$ where $\ddd^\dagger \defeq (d_n, \dots, d_1)$.

\begin{ex}
\label{ex. wheeled prop OBD}
The full subcategory $\mathsf W^{\DDD} \subset   \DiCBD$ on objects of the form $\uparrow ^{\ddd}$ is canonically a $\DDD$-coloured wheeled prop. But it is not compact closed, since $\mathsf W^{\DDD}$ does not admit duals. 

Applying the Int construction \cite{JSV96} to $ \mathsf W^{\DDD}$ results in the category $ W \CBD[\DDD]$ of \emph{$\DDD$-coloured walled Brauer diagrams}. This is the full subcategory of $\DiCBD$ on objects of the form $ \uparrow ^{\ccc}\downarrow ^{\ddd} $, $\ccc, \ddd \in \listm (\DDD)$. The inclusion $ W \CBD[\DDD] \hookrightarrow \DiCBD$ is an equivalence of categories since every object of $\DiCBD$ is isomorphic -- via a canonical shuffle permutation -- to a unique object of $W\CBD [\DDD]$ (see \cref{fig. composing directions}). (\emph{Walled Brauer algebras} were introduced independently in \cite{Koi89, Tur90}.) 
\end{ex}

\begin{figure}
	[htb!]
		\centerfloat
	\begin{tikzpicture}
			\node at (2.3,-2){(a)};
		\node at (0,0) {
			
			\begin{tikzpicture}
				\node at (9,-4){\begin{tikzpicture}[scale = .4]
						\begin{pgfonlayer}{above}
							\node [dot ]  (12) at (-2, 15) {};
							\node  [dot ](13) at (-1, 15) {};
							\node  [dot ] (14) at (0, 15) {};
							\node  [dot ](15) at (1, 15) {};
							\node [dot ] (16) at (2, 15) {};
							\node  [dot ] (17) at (3, 15) {};
							\node  [dot ] (18) at (4, 15) {};
							
							\node at (-4, 16.5){$f\colon$};
							\node  [dot ] (20) at (-0.5, 18) {};
							\node [dot ](21) at (1, 18) {};	
							\node  [dot ](22) at (2.5, 18) {};
							\node at (-4, 13.5){$g\colon$};		
							\node  [dot ](28) at (-2, 12) {};
							\node   [dot ](29) at (-0.5, 12) {};
							\node   [dot ](30) at (1, 12) {};
							\node   [dot ](31) at (2.5, 12) {};
							\node  [dot ](32) at (4, 12) {};

						\end{pgfonlayer}
						\begin{pgfonlayer}{background}
							\draw [bend left=300, looseness=0.75, red, ->-=.5] (20.center) to (12.center);
							\draw [ in=90, out=-150, red, ->-=.5] (21.center) to (14.center);
							\draw [bend left=60, red, -<-=.5] (22.center) to (18.center);
							\draw [bend left=90, looseness=1.5, red, ->-=.5] (16.center) to (17.center);
							\draw [ bend left=60, looseness=1.5,  red, ->-=.6] (13.center) to (15.center);
							\draw [ bend right=75, looseness=1.5, red, ->-=.5] (15.center) to (16.center);
							\draw [ in=135, out=-60, red, ->-=.5] (14.center) to (31.center);
							\draw [bend left=75, , red, -<-=.48] (28.center) to (29.center);
							\draw [bend left=60, looseness=1,  red, ->-=.5] (30.center) to (32.center);
							\draw [ bend right = 75, looseness=1, red, -<-=.5] (13.center) to (17.center);
							\draw [bend right = 90, looseness=1, red, ->-=.5] (12.center) to (18.center);
						\end{pgfonlayer}
					\end{tikzpicture}
				};
			\end{tikzpicture}
		};
			\draw [thick, dashed, -> ] (2.5,0)--(3.2,0);
		\node at (4.5,0){\begin{tikzpicture}[scale = .4]
				\begin{pgfonlayer}{above}
					
					\node  [dot]  (20) at (-0.5, 16) {};
					
					\node  [dot](21) at (1, 16) {};
					
					\node  [dot] (22) at (2.5, 16) {};

					\node  [dot](28) at (-2, 12) {};
					
					\node   [dot](29) at (-0.5, 12) {};
					
					\node  [dot] (30) at (1, 12) {};
				
					\node   [dot](31) at (2.5, 12) {};
					
					\node  [dot]  (32) at (4, 12) {};

				\end{pgfonlayer}
				\begin{pgfonlayer}{background}
					\draw [red, ->-=.6, in=90, out=-150] (21.center) to (31.center);
					\draw [red, ->-=.6, bend right=90, looseness=1.5] (20.center) to (22.center);
					\draw [red, -<-=.5, bend left=75] (28.center) to (29.center);
					\draw [red, ->-=.6, bend left=60, looseness=0.75] (30.center) to (32.center);
					\draw [red] (4, 14) circle (.5cm);
				\end{pgfonlayer}
			\end{tikzpicture}
		};
			\node at (11,-2){(b)};
		\node at (9, 0) {
			\begin{tikzpicture}[scale = .4]
				\begin{pgfonlayer}{above}
					\node [dot, blue] (0) at (0, 0) {};
					\node [dot, blue] (1) at (-1, 0) {};
					\node [dot, blue] (2) at (-2, 0) {};
					\node [dot, blue] (3) at (-3, 0) {};
					\node [dot, red] (4) at (1, 0) {};
					\node [dot, red] (5) at (2, 0) {};
					\node [dot, red] (6) at (3, 0) {};
					\node [dot, blue] (7) at (-2.25, 3) {};
					\node [dot, blue] (8) at (-0.75, 3) {};
					\node [dot, red] (9) at (1, 3) {};
					\node [dot, blue] (10) at (0, -3) {};
					\node [dot, blue] (11) at (-1.5, -3) {};
					\node [dot, blue] (12) at (-3, -3) {};
					\node [dot, red](13) at (1.5, -3) {};
					\node [dot, red] (14) at (3, -3) {};
				\end{pgfonlayer}
				\begin{pgfonlayer}{background}
					\draw[dashed, gray](.5, 3.5)--(0.5, -3.5);
					\draw [bend right=15] (7.center) to (3.center);
					\draw [bend right] (8.center) to (2.center);
					\draw [bend left] (9.center) to (6.center);
					\draw [bend right=15, looseness=1.25] (2.center) to (11.center);
					\draw [bend right=90] (0.center) to (4.center);
					\draw [bend left=90, looseness=0.75] (10.center) to (14.center);
					\draw [bend right=90, looseness=0.75] (1.center) to (5.center);
					\draw [in=270, out=-45, looseness=0.75] (3.center) to (6.center);
					\draw [bend left=90, looseness=0.75] (0.center) to (5.center);
					\draw [bend left=90, looseness=0.75] (1.center) to (4.center);
					\draw [bend left=90, looseness=0.75] (12.center) to (13.center);
				\end{pgfonlayer}
			\end{tikzpicture}
		};
			\draw [thick, dashed, -> ] (10.8,0)--(11.5,0);
			\node at (13, 0) {
			\begin{tikzpicture}[scale = .4]
				\begin{pgfonlayer}{above}
					\node [dot, blue] (7) at (-2.25, 2) {};
					\node [dot, blue] (8) at (-0.75, 2) {};
					\node [dot, red] (9) at (1, 2) {};
					\node [dot, blue] (10) at (0, -2) {};
					\node [dot, blue] (11) at (-1.5, -2) {};
					\node [dot, blue] (12) at (-3, -2) {};
					\node [dot, red](13) at (1.5, -2) {};
					\node [dot, red] (14) at (3, -2) {};
				\end{pgfonlayer}
				\begin{pgfonlayer}{background}
					\draw[dashed, gray](.5, 2.5)--(0.5, -2.5);
				
					\draw [bend right] (8.center) to (11.center);

					\draw [bend left=90, looseness=0.75] (10.center) to (14.center);
					\draw [bend right=90, looseness=0.75] (7.center) to (9.center);
				
					\draw [bend left=90, looseness=0.75] (12.center) to (13.center);
					\draw (.5, 0) circle (.5cm);
				\end{pgfonlayer}
			\end{tikzpicture}
		};

	\end{tikzpicture}
	\caption{(a) Composing oriented Brauer diagrams.  
		(b) Up to a shuffle permutation, this is equivalent to a composition of walled Brauer diagrams, where horizontal arrows go from left to right.} 
	\label{fig. composing directions}
\end{figure}
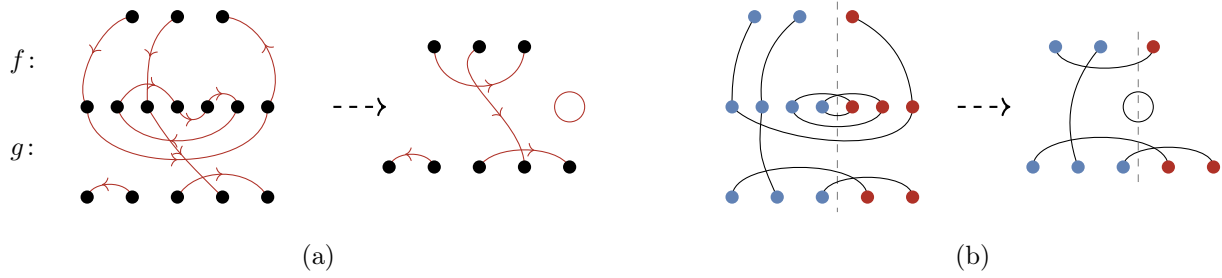

In fact, the category $\DiCBD$ classifies $\DDD$-coloured wheeled props: 
\begin{prop}\label{prop: wheeled prop oriented brauer}
There is an equivalence of categories between the category $\WP_{\V}^{\DDD}$ of $\DDD$-coloured wheeled props in a symmetric monoidal category $(\V, \otimes, I)$ and the category $[\DiCBD[\DDD] ,\V]_{\mathrm{lax}}$ of symmetric monoidal functors $\DiCBD[\DDD] \to \V$ and natural transformations that commute with the structure maps.
\end{prop}

\begin{proof}

Let $ (\alg, \pi, \eta) \colon \DiCBD[\DDD] \to \V$ be a symmetric monoidal functor. This describes a $\DDD$-coloured wheeled prop $(P_{\alg}, \otimes^P, \varnothing _{\DDD}, tr^P)$ as follows:

For $ \ccc, \ddd \in \listm (\DDD)$, $P_{\alg} (\ccc, \ddd) \defeq \alg (\uparrow ^{\ccc} \downarrow^{ \ddd^\dagger})$. The symmetric structure of $(\alg, \pi)$ induces symmetry isomorphisms in $P_{\alg}$. The monoidal (horizontal) composition $\otimes_P$ on morphisms in $P_{\alg}$ is obtained from $\pi$ by composition with the appropriate symmetry (shuffle) isomorphism and has monoidal unit $1_{\varnothing} = \eta\colon I_{\V} \to P_{\alg}(\varnothing, \varnothing) = \alg(\varnothing)$.

For each $ d \in \DDD$, the identity $ 1_d \colon I \to P_{\alg}(d,d)$ is given by the composition 
\[I \xrightarrow {\eta} \alg (\varnothing) \xrightarrow {\alg (\cup_{\uparrow ^d})} \alg( \uparrow^d \downarrow ^d ) = P_{\alg}(d,d).\] 

Categorical (vertical) composition $P_{\alg} (\bbb, \ccc) \otimes P_{\alg} (\ccc, \ddd) \to P_{\alg} (\bbb, \ddd)$ in $P_{\alg}$ is defined by 
\[\alg (id_{\uparrow ^{\bbb}} \otimes \cap_{\uparrow ^{\ccc}} \otimes id_{\downarrow ^{\ddd^\dagger}})\circ  \pi \colon  \alg (\uparrow ^{\bbb} \downarrow^{ \ccc^\dagger}) \otimes \alg (\uparrow ^{\ccc} \downarrow^{ \ddd^\dagger}) \to \alg (\uparrow ^{\bbb} \downarrow^{ \ddd^\dagger}),  \] and likewise the trace 
is given by 
\[ \xymatrix{
	P_{\alg} (\ccc \bbb,\ddd \bbb) \ar[rrr]^{tr_{\ccc,\ddd}^{\bbb}}\ar@{=}[d] &&& P_{\alg} (\ccc, \ddd) \\
		\alg (\uparrow^{\ccc} \uparrow^{\bbb}\downarrow^{\bbb^\dagger}\downarrow^{\ddd^\dagger}) 
		\ar[rrr]_{\alg (id_{\uparrow ^{\ccc}} \otimes \cap_{\downarrow ^{\bbb^\dagger}} \otimes id_{\downarrow ^{\ddd^\dagger}})} &&&\alg (\uparrow^{\ccc}\downarrow^{\ddd^\dagger}) \ar@{=}[u]
}\]

It follows immediately from the relations in $\DiCBD[\DDD]$ that $P_{\alg}$ satisfies the axioms for traced monoidal categories (\cref{ssec: duality}).

Conversely, let $(P, \otimes_P,\varnothing_{\DDD}, tr)$ be a $\DDD$-coloured wheeled prop in $\V$. Define $(\alg_P, \pi_P, \eta_P) \colon \DiCBD[\DDD] \to \V$ by $\alg_P (\uparrow ^{\ccc} \downarrow ^{\ddd}) \defeq P (\ccc, \ddd^\dagger)$. The symmetric action on $\alg_P$ is induced by symmetry in $P$. Permutations that shuffle $\uparrow ^{\ccc}$ with $ \downarrow ^{\ddd^\dagger}$ act trivially on $\alg_P$.

Horizontal composition $\otimes _P$ in $P$ induces a lax multiplication $\pi_P$ on $\alg$ with a lax unit for $\alg_P$ described by the unit morphism $ 1_{\varnothing} \colon I_V \to P(\varnothing, \varnothing)  = \alg_P(\varnothing)$. 

For each $d \in \DDD$, \[\alg (\cap_{\downarrow^d}) = tr^d: \alg(\uparrow^d \downarrow^d)\to \alg (\varnothing)\] and $\alg (\cup_{\downarrow^d}) $ is given by
\[\xymatrix{\alg (\varnothing_{\CCC}) = P(0,0) \ar[rr] \ar[dr]_{tr(id_0)}&&P(d,d) = \alg(\uparrow^d \downarrow^d).\\
	& I \ar[ur]_{1_d}&} \]

Since $P$ satisfies the wheeled prop axioms (\cref{defn WP} \&~ \cite{JSV96}), $\alg_P$ satisfies the relations in $\DiCBD$, and hence defines a symmetric lax functor from $\DiCBD$. 

The assignments $\alg \to P_{\alg}$ and $ P \to \alg_P$ preserve all defining structure and are each others' inverses up to shuffle isomorphisms in $\DiCBD$. Hence $\WP_{\V}^{\DDD}\simeq[\DiCBD[\DDD] ,\V]_{\mathrm{lax}}$.
\end{proof}

By \cite{DHR20}, this result will also follow from \cref{thm. lax functor ca}.

	\subsection{Representations of $\BD$ and $\DiBD$}\label{ssec. representations of BD}\label{ex: Brauer category}

This short section reviews some known results in the representation theory of (oriented) Brauer diagrams.

	Let $R$ be a commutative ring, and $\modR$ its category of modules. For $\delta \in R$, let $\Br = \Br^R$ be the $\modR$-enriched \textit{Brauer category (with specialisation $\delta$)} defined in \cite{LZ15}, whose objects are natural numbers $n \in \N$ and, 
	for all $m, n \in \N$, $\Br(m,n) $ is the free $R$-module (finitely) generated by the open Brauer diagrams $\tau \in \BDop(m,n)$. If $\tau_f \in \Br(k,m)$ and $ \tau_g \in \Br(m,n)$ are generating morphisms, then their composition in $\Br$ is defined by 
$ \tau_g \tau_f =\delta^{\kclf[gf]}\tau_{gf}\in \Br(k,n).$ 
	In particular, $ \Br(0,0)  = \langle \delta \rangle \subset   R$ is the ideal generated by $\delta$.

		Let $\underline{\BD}_R$ be the free $\modR$-category on $\BD$. So, for each pair $m, n $ of natural numbers, $\underline{\BD}_R(m,n)$ is the free $R$ module (infinitely) generated by $\BD(m,n)$. There is a canonical isomorphism $ \underline{\BD} \cong \Br[t]^{R[t]}$ of $\modR$-enriched categories given by $(\tau, \kcl)\leftrightarrow t^{\kcl}\tau$. 
		For each $\delta \in \R$, let $T_\delta \colon \BD \to \Br$ be the obvious identity-on-objects symmetric (strict) monoidal functor such that $\bigcirc = \cap \circ \cup \mapsto \delta$. This factors through the symmetric strict monoidal ($\modR$)-enriched \emph{specialisation} functor $  \Br[t]^{R[t]}\to  \Br$ induced by $t \mapsto \delta$.

	In the oriented case, let $\DiBr$ be the \emph{oriented Brauer category (with specialisation $\delta$)} defined similarly to $\Br$ (but with oriented Brauer diagrams). 
	In particular, 
	the free $\modR$ category $ \underline{\DiBD}$ on $ \DiBD$ is isomorphic to $\DiBr[t]^{R[t]}$. As in the unoriented case, for each $\delta \in \R$, the obvious identity-on-objects symmetric (strict) monoidal functor $ \DiBD 
	\to \DiBr$ such that $\bigcirc \mapsto \delta$ is denoted by $T_\delta$.   

If $(\alg, \pi, \eta) \colon \BD \to \modR$ is a symmetric monoidal functor, then $ \alg(0)$ is an $R$-algebra with unit $ \eta $ 
and algebra multiplication $\pi$. For $r \in R$, it is convenient to denote $\eta (r)\in \alg (0)$ simply by $r$.
	\begin{lem}
		\label{lem universal property}
	A symmetric monoidal functor $\alg \colon \BD \to \modR$ factors through $T_\delta$ if and only if  $\alg (\bigcirc)= \delta $. 
		An identical statement -- with $\BD$ replaced by $\DiBD$ -- holds in the oriented case. 
	\end{lem}
	\begin{proof} If $ \alg$ factors through $T_\delta$, then clearly $\alg (\bigcirc) = \delta$. For the converse, let $\alg \colon \BD \to \modR$ be a symmetric monoidal functor such that $\alg (\bigcirc) = \delta$. Define a symmetric monoidal functor $ \underline \alg \colon \Br \to \modR$ by $\underline \alg (\tau') = \alg (\tau', 0)$ for each generator $\tau' \in \Br (m,n)$. Since $\alg$ is lax monoidal, for all morphisms $f = (\tau, \kcl) = (\tau, 0) \oplus (\emptyset, \kcl)\in \BD(m,n)$,
		\[\alg (f) = \delta^{\kcl}\alg (\tau) = \delta^{\kcl}\underline \alg (\tau) = \underline \alg (T_\delta (f)).\]
		Hence $ \alg = \underline \alg \circ  T_\delta \colon \BD \to \Br \to  \modR.$ 
		The proof is unchanged for the oriented case.
	\end{proof}

	For fixed $\delta \in R$ and $n \in \N$, the endomorphism algebras $\Br(n,n)$ coincide with \emph{Brauer algebras}, introduced by Brauer in \cite{Bra37} to study of representations of the finite dimensional orthogonal and symplectic groups $O_d$ and $Sp_k$ ($d, k \in \N$).
	
	Let $\Bbbk$ be a field of characteristic $0$ and let 
	$V$ be a $d$-dimensional vector space equipped with a nondegenerate bilinear form $\theta \colon V \otimes V \to V$ that is either symmetric or skew-symmetric (in which case, nondegeneracy implies that $d = 2k$ for some $k$).  Since $\theta$ is nondegenerate, it defines an isomorphism $v \mapsto \theta (v, -)$ of $V$ with its dual $V^*$.  Fix $\delta = d $ if $\theta$ is symmetric, and $\delta = -k = -d/2$ if $\theta$ is skew-symmetric. 
	
The isometry group $G = \{g \colon \theta (gv,gw) = \theta(v,w) \text{ for all } v,w, \in V \} \subset   GL(V)$  
 of $\theta$ is \begin{itemize}
	\item  the orthogonal group $ O(V, \theta) \cong O_d $ when $ \theta$ is symmetric
	\item the symplectic group $ Sp(V, \theta) \cong Sp_k$ when $ \theta$ is skew-symmetric.
\end{itemize}

	Brauer \cite{Bra37} extended the Schur-Weyl duality between representations of the symmetry and general linear groups to prove that, for $n \geq |\delta|$,  
	representations of $\Br^{\Bbbk}(n,n)$ in $V^{\otimes n}$ are in one-to-one correspondence with degree $n$ representations of $G$.

	Categorified versions of these results were established in \cite[Theorems~3.4,~4.6,~4.8,~5.9,~6.10]{LZ15}:

View the endomorphism prop $ T(V)$ (see \cref{ex. endo prop}) as a full sub-category of $\vect$ with objects $V^{\otimes k}, k \in \N$ (by convention, $V ^{\otimes 0} = \Bbbk$). Note that objects of $T(V)$ have a $G$-module structure induced by the factorwise action $ g\cdot (v_1, \dots, v_n) = (g (v_1), \dots, g(v_n))$ on each $V^{\otimes n}$. 

Let $ T_G(V) \subset   T(V)$ be the subprop of $G$-equivariant morphisms.
By definition $ \theta \colon V^{\otimes 2} \to \Bbbk$ is in $T_G(V)$ and hence, for all $n \in \N$ and $1 \leq i < j \leq n +2$, so are the \emph{``contraction''} maps 
$\theta^{i \ddagger j} \colon V^{\otimes (n + 2)} \to V^{\otimes n}$ induced by applying $\theta $ to the $i^{\mathrm{th}}$ and $j^{\mathrm{th}}$ factors. 

Recall that $ \Bbbk[\Sigma] \defeq \bigoplus_{n \in \N} \Bbbk[\Sigma_n]$ describes a monochrome $\vect$-prop. The canonical levelwise action of $\Sigma$ on $T(V)$ by permuting factors (see \cref{ex. endo prop}) extends linearly to a functor $ \Bbbk[\Sigma] \hookrightarrow T(V)$. 

For all $k \geq 0$, define \begin{equation}\label{eqn e(k)}
e(k) \defeq \sum_{\sigma \in \Sigma _k} \mathrm{sgn}(\sigma)\sigma \in \Bbbk[\Sigma_k] 
\end{equation} where $\mathrm{sgn}(\sigma)$ is the sign of a permutation $\sigma \in \Sigma _k$. Since $ \Bbbk[\Sigma]  \subset   \Br$ for all $ \delta \in \Bbbk$, for each $m,n \in \N$, we may define $\langle e(k) \rangle_{m,n} \subset   \Br(m,n)$ to be the subspace generated by $ e(k)$ under horizontal and vertical composition in $\Br$.

	\begin{thm}\label{thm lz15}
		[Lehrer-Zhang, 2015] 
		There is a unique symmetric strict monoidal (tensor) functor $ \Br\to \vect$ such that $1 \mapsto V$, $\cap \mapsto \theta$. This factors through the inclusion $T_G(V) \hookrightarrow \vect$. Let $F_G \colon \Br\to T_G(V) $ denote the corresponding (corestriction) functor.
		
		Let $\sigma_\two \in \Sigma_2 \subset  \Br(2,2)$ be the unique non-identity permutation. For all $v \otimes w \in V^{\otimes 2}$, \[F_G(\sigma_\two)(v\otimes w) = \left \{ \begin{array}{ll}
		w\otimes v& \text{ when } \theta \text { is symmetric,}\\
			-w\otimes v & \text{ when } \theta \text { is skew-symmetric.} 
		\end{array}\right . \]

The functor $F_G$ is full. Its restriction $\Br(m,n) \to T_G(V) (V^{\otimes m}, V^{\otimes n})$ is injective when $m+n \leq 2|\delta|$. When $ m +n > 2|\delta|$, 
its kernel is $\langle e(|\delta|+1) \rangle _{m,n}$.  
	\end{thm}

	\begin{rmk}
		\label{rmk fft and sft}
		The statement that $F_G \colon \Br\to T_G(V)$ is full is one formulation of the first fundamental theorem of invariant theory for the orthogonal and symplectic groups. In particular, it implies that, since $\Br(m,n) = 0$ when $ m +n $ is odd, so also $T_G(V)(V^{\otimes m}, V^{\otimes n}) = 0$ when $m +n$ is odd. 
		
	The second fundamental theorem 
	is given by the description of the kernels of the maps \newline$F_G(m,n) \colon \Br(m,n) \to T_G(V)(V^{\otimes m}, V^{\otimes n})$.
	\end{rmk}

Weyl's first and second fundamental theorems of invariant theory of the finite dimensional general linear groups are obtained from an oriented version of \cref{thm lz15}:

If $V$ is a finite $d$-dimensional vector space, then the general linear group $ GL = GL(V)$ (left) acts on $V$ by the standard representation $(g,v) \mapsto g(v)$, and (right) acts on $V^*$ by the dual representation $ (g, \alpha) \mapsto \left (v \mapsto \alpha (g^{-1} (v))\right ) $. As above, let $T_{GL}(V) \subset  T(V)$ be the subcategory of 
subcategory of $GL$-equivariant morphisms. In particular, the trace on $T(V)$ is $GL$-equivariant, as is the monoidal product of $GL$-equivariant morphisms in $T(V)$, so $T_{GL}(V)$ inherits a wheeled prop structure from $T(V)$.

For $k \in \N$, let $e (k) \in \Bbbk[\Sigma]$ be defined as above (\ref{eqn e(k)}) and let $\langle e(k) \rangle_{n,n}^{\Bbbk[\Sigma]} \subset   {\Bbbk[\Sigma_n]}$ be the subspace generated by $ e(k)$ under horizontal and vertical composition in $\Bbbk[\Sigma]$.

\begin{thm}
	[Weyl, \cite{Wey97}]\label{thm Weyl FFT SFT} 

The category $T_{GL}(V)$ of $GL$-equivariant morphisms in $T(V)$ is a $\vect$-groupoid such that $T_{GL} (V)(m,n) = 0$ when $m \neq n$. 
	
For $n \leq d$, $T_{GL}(V)(n,n)\cong \Bbbk[\Sigma_n]$ and for $n >d$,   $T_{GL}(V)(n,n)\cong \Bbbk[\Sigma_n]/ \langle e(d+1) \rangle^{\Bbbk[\Sigma]} _{n,n}$.
	
\end{thm}

By \cref{prop: wheeled prop oriented brauer}, this can be reformulated almost identically to \cref{thm lz15}:

\begin{cor}
	\label{thm GL OBD}
	There is a unique symmetric strict monoidal (tensor) functor $ \DiBr[d]\to \vect$ such that $(\uparrow) \ \mapsto V$, $(\downarrow) \ \mapsto V^*$ and $\cap \longmapsto \left((\alpha, v) \mapsto \alpha (v) \colon V^* \otimes V \to \Bbbk \right).$ This factors through the inclusion $T_{GL}^{\dipal}(V) \hookrightarrow \vect$.  
	
	 The corresponding (corestriction) functor $F_{GL} \colon \DiBr[d]\to T^{\dipal}(V)_{GL}$ is full. For $m,n \in \N$, its restriction  $\DiBr(m,n) \to T^{\dipal}(V)_{GL}(V^{\otimes m}, V^{\otimes n})$ is injective when $m+n \leq 2d$. When $ m +n > 2d$, 
	its kernel is $\langle e(d+1) \rangle _{m,n} \subset   \DiBr[d](m,n) $. 
	
\end{cor}

By \cite{DM23}, there is an equivalence of categories between algebras over $GL_d$ and wheeled props for which $\bigcirc = d$ and $e(d+1) = 0$. The comparison of Theorems \ref{thm Weyl FFT SFT} and \ref{thm lz15} is used in \cref{sec. invariants} to prove similar results -- in terms of unoriented circuit algebras -- for the categories of $O_d$ and $Sp_k$ algebras.

\begin{rmk}\label{rmk stable rep}
	Given a sequence of groups $(G_d)_d$ such that $G_d \hookrightarrow G_{d+1}$ for all $d \geq 0$, let $ G_\infty \defeq \bigcup G_d$ denote the colimit. A representation $W$ of $G_\infty$ is the colimit of a sequence of representations $(W_d)_{d }$ of the sequence of groups $(G_d)_d$ with inclusions $W_d \hookrightarrow W_{d+1}$ induced by the inclusions $G_d \hookrightarrow G_{d+1}$.
	
	For example, for all $d \geq 1$, the $d$-dimensional general linear group $ GL_d$ is naturally a subgroup of $GL_{d+1}$ under the inclusion induced by $\Bbbk^d \hookrightarrow \Bbbk^{d+1}  = \Bbbk^d \times \Bbbk$. The infinite general linear, orthogonal and symplectic groups $GL _\infty, O_\infty$ and $Sp_\infty$ are the colimits of the induced sequences $ (GL_d)_d,(O_d)_d$ and $(Sp_k)_k$. Let $\bm V \defeq \bigcup_{j = 0}^\infty \Bbbk^{\otimes j}$ be the standard representation.

The triangle identities (\ref{eq. n triangle}) in $\BD$ imply that, if $F \colon \BD \to \vect$ is a strict (or strong) symmetric monoidal functor with $F(1)  = V$, then $\theta = F(\cap)$ induces an isomorphism $\theta^* \colon V \xrightarrow{\cong} V^*$ and $V$ must be finite dimensional. Hence, there is no strict monoidal functor $ BD \to \vect$ such that $1 \mapsto {\bm V}$. 

However, if the sequence $(\theta_d  )_d$ of nondegenerate symmetric or skew-symmetric forms induces the sequence of orthogonal or symplectic groups $(G_d)_d$, there is a unique form ${\bm \theta} \defeq \mathrm{colim}_{d} \theta_d$ on $\bm V$ and a unique strict monoidal functor $ {\bm F}\colon\BDd \to \vect$, $1 \mapsto {\bm V}$ and $\cap \mapsto \bm \theta$, the image of which is the colimit of the (image of the) functors $F_{G_d}$ described in \cref{thm lz15}. 

Sam and Snowden \cite{SS15} established a contravariant equivalence between the categories of finite length functors $\BDd \to \vect$ (respectively $\mathsf{dOBD} \to \vect$) and representations of the infinite orthogonal and symplectic groups (respectively algebraic representations of the infinite general linear group). See also \cref{rmk nonunital CA} and \cref{ssec nonunital CA}, where a related result, \cref{thm. Ginfty}, is proved by extending the methods of \cite{DM23}.
\end{rmk}

\section{Wiring diagrams and circuit algebras}\label{sec. CA}

A circuit algebra is a given by a family of objects, indexed by some free commutative monoid (see \cref{ssec. wd and ca}), with operations that are governed by wiring diagrams. These are, essentially, non-planar versions of Jones's planar diagrams \cite{Jon99}. Wiring diagrams are commonly described by partitioning boundaries of 1-manifolds (e.g.,~{\cite{BND17, DHR20,DHR21}}). However, they 
admit a straightforward definition in terms of 
Brauer diagrams. This paper takes the latter approach.

\subsection{Operadic preliminaries} \label{ssec operads}

This section summarises the basic theory of (coloured) operads. See \cite{Lei04} and \cite{BM07} for more details. 

 A (symmetric) $\DDD$-coloured operad $\op$ (in the category of sets) is given by a $(\listm (\DDD) \times \DDD)$-graded set $ (\op(\ccc; d))_{(\ccc; d)}$, and a family of \textit{composition morphisms},
\[ \gamma \colon \op(\ccc; d) \times \left( \prod_{ i = 1}^m \op(\bbb_i; c_i)\right) \to \op(\bbb_1 \dots \bbb_m;d),\] 
 defined for each $d \in \DDD$, $ \ccc  = (c_i)_{i = 1}^m \in \listm(\DDD)$ and  $\bbb_i  \in \listm (\DDD)$, for $1 \leq i \leq m$.
 
If $\phi \in \op(c_1, \dots, c_m; d)$, then $d$ is called the \textit{output} of $\phi$ and each $c_i$ is an \textit{input} of $\phi$. The symmetric groupoid $\Sigma$ acts on $\op$ by permuting the inputs: each $\sigma \in \Sigma_m$ induces isomorphisms $\op(c_{\sigma 1}, \dots, c_{\sigma m}; d)  \xrightarrow{\cong}\op(c_1, \dots, c_m; d).$
The composition $\gamma$ is required to be associative and equivariant with respect to the $\Sigma$-action on $\op$. 

Moreover, for all $d \in \DDD$, there is an element $\nu_d \in \op (d;d)$ that acts as a 2-sided unit for $\gamma$:  
for all $\ccc  = (c_1, \dots, c_m) \in \listm (\DDD)$, the composite morphisms 
\[\op (\ccc;d)\xrightarrow{ \cong } I \times \op (\ccc;d)\xrightarrow{(\nu_d, id)}  \op(d; d) \times \op (\ccc;d) \xrightarrow{\gamma} \op (\ccc;d), \]

\[
\op (\ccc;d)\xrightarrow{ \cong }\op (\ccc;d)\times I \\ \xrightarrow{(id, \bigotimes_{i =1}^m \nu_{c_i})}  \op (\ccc;d)\times \left(\bigotimes_{ i =1}^m \op (c_i; c_i) \right)\xrightarrow{\gamma}   \op (\ccc;d)\] are the identity on $\op (\ccc;d)$.  

	Let $ (\CCat, \oplus, 0)$ be a small permutative category with object set $\CCat_0$.
	
	\begin{defn}
		\label{def. operad of a perm cat}
		The $\CCat_0$-coloured \emph{operad $\op^{\CCat}$ underlying $(\CCat, \oplus, 0)$} is 
		defined by 
		\[\op^{\CCat}(x_1, \dots, x_n; y)  \defeq \CCat(x_1 \oplus \dots \oplus x_n, y ),\] 
		with operadic composition $\gamma $ in $\op^{\CCat}$ induced by composition in $\CCat$ as follows:
		
		Let the operation $\overline g \in \op^{\CCat}(x_1, \dots, x_n;y)$ correspond to the morphism $g \in \CCat (x_1 \oplus \dots \oplus x_n, y)$ and, for $1 \leq i \leq n$, let $\overline f_i \in \op^{\CCat}(w_{i,1}, \dots, w_{i, m_i};x_i)$ correspond to $f_i \in \CCat (w_{i,1}\oplus \dots \oplus w_{i,m_i}, x_i)$. 
Then,
		\[ \gamma \left(\overline g, (\overline f_i)_i\right)  \defeq  \overline{\left(g \circ  (f_1 \oplus \dots \oplus f_n)\right)}. \]

	\end{defn}
(In fact, any small cocomplete symmetric monoidal category $\V$ has an underlying operad by defining, for $x_1, \dots, x_n \in \V$, the object $x_1 \otimes \dots \otimes x_n \in \V$  as in \cref{ex. endo operad}.) 
	
	Observe that, if $\op^{\CCat}$ is the operad underlying a small permutative category $\CCat$, then, for all $\overline f_1 \in \op^{\CCat}(x_{1,1},\dots, x_{1,m};y_1)$ and $\overline f_2  \in \op^{\CCat}(x_{2,1},\dots, x_{2,n};y_2)$, there is an operation
	\begin{equation}
		\label{eq. op monoid}
		\overline f_1 \oplus \overline f_2 \defeq \gamma \left (\overline {id}_{y_1 \oplus y_2}, (\overline f_1, \overline f_2) \right) \in 
		\op^{\CCat}(x_{1,1},\dots, x_{1,m},x_{2,1},\dots, x_{2,n};y_1\oplus y_2).
	\end{equation}
	
By definition, $  \op^{\CCat}(-;y) \cong \op^{\CCat}(0;y)$ canonically for all $y$. In particular, there is a canonical isomorphism $ \op^{\CCat}(-;0) \xrightarrow \cong \op^{\CCat}(0;0) = \CCat(0,0)$. Let $\overline{ id_0} \in \op^{\CCat}(-;0)$  be the preimage of $id_0 \in \CCat(0,0)$ under this isomorphism. Then, for all $ (x_1, \dots, x_k, y)$, precomposition with $(\bigotimes_{i = 1}^k id_{x_k}, \overline{ id_0})$ induces an isomorphism $   \op^{\CCat}(x_1, \dots, x_k, 0; y) \xrightarrow \cong \op^{\CCat}(x_1, \dots, x_k; y)$.

	\medspace

For $ i \in \{1,2\}$, let $(\op^i, \gamma^i, \nu^i)$ be a $\DDD_i$-coloured operad. A morphism $\mathcal F \colon (\op^1, \gamma^1, \nu^1) \to  (\op^2, \gamma^2, \nu^2)$ of (coloured) operads is given by 
a map of sets $f \colon \DDD_1 \to \DDD_2$, and a $(\listm (\DDD)_1 \times \DDD_1)$-indexed family of maps 
\[\mathcal F_{(c_1, \dots, c_k; d)} \colon \op^1(c_1, \dots, c_k; d) \to \op^2(f (c_1), \dots,f (c_k); f(d) )\] that respect units and composition, and are equivariant with respect to the symmetric action.

If $f = id_{\DDD}$ (with $\DDD = \DDD_1= \DDD_2 $), then $\mathcal F \colon \op^1 \to \op^2$ is called \textit{colour-preserving}. The category of $\DDD$-coloured operads and colour-preserving morphisms is denoted by $\mathsf{Op}^{\DDD}$.

In the remains of this section, $(\V, \otimes, I)$ is a symmetric monoidal category with all finite colimits, and $(\CCat, \oplus, 0)$ is a small permutative category with object set $\CCat_0$.

\begin{ex}
	\label{ex. endo operad}
	For any $n$-tuple $(x_1, \dots, x_n)$ of objects in $\V$, define $x_1 \otimes \dots \otimes x_n$ to be the colimit, under associator isomorphisms in $\V$, of all ways (indexed by planar binary rooted trees) of tensoring $x_1, \dots, x_n$. Given a set $\DDD$ and a $\DDD$-indexed object $A = (A_c)_{c \in \DDD}$ in $\V$, the $\DDD$-coloured \textit{endomorphism operad $End^A$} is defined by
	\[End^A(c_1, \dots, c_k;d)\defeq \V\left(A_{c_1} \otimes \dots \otimes A_{c_k} , A_{d} \right),\] together with the obvious composition and units induced by composition and identities in $\V$.
	
\end{ex}

\begin{defn}\label{def operad alg}
	A \emph{$\V$-algebra for a $\DDD$-coloured operad $\op$} is a $\DDD$-indexed object $(A_c)_{c \in \DDD}$ in $\V$, together with a morphism $\alg \colon \op \to End^A$ of $\DDD$-coloured operads. 
	
	The category $\mathsf{Alg}_{\V}(\op)$ of $\V$-algebras for $\op$ is the subcategory of the slice category $\op \ov \mathsf{Op}^{\DDD}$ whose objects are $\V$-algebras for $\op$. 
	Morphisms in $\mathsf{Alg}_{\V}(\op) \left((A, \alg) , (B, \mathcal B)\right)$ are of the form $(g, (g_c)_c)$ where $ g \colon \alg \to \mathcal B$ in $\op \ov \mathsf{Op}^{\DDD}$ and, for all $c \in \DDD$, $ g_c \in \V( A_c , B_c)$ such that, if $\phi \in \op (c_1, \dots, c_k;d)$, then the following diagram commutes in $\V$:
	
	\[ \xymatrix{ A_{c_1}\otimes \dots \otimes A_{c_k} \ar[d]_{\alg(\phi)}\ar[rr]^-{g_{c_1} \otimes \dots \otimes g_{c_k}}&& 
		B_{c_1} \otimes \dots \otimes B_{c_k} \ar[d]^{g\alg(\phi)}\\ 
		A_{d}\ar[rr]_-{g_d} && B_{d}.}\]
\end{defn}

\begin{rmk}
Observe that 
\cref{def operad alg}, though it relies on the symmetric monoidal structure on $\V$, is concerned with operads in the category of sets and does not involve operads enriched in a (closed) symmetric monoidal category. It therefore diverges slightly from the usual definition of an operad algebra (as in \cite{BM07}).
\end{rmk}

	Let $(\CCat, \oplus, 0)$ be a small permutative category and $ (\V, \otimes, I)$ a cocomplete symmetric monoidal category 
	 and let $ [\CCat, \V]_{\mathrm{lax}}$ denote the category of symmetric monoidal functors $\alg \colon(\CCat, \oplus, 0) \to (\V, \otimes, I)$.

\begin{lem}\label{rmk. lax functors and algebras}\label{prop. lax functors and algebras}
The categories $\mathsf{Alg}_{\V} (\op^{\CCat})$ of $\V$-algebras for the operad $\op^{\CCat}$ underlying 
$(\CCat, \oplus, 0)$ and $[\CCat, \V]_{\mathrm{lax}}$ are canonically isomorphic. 

\end{lem}
\begin{proof}

If $\op = \op^{\CCat}$ is the $\CCat_0$-coloured operad underlying $\CCat$, and $(\alg, \pi, \eta)\colon \CCat \to \V$ is a symmetric monoidal functor, then $ (\alg(x))_{x \in \CCat_0}$ has an $\op$-algebra structure as follows: For $k \geq 1$ and all $\overline f \in \op(x_1, \dots, x_k;y)$ induced by $f \in \CCat (x_1\dots x_k,y)$,  \[\alg_{x_1, \dots, x_k;y} (\overline f ) = \alg(f) \circ \pi_{x_1, \dots, x_n} \in \V (\alg(x_1)\otimes \dots \otimes \alg(x_k), \alg(y)).\] 
(Here $\pi_{x_1, \dots, x_n} \colon \alg(x_1) \otimes \dots \otimes \alg(x_n) \to \alg(x_1 \dots x_n)$ is the universal map  from the colimit.)

When $k = 0$, and $\overline f \in \op (-;y)$ is induced by $ f \in \CCat (0,y)$, 
\[\alg_{-;y} (\overline f ) \defeq \alg(f) \circ\eta \in \V (I, \alg(y)).\]

Conversely, a $\V$-algebra $(A, \hat \alg)$ for $\op$ induces a functor $\alg \colon \CCat \to \V$ 
described by $x \mapsto  A_x$ for all $x \in \CCat$. If $\overline f \in \op (x;y)$ is induced by $f \in \CCat (x,y)$, then  $  f \mapsto \hat \alg (\overline f) \in \V( A_x, A_y)$.
This has symmetric lax monoidal structure 
$\pi_{\alg} \colon A_x \otimes A_y \to A_{x\oplus y}$ induced by $ \overline {id_{x \oplus y}} \in \op (x,y;x\oplus y)$ and $\eta_{\alg} \colon I \to A_0$ induced by $\overline {id_0}\in \op (-;0)$. 
It follows from the definitions that the assignments $(\alg, \pi, \eta) \mapsto ((\alg(x))_x, \alg)$ and $(A, \hat \alg) \mapsto (\alg, \pi_{\alg}, \eta_{\alg})$ extend to mutually inverse functors $\mathsf{Alg}_{\V} (\op^{\CCat}) \leftrightarrows [\CCat, \V]_{\mathrm{lax}}$.

(For more details, see e.g.,~ \cite[Chapters~2-3,]{Lei04}: Example~2.1.10 and Section~3.3 in particular.)
\end{proof}

\begin{defn}
	\label{defn operad alg ideal}
	Let $(A , \alg)$ be an algebra over a $ \DDD$-coloured operad $\op$. An \emph{ideal of $(A, \alg)$} is an $\op$-subalgebra $ (I, \mathcal I) \subset   (A , \alg)$ such that, for all $n\in \N$, $(c_1, \dots, c_n ) \in \DDD^{n}$, $d \in \DDD$ and $x_i \in A_{c_i}$, and all 
	$\phi \in \op (c_1, \dots, c_n;d)$, if $x_j \in I_{c_j}$ for some $1 \leq j \leq n$, then $\alg (\phi)(x_1 , \dots  , x_n) \in \mathcal I(d)$.
\end{defn}

Equivalently, $(I, \mathcal I) \subset   (A , \alg)$ is an ideal precisely if the quotient $ (A/I, \alg / \mathcal I)$ inherits an $\op$-algebra structure from $(A, \alg)$.

	If $\op^{\CCat}$ is the operad underlying a monoidal category $ \CCat$, and $ (\alg, \pi,\eta )\colon \CCat \to \V$ is a symmetric monoidal functor as in \cref{prop. lax functors and algebras}, then an ideal $(I , \mathcal I)$ of the operad algebra corresponding to $\alg$ is a symmetric monoidal subfunctor $  \mathcal I \hookrightarrow \alg$ such that, for all $x, y \in \CCat$, the restrictions of  $\pi_{x,y} \colon \alg(x) \otimes \alg (y) \to \alg (x \oplus y)$ to $ \mathcal I (x) \otimes \alg (y) $ and $ \alg (x ) \otimes \mathcal I (y)$ describe morphisms to $\mathcal I (x \oplus y)$.

\subsection{Wiring diagrams and circuit algebras}\label{ssec. wd and ca}

As in \cite{BND17, DHR20, DHR21}, circuit algebras will be defined as algebras over an operad of \emph{wiring diagrams}.

\begin{defn}\label{defn. operad of wiring diagrams}\label{def. oriented WD}

For a given palette $(\CCC, \omega)$, and each  $(\ccc_1, \dots, \ccc_k;\ddd) \in \listm^2 (\CCC) \times \listm(\CCC)$, a \emph{wiring diagram of type $(\ccc_1, \dots, \ccc_k; \ddd)$} is an element of the set
\[ \CWD (\ccc_1, \dots, \ccc_k; \ddd) \defeq \CBD(\ccc_1 \oplus \dots \oplus \ccc_k; \ddd).\]
The \emph{$\listm (\CCC)$-coloured operad of $(\CCC, \omega)$-wiring diagrams} is the {operad $\CWD \defeq \op^{{\CBD}}$ underlying $\CBD$}.

In particular, for $ (\CCC, \omega) = \dipal$, $\oWD \defeq \CWD$ is the operad of \emph{(monochrome) oriented wiring diagrams}, and for a set $\DDD$, the \emph{operad $\oCWD$ of $\DDD$-coloured oriented wiring diagrams} is the operad underlying the category $\DiCBD = \BD^{\DDD\times \dipal}$ of $\DDD$-coloured oriented Brauer diagrams. 

When $\CCC$ is the singleton set, the $\N $-coloured \emph{operad of (monochrome) wiring diagrams} $\CWD  \defeq \op^{\BD}$ is denoted by $ \WD$.

\end{defn}

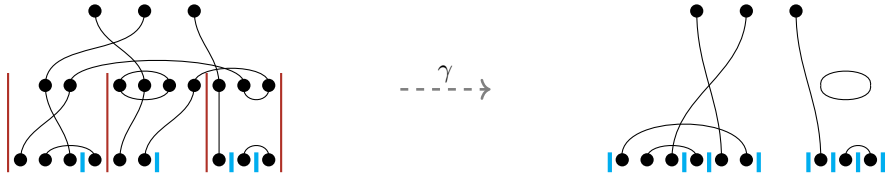
\begin{figure}
[htb!]
\begin{tikzpicture}
	\node at (0,0){
\begin{tikzpicture}[scale = .33]
	\begin{pgfonlayer}{above}
		\draw[thick,red] (.5,-.5)--(.5,3.5);
		\node [dot] (1) at (1,0) {};
		\node [dot] (2) at (2,0) {};
		\node [dot] (3) at (3,0) {};
		\draw[ultra thick,cyan] (3.5,-.5)--(3.5,.3);
		\node [dot] (4) at (4,0) {};
		\draw[thick,red] (4.5,-.5)--(4.5,3.5);
		\node [dot] (5) at (5,0) {};
		\node [dot] (6) at (6,0) {};
		\draw[ultra thick,cyan] (6.5,-.5)--(6.5,.3);
		\draw[thick,red] (8.5,-.5)--(8.5,3.5);
		\node [dot] (9) at (9,0) {};
		\draw[ultra thick,cyan] (9.5,-.5)--(9.5,.3);
			\node [dot] (10) at (10,0) {};
		\draw[ultra thick,cyan] (10.5,-.5)--(10.5,.3);
			\node [dot] (11) at (11,0) {};
		\draw[thick,red] (11.5,-.5)--(11.5,3.5);
		
		\node [dot] (12) at (2,3) {};
			\node [dot] (13) at (3,3) {};
			\node [dot] (15) at (5,3) {};
			\node [dot] (16) at (6,3) {};
			\node [dot] (17) at (7,3) {};
			\node [dot] (18) at (8,3) {};
				\node [dot] (19) at (9,3) {};
			\node [dot] (20) at (10,3) {};
				\node [dot] (21) at (11,3) {};
				
				\node [dot] (24) at (4,6) {};
					\node [dot] (26) at (6,6) {};
					\node [dot] (28) at (8,6) {};
		\end{pgfonlayer}
	\begin{pgfonlayer}{background}
			\draw [in=-90, out=90] (1.center) to (13.center);
				\draw [in=-90, out=90] (3.center) to (12.center);
					\draw [in=-90, out=90] (5.center) to (16.center);
						\draw [in=-90, out=90] (6.center) to (18.center);
							\draw [in=-90, out=90] (9.center) to (19.center);
								\draw [in=-90, out=90] (12.center) to (26.center);
									\draw [in=-90, out=90] (16.center) to (24.center);
						\draw [bend left=100, looseness=1] (2.center) to (4.center);
							\draw [bend left=100, looseness=2] (10.center) to (11.center);
						\draw [in=-90, out=90] (19.center) to (28.center);
							\draw [bend right=100, looseness=1] (15.center) to (17.center);
								\draw [bend right=100, looseness=2] (20.center) to (21.center);
							\draw [bend left=100, looseness=1] (15.center) to (17.center);
								\draw [bend left=100, looseness=.5] (13.center) to (20.center);
								\draw [bend left=100, looseness=.8] (18.center) to (21.center);
		\end{pgfonlayer}
\end{tikzpicture}};
\draw[gray, dashed, ->, line width = 1](3.4,0)--(4.6,0);
\node at (4,.2){$\gamma$};
	\node at (8,0){
	\begin{tikzpicture}[scale = .33]
		\begin{pgfonlayer}{above}
			\draw[ultra thick,cyan] (.5,-.5)--(.5,.3);
			\node [dot] (1) at (1,0) {};
			\node [dot] (2) at (2,0) {};
			\node [dot] (3) at (3,0) {};
			\draw[ultra thick,cyan] (3.5,-.5)--(3.5,.3);
			\node [dot] (4) at (4,0) {};
			\draw[ultra thick,cyan] (4.5,-.5)--(4.5,.3);
			\node [dot] (5) at (5,0) {};
			\node [dot] (6) at (6,0) {};
			\draw[ultra thick,cyan] (6.5,-.5)--(6.5,.3);
			\draw[ultra thick,cyan] (8.5,-.5)--(8.5,.3);
			\node [dot] (9) at (9,0) {};
			\draw[ultra thick,cyan] (9.5,-.5)--(9.5,.3);
			\draw[ultra thick,cyan] (10.5,-.5)--(10.5,.3);
			\draw[ultra thick,cyan](11.5,-.5)--(11.5,.3);
				\node [dot] (10) at (10,0) {};
					\node [dot] (11) at (11,0) {};
					
			\node [] (12) at (2,3) {};
			\node [] (13) at (3,3) {};
			\node [] (15) at (5,3) {};
			\node [] (16) at (6,3) {};
			\node [] (17) at (7,3) {};
			\node [] (18) at (8,3) {};
			\node [] (19) at (9,3) {};
			\node [] (20) at (10,3) {};
			\node [] (21) at (11,3) {};
			
			\node [dot] (24) at (4,6) {};
		\node [dot] (26) at (6,6) {};
		\node [dot] (28) at (8,6) {};
	\end{pgfonlayer}
\begin{pgfonlayer}{background}
\draw [in=-90, out=90] (3.center) to (26.center);
\draw [in=-90, out=90] (5.center) to (24.center);
\draw [in=-90, out=90] (9.center) to (28.center);
\draw [bend left=100, looseness=1] (2.center) to (4.center);
\draw [bend left=100, looseness=1] (1.center) to (6.center);
	\draw [bend left=100, looseness=2] (10.center) to (11.center);
\draw [bend right=100, looseness=1] (19.center) to (21.center);
\draw [bend left=100, looseness=1] (19.center) to (21.center);
\end{pgfonlayer}
\end{tikzpicture}};
\end{tikzpicture}

\caption{Composition in $\WD$. (See also \cref{fig. pictorial}.) }
\label{fig. CWD comp}
\end{figure}

\begin{defn}\label{def CA}A \emph{$(\CCC,\omega)$-coloured $\V$-circuit algebra} is a $\V$-valued algebra for the operad $\CWD$ of $(\CCC, \omega)$-coloured wiring diagrams. 
	The full subcategory of $\V$-circuit algebras in $\mathsf{Alg}(\CWD)$ is denoted by $\V\mathdash\CCA$. When $\V = \Set$, $\V\mathdash\CCA$ is denoted simply by $\CCA$. 
	
	If $(\CCC, \omega) = \{*\}$ is trivial, then $\V\mathdash\CA \defeq \V\mathdash\CCA[*]$ is the category of monochrome $\V$-circuit algebras.

	\emph{Oriented} (respectively \emph{non-oriented}) \emph{circuit algebras} are algebras over operads of oriented (respectively non-oriented) wiring diagrams.

	\end{defn}
	
	\begin{rmk}\label{rmk classic wd}
		Though \cref{defn. operad of wiring diagrams} is already observed in \cite[Definition~2.9]{BND17}, wiring diagrams are commonly described (for example in {\cite{DHR20,DHR21}}) as isotopy classes of immersions of compact 1-manifolds in punctured 2-discs that are injective on boundaries and preserve boundaries and interiors.
		
	 In this representation, composition is defined by inserting discs into the punctures in such a way that the boundaries agree. \cref{fig. disc rep} provides a punctured disc representation of the same composition of wiring diagrams as \cref{fig. CWD comp}. In the coloured case, 1-manifolds are coloured according to \cref{ex. manifold components} and \cref{def. colouring} to define $(\CCC, \omega)$-coloured wiring diagrams. For the operadic composition in $\CWD$ the colours on the disc boundaries are required to match. 
		
		The punctured disc representation of wiring diagrams provides a clear visualisation of the relationship of wiring diagrams (and hence circuit algebras) to planar diagrams and algebras \cite{Jon99} and tangle categories \cite{Tur89}. It also exhibits the operad of monochrome wiring diagrams as a suboperad of the Spivak's \emph{operad of wiring diagrams} \cite{Spi13}. Moreover, the disc representation of wiring diagrams is highly suggestive of the relationship between circuit algebras and modular operads (c.f.,~\cref{ssec ca mo}), and the graphical construction of circuit algebras that is developed in the sister paper \cite{RayCA2}. 
		On the other hand, the definition in terms of Brauer diagrams is obviously combinatorial and reveals connections between circuit algebras and representations of classical groups (c.f., Sections \ref{ssec. representations of BD}, \ref{sec. invariants}).

	\end{rmk}
	
		\begin{figure}[htb!]
		\includegraphics[width=0.95\textwidth]{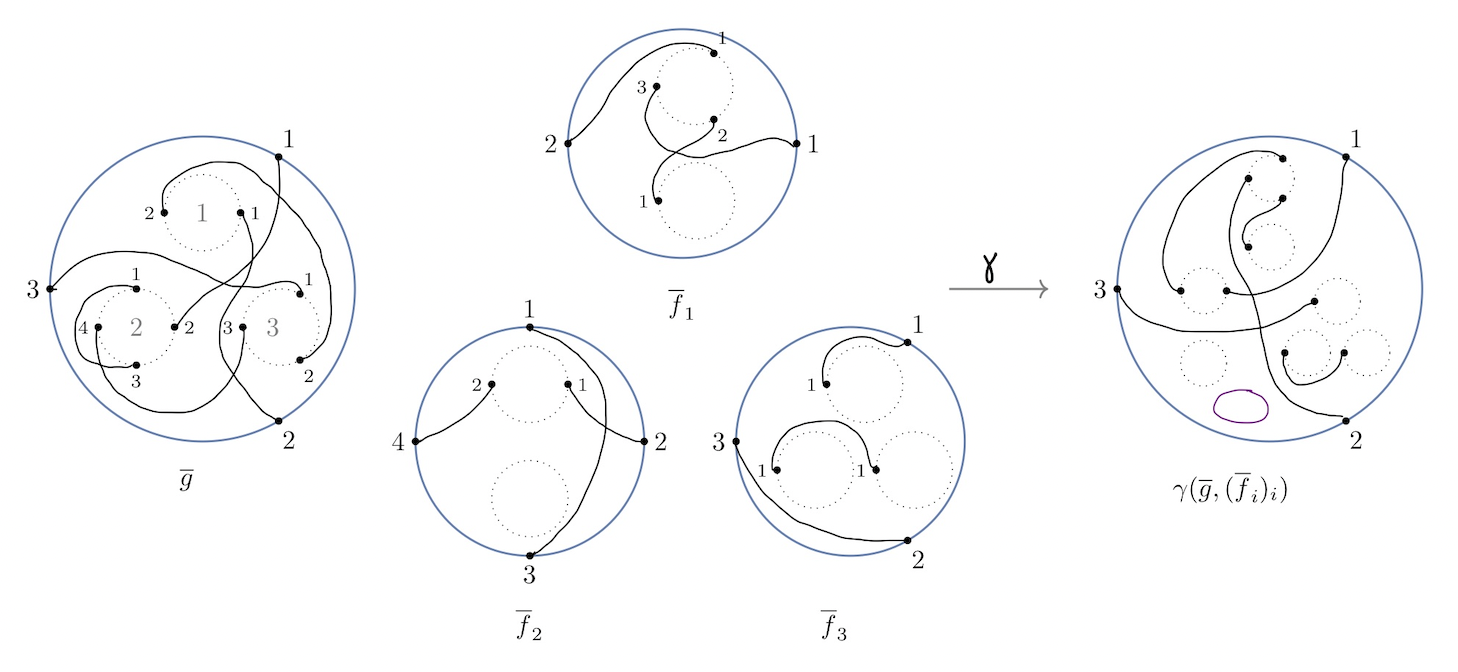}

		\caption{Disc representation of the wiring diagram composition in \cref{fig. CWD comp}.}
		\label{fig. disc rep}\label{fig. pictorial}
	\end{figure}

Let $(\V, \otimes, I)$ be a cocomplete symmetric monoidal category. As an algebra for the operad $\CWD$, a $(\CCC, \omega)$-coloured $\V$-circuit algebra consists of objects $(\alg(\ccc))_{\ccc \in \listm (\CCC)}$ and, for each $(\ccc_1, \dots, \ccc_k; \ddd)  \in \listm^2 \CCC \times \listm (\CCC)$, a set of $\V$-morphisms $\alg{(\overline f, \lambda)}\colon \bigotimes _{ i = 1} ^k \alg (\ccc_i) \to \alg(\ddd)$ indexed by Brauer diagrams $(f, \lambda) \in \CBD(\ccc_1 \oplus \dots \oplus \ccc_k , \ddd)$. These satisfy:
\begin{itemize}
\item for all $\ccc \in \listm (\CCC)$, $\alg({ \overline {id}_{\ccc}} )= id_{\alg(\ccc)} \in \V(\alg(\ccc),\alg(\ccc))$; 
\item the morphisms $\alg{( \overline f, \lambda)}$ are equivariant with respect to the $\Sigma$-action on $\listm (\CCC)$ and on $\CWD$;
\item given wiring diagrams $(\overline f, \lambda) \in \CWD (\ccc_1, \dots, \ccc_k; \ddd)$, and, for all $1 \leq i \leq k $, $(\overline f^i, \lambda^i)\in \CWD (\bbb_{i, 1}, \dots, \bbb_{i, k_i};\ccc_i)$, the following diagram commutes in $\V$:
\begin{equation}\label{eq. algebra maps} \xymatrix{ \bigotimes_{i = 1}^k \bigotimes_{j = 1}^{k_i} \alg(\bbb_{i,j}) \ar[rrrrd]_-{ \alg{ \gamma\left(( \overline f, \lambda),( \overline f^i, \lambda^i)_{ i}\right)}}\ar[rrrr]^{\bigotimes_{i = 1}^k \alg{(\overline f^i, \lambda^i)}} &&&& \bigotimes_{i = 1}^k \alg(\ccc_i) \ar[d]^{ \alg{(\overline f, \lambda )}}\\ &&&& \alg(\ddd)}\end{equation}

\end{itemize}

The following is immediate from \cref{rmk. lax functors and algebras}:
\begin{thm}\label{thm. lax functor ca}
The category $\V\mathdash\CCA$ of $(\CCC, \omega)$-coloured $\V$-circuit algebras is isomorphic to the category of symmetric monoidal functors $ \CBD \to \V$. (See e.g.,~\cite[Example~2.1.10~\&~Section~3.3]{Lei04} for more details.)

\end{thm}

\begin{ex}		\label{ex. wheeled props}

	Oriented circuit algebras are described in detail in \cite{DHR20, DHR21}. \cref{prop: wheeled prop oriented brauer} and \cref{thm. lax functor ca} provide another proof of the result, established in \cite{DHR20}, that $\DDD$-coloured oriented circuit algebras are equivalent to $\DDD$-coloured wheeled props. 

\end{ex}

\begin{ex}\label{ex initial CA} 
	
	Let $\inica \colon \BD \to \Set$ be the circuit algebra defined by $\inica(n )\defeq \BD(0,n)$ and for all $ g \in \BD (m,n)$, $\inica(g)(f)=g \circ f$. This is initial in the category of (monochrome) circuit algebras in $\Set$: For any such  $(\alg, \pi, \eta)$, there is a unique morphism $\uniA \colon \inica \to \alg$ such that $\uniA(f) =( \alg (f) \circ \eta)(1) \in \alg(n)$ for all $ f \in \BD (0,n)$.

	As in \cref{ssec. representations of BD}, for a fixed commutative ring $R$, let $\underline {\BD} $ be the free $\modR$-category on $\BD$. Let $\kinica = \kinica[R]  $ be the free $\modR$-circuit algebra on $\inica$, defined by $ \kinica(n)\defeq \underline {\BD} (0,n)$. This is initial in the category of (monochrome) $\modR$-circuit algebras.

For a palette $(\CCC, \omega)$, the initial $(\CCC,\omega)$-coloured circuit algebra $\Cinica$ (and $(\modR)$-circuit algebra $\Ckinica$) with $ \Cinica(\ccc) \defeq \CBD(\varnothing, \ccc)$  may be similarly defined.
In particular, by \cref{prop: wheeled prop oriented brauer}, $ \Ckinica[\dipal]$ describes the initial monochrome $\modR$-wheeled prop $U$ (called $\mathcal Z$ in \cite{DM23}) with $U(m,n)  =\underline{\BD}^{\dipal}(\varnothing, \uparrow^m \downarrow^n)$. 

\end{ex}

In the following examples, $\vect$ is always the category of vector spaces over a field $\Bbbk$ of characteristic 0 and $V$ is a (finite) $d$-dimensional vector space that generates the full subcategory $T(V)\subset \vect$ on objects of the form $V^{\otimes n}$, $n \in \N$.  
\begin{ex}\label{ex lz ca}

If $\theta$ is a symmetric or skew-symmetric nondegenerate bilinear form on $V$ with isometry group $G$, then let 
\[ \delta  = \left\{ \begin{array}{ll}
	 d & \text{ when } \theta \text{ is symmetric, in which case } G \cong O_{\delta}\\
	-\frac d2 & \text{ when } \theta \text{ is skew-symmetric, , in which case } G \cong Sp_{|\delta|}.
	\end{array}\right. \]
	
By \cref{thm lz15}, there is a unique symmetric strict monoidal functor $ \vGca  \colon \BD\to \vect$ such that $1 \mapsto V$ and $\cap \mapsto \theta$, and this factors through the symmetric strict monoidal functor $F_G \colon \Br\to T_G(V)$ where 
$T_G(V) \subset  T(V)$ is the subprop of $G$-equivariant morphisms and, as in \cref{ssec. representations of BD},  $\Br$ is the Brauer category with specialisation $\delta \in \Bbbk$. \cref{thm lz15} implies, moreover, that the kernel of the unique $\vect$-circuit algebra morphism $\uniG \colon \kinica \to \vGca$ is the circuit algebra ideal  $\mathcal I_{\theta}\subset  \kinica$ generated by $\bigcirc - \delta \in \kinica(0)$ and $ \coev{e(|\delta|+1)}\in \kinica(2(|\delta| + 1))$, where \begin{equation}
	\label{eq coev e(k)}
	\coev{e(k)} \defeq \sum_{\sigma \in \Sigma_k} \mathrm{sgn}(\sigma)\coev{\sigma} 
\end{equation}is the element of $\kinica(2k)$ obtained by linear coevaluation of the components of $e(k) \in \underline{\BD}(k,k)$ (\ref{eqn e(k)}).

In particular, if $\vGca^G\subset  \vGca$ is the $G$-invariant sub-circuit algebra, then $\kinica / \mathcal I_{\theta} \cong \vGca^G$.

In \cref{sec. invariants}, it is proved that there is an equivalence between algebras over the orthogonal (and symplectic) groups and circuit algebras $\alg$ such that  $\coev{e(|\delta|+1)}$ and  $\bigcirc - \delta $ are in the kernel of the unique morphism $ \uniA  \colon \kinica \to \alg$.

\end{ex}

\begin{ex}
	\label{ex weyl WP}

	As in \cref{ssec. representations of BD}, $T_{GL}(V)\subset  T(V)$ is the sub-wheeled prop of $GL(V)$-equivariant morphisms. 
	By \cref{thm Weyl FFT SFT}, the kernel of the unique morphism $a_V \colon U \to T_{GL}(V)$ of wheeled props is generated by $e(d+1) \in U(d+1,d+1)$ and $ \bigcirc - d \in U(0,0)$. Equivalently, the kernel of the unique $\vect$-valued oriented circuit algebra morphism $\Ckinica[\dipal]\to \vca$ is generated by $\coev{e(d+1)}$ and $\bigcirc -d$.

\end{ex}

\begin{defn}\label{ex. free ca}
	
Given any $\listm (\CCC)$-graded set $S = (S_{\ccc})_{\ccc \in \listm (\CCC)}$, the \textit{free ($\Set$-valued) circuit algebra $\mathcal{F}^{(\CCC, \omega)}\langle S \rangle$ on $S$} is defined as follows:
	
	The collection $(F^{(\CCC, \omega)}\langle S \rangle_{\ddd})_{\ddd \in \listm (\CCC)}$ of $(\CCC,\omega)$-coloured wiring diagrams \textit{decorated by $S$} is defined by
	\[ \begin{array}{ll}
		F^{(\CCC, \omega)}\langle S \rangle_{\ddd}& = \coprod_{(\ccc_1, \dots, \ccc_k) \in \listm^2 \CCC}\left( \CWD (\ccc_1, \dots, \ccc_k; \ddd) \times \prod_{i = 1}^k S(\ccc_i)\right)\\
		& = \coprod_{((\ccc\oplus\dots\oplus\ccc_k),(f, \lambda)))\in \CBD \ov \ddd} \left ( \prod_{i = 1}^k S(\ccc_i)\right).
	\end{array}\]

	For each $(\overline f, \lambda )\in \CWD(\ccc_1, \dots, \ccc_k; \ddd)$, the morphism $\mathcal F^{(\CCC, \omega)} \langle S \rangle(\overline f, \lambda) \colon F^{(\CCC, \omega)}\langle S \rangle_{\ccc_1}\otimes \dots \otimes  F^{(\CCC, \omega)}\langle S \rangle_{\ccc_k} \to F^{(\CCC, \omega)}\langle S \rangle_{\ddd}$ is described by 
	\[ \prod_{i = 1}^k \left( (\overline f^i, \lambda^i ), (x_{j_i}^i)_{j _i= 1}^{m_i}\right) \mapsto \left(\gamma\left((\overline f, \lambda) , \left((\overline f^i, \lambda^i )_{i = 1}^k\right)\right), (x_{j_i}^i)_{{\overset{1 \leq j_i \leq m_i}{1 \leq i \leq k}}} \right).\]

	For a fixed commutative ring $R$, let $\Ckinica\langle S \rangle$ be the $\modR$-circuit algebra freely generated by $\mathcal F^{(\CCC, \omega)}\langle S \rangle$. So, for all $\ccc$, 
	$\Ckinica\langle S \rangle(\ccc)$ is the free $R$-module on $F^{(\CCC, \omega)}\langle S \rangle_{\ccc}$. 
	
	When $(\CCC, \omega) = \{*\}$ is trivial, write $\mathcal F\langle S \rangle  = \mathcal F^{(\CCC, \omega)}\langle S \rangle $ and $\kinica\langle S \rangle = \Ckinica\langle S \rangle$.

\end{defn}
Note that, when $ S = \emptyset$, $ \mathcal F^{(\CCC, \omega)}\langle S \rangle  = \Cinica$ (and likewise $\Ckinica\langle S \rangle  = \Ckinica$) 
 is just the initial ($\modR$-) $(\CCC, \omega)$-coloured circuit algebra.  

Circuit algebras, like operads, 
admit \textit{presentations} in terms of \textit{generators} and \textit{relations} (see \cite[Remark~2.6]{DHR21}): A ($\modR$-) circuit algebra $\alg = (A, \alpha)$ may be obtained as a quotient of the free ($\modR$-) circuit algebra $ \mathcal F^{(\CCC, \omega)}\langle A \rangle $ (or $\Ckinica  \langle A \rangle $  on its underlying symmetric graded set $A$. 

In the remainder of this paper, we will always take $R = \Bbbk$, a field of characteristic 0 and so $\modR = \vect$.

\begin{ex}\label{ex. tangles}Let $T_4 = \{  \xover, \xunder  \}$ and $T_n = \emptyset$ for $n \neq 4$. Then, $\mathcal F\langle T \rangle  (n) = \emptyset$ when $n$ is odd and $\mathcal F\langle T \rangle  (2m) $ is the set of diagrams (planar representations) of virtual tangles on $m$ unoriented strands. 
The \emph{circuit algebra of virtual tangles} $\mathcal T$ is the quotient of $\mathcal F\langle T \rangle $ by the (ordinary) Reidemeister relations since the virtual Reidemeister relations of \cite{Kau99} are a consequence of the relations in $\BD$.
	The \emph{oriented virtual tangle circuit algebra} $\mathcal{OT}$, with generating set $  \{\oxover, \oxunder\} \subset  \mathcal{OT} (\uparrow^2\downarrow^2)$ is defined similarly. This is explained in detail in \cite[Section~4.2]{DHR20}. 
	
	More generally, we may consider circuit algebras of $(\CCC, \omega)$-coloured virtual tangles. 
	This includes, for example, circuit algebras of embedded tangles of mixed dimensions.

\end{ex}

\begin{ex}
	\label{ex. skeleton} Given a (virtual) tangle with $2m$ labelled boundary points, its \emph{skeleton} \cite{BND17} is the virtual tangle obtained by replacing each over- and under-crossing with a virtual (symmetric) crossing. This is an element of $\BD(0,2m)$. 
	
	In \cite{BND17} and \cite{DHR21}, a \emph{circuit algebra with skeleton} is a circuit algebra $\mathcal S$ indexed by Brauer diagrams rather than lists of colours. More formally, $\mathcal S$ is a circuit algebra together with a surjective circuit algebra morphism $\mathcal S\to \inica$ . Equivalently, this is a symmetric monoidal functor from the slice category $ (0 \ov \BD, \oplus, id_0)$ (see \cref{ex. slice}). Oriented circuit algebras with skeleton may be similarly defined as symmetric monoidal functors from $ (0 \ov \DiBD, \oplus, id_0)$. 

\end{ex}

\section{Circuit algebras are modular operads}\label{sec: definitions}

Modular operads \cite{HRY19a, HRY19b, Ray20} are symmetric graded objects that admit two operations -- contraction and multiplication -- such that certain axioms are satisfied. They were introduced in the study of moduli spaces of higher genus curves \cite{GK98}. 

In \cref{ssec. ca axioms}, an axiomatic (biased) description of circuit algebras is given in terms of operations on the underlying graded symmetric monoid and in \cref{ssec ca mo}, this is shown to satisfy the modular operad axioms.

	\subsection{Axioms for circuit algebras}\label{ssec. ca axioms}
	By \cref{thm. lax functor ca}, the combinatorics of a 
	$(\CCC, \omega)$-coloured 
	circuit algebra are completely described by $\CBD$. This enables an axiomatic (biased) description of circuit algebras in terms of their underlying symmetric monoids.

Let $(\CCC, \omega)$ be a palette. For $1 \leq i \leq n$ and $ \ccc = (c_1, \dots, c_n) \in \CCC^n$, let $\ccc_{\hat i} \defeq (c_1, \dots, c_{i-1}, c_{i +1},\dots, c_n) \in \CCC^{n-1}$ be the tuple obtained by ``forgetting'' $c_i$. More generally, for distinct $1 \leq j_1, \dots, j_k \leq n$, the tuple $ \ccc_{\widehat{j_1, \dots, j_k}} \in \CCC^{n-k}$ is obtained from $ \ccc$ by forgetting $ c_{j_1}, \dots, c_{j_k}$. 

Let $S =( S (\ccc))_{\ccc}$ be a $\listm (\CCC)$-graded symmetric object in $ \V$. 

		\begin{defn}
		\label{def contraction}	\label{defn: multiplication}\label{coloured mult cont} 		\label{defn: formal connected unit}
		A \emph{contraction} $\zeta$ on $S$ is a collection of morphisms 	$  \zeta^{i \ddagger j}_{\ccc} \colon S (\ccc) \to S (\ccc_{\widehat{i,j}})$ in $ \V$ defined for all $\ccc = (c_1, \dots, c_n) \in \listm (\CCC)$ such that $ c_i = \omega c_j$. 
		
			A \emph{multiplication} $\diamond$ on $S$ is 
			is a family of maps 
			\begin{equation}\label{eq. coloured mult}- \diamond^{i \ddagger j }_{\ccc , \ddd } \colon S_{\ccc } \otimes S_{\ddd} \to S_{(\ccc_{\hat i}  \ddd_{\hat j})}\end{equation}
			defined for all $\ccc \in \CCC^m, \ddd \in \CCC^n$ and $1 \leq i \leq m, 1 \leq j \leq n$ such that $c_i = \omega d_j$.
		
	\end{defn}
	
	A contraction or multiplication that commutes with the $\Sigma$-action on $S$ is 
	\emph{ $\Sigma$-equivariant. }
	
	A multiplication $\diamond$ is \emph{commutative} if, for all $\ccc, \ddd$ as above, the following diagram commutes in $\E$:\[
	\xymatrix@R = .5cm{ S_{\ccc} \otimes S_{\ddd}\ar[rr]^-{\diamond_{\ccc, \ddd}^{i\ddagger j}}\ar[d]_{\cong}&&  S_{\ccc_{\hat i}\ddd_{\hat j}} \ar[d]^{\cong}\\
		S_{\ddd} \otimes S_{\ccc} \ar[rr]_-{ \diamond_{\ddd, \ccc}^{j\ddagger i}}&& S_{\ddd_{\hat j}\ccc_{\hat i}}.}\]
	A \emph{unit} $\epsilon$ for a commutative multiplication $\diamond$ on $S$ is a choice, for each $c \in \CCC$, of distinguished morphism $\epsilon_c\colon I \to S_{c,\omega c}$ in $\V$, such that for all $\ccc = (c_1, \dots, c_n) \in \listm (\CCC)$ and $1 \leq i \leq n$ such that $c_i = c$, the compositions
	\[ S_{\ccc} \xrightarrow{\cong} I \otimes S_{\ccc} \xrightarrow {\epsilon_c \otimes id_{S_{\ccc}}} S_{c, \omega c} \otimes S_{\ccc} \xrightarrow { \diamond_{(c, \omega c), \ccc}^{2 \ddagger i}} S_{\ccc}\] and 
	\[ S_{\ccc} \xrightarrow{\cong} I \otimes S_{\ccc} \xrightarrow {\epsilon_{\omega c} \otimes id_{S_{\ccc}}} S_{c, \omega c} \otimes S_{\ccc} \xrightarrow { \diamond_{(c, \omega c), \ccc}^{1 \ddagger i}} S_{\ccc}\] are equal to the identity on $S_{\ccc}$.
	
	By \cite[Lemma~1.13]{Ray20}, if a multiplication $\diamond$ on $S$ admits a unit $\epsilon$, then it is unique.

	Observe in particular that, if $\left((S_{\ccc})_{\ccc}, \boxtimes, \eta, \zeta \right)$ 
	is a symmetric $\listm (\CCC)$-graded monoid with contraction, then $S$ admits a commutative equivariant multiplication given by:
	\begin{equation}\label{eq CA mult}
		\diamond^{i \ddagger j}_{\ccc,\ddd} \defeq \zeta^{i \ddagger m+j}_{\ccc\ddd}\circ \boxtimes_{\ccc,\ddd}\colon S(\ccc)\otimes S(\ddd) \to S(\ccc_{\hat i}\ddd_{\hat j}),
	\end{equation} 
	defined for all $\ccc = (c_1, \dots, c_m), \ddd = (d_1, \dots, d_n)$ and all $1 \leq i \leq m, 1 \leq j \leq n$ such that $c_i = \omega d_j$.

	\begin{prop}	\label{prop. product and contraction prop}
		A $\listm (\CCC)$-graded symmetric object $( A_{\ccc})_{\ccc}$ in $ \V$ describes a $\V$-circuit algebra if and only if  
		it is has the structure of a symmetric graded monoid $(A, \boxtimes , \eta)$ in $\V$ and is equipped with an equivariant contraction $\zeta$ and, for each $c \in \CCC$, a distinguished \emph{unit} morphism $\epsilon_c \colon I \to  A_{{(c, \omega c)}}$, 
		such that the following conditions (illustrated in \cref{fig. CO axioms}) hold:
		\begin{enumerate}[(c1)]
			\item the graded monoidal product $\boxtimes$ on $(A_{\ccc})_{\ccc \in \listm (\CCC)}$ is associative up to associators in $\V$;
			\item contractions commute (see also (m1) \cref{defn: Modular operad}): 
			\[\zeta^{i' \ddagger j'}_{\ccc_{\widehat{k,m}}} \circ \zeta^{k \ddagger m}_{\ccc}= \zeta^{k' \ddagger m'}_{\ccc_{\widehat{i,j}}} \circ \zeta^{i \ddagger j}_{\ccc} \colon A_{\ccc}\to A_{\ccc_{\widehat{i,j,k,m}}} \ \text{ wherever defined};\] 
			\item contraction commutes with the monoid operation:
			\[ \zeta^{ i \ddagger j}_{\ccc\ddd}\  \circ\  \boxtimes_{\ccc, \ddd}  = \boxtimes_{\ccc_{\widehat{i,j}}\oplus \ddd} \ \circ \ (\zeta^{i \ddagger j}_{\ccc} \otimes id_{\ddd}) \colon A_{\ccc}\otimes A_{\ddd} \to A_{\ccc_{\widehat{i,j}}\ddd}\]
			for all $\ddd \in \listm (\CCC)$ and $\ccc = (c_1, \dots, c_m) \in \listm (\CCC)$ with $c_i  = \omega c_j$, $1 \leq i < j \leq m$.
		\end{enumerate}

		\begin{enumerate}[(e1)]
			\item the distinguished morphisms $(\epsilon_c)_c$ provide units for the multiplication $\diamond$ induced, as in (\ref{eq CA mult}), by $\boxtimes$ and $\zeta$:
		\[\begin{array}{ll}
		id_{A_{\ccc}} & = \zeta^{2 \ddagger 2+j}_{(c,\omega c)\ccc}\circ \boxtimes_{(c,\omega c)\ccc}\circ (\epsilon_c \otimes id_{\ccc}) \\ [2pt]
		&= \zeta^{1 \ddagger 2+j}_{(c,\omega c)\ccc}\circ \boxtimes_{(c,\omega c)\ccc}\circ (\epsilon_{\omega c} \otimes id_{\ccc}).
		\end{array}\]
	\end{enumerate}
	
	A morphism of $(\CCC, \omega)$-coloured circuit algebras in $\V$ is precisely a morphism of the underlying graded symmetric objects in $\V$ that preserves the 
	monoid operation, contraction and multiplicative units.
\end{prop}

\begin{figure}[htb!]
	\centerfloat
\begin{tikzpicture}
	\node at (-6,3){	\includestandalone[width = .39\textwidth]{c1new}};
	\node at (0,3){	\includestandalone[width = .39\textwidth]{c2new}};
	\node at (6,3){	\includestandalone[width = .39\textwidth]{c3new}};
	\node at (-6,3){(c1)};
	\node at (0,3){(c2)};
	\node at (6,3){(c3)};
\end{tikzpicture}
\caption{Circuit algebras satisfy the conditions (c1)-(c3).}\label{fig. CO axioms}
\end{figure}
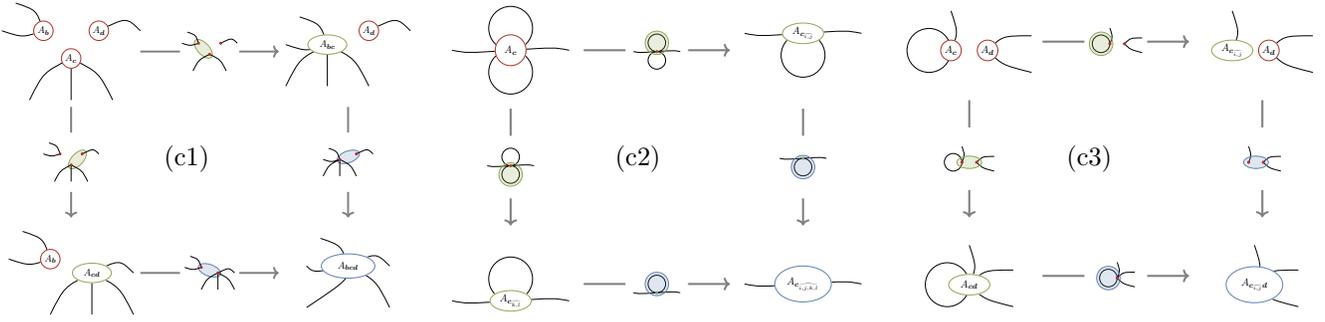

\begin{proof}

By \cref{thm. lax functor ca}, a $(\CCC, \omega)$-coloured $\V$-circuit algebra  is given by a symmetric monoidal functor $(\alg, \pi, \eta) \colon \CBD \to \V$. Since $ \Sigma^{\CCC} \subset  \CBD$, $ (\alg, \pi, \eta)$ describes a symmetric graded monoid in $\V$ and so satisfies (c1). 

Let $(\tau^{i \ddagger j}_{\ccc}, \emptyset ) \in \BDd^{(\CCC, \omega)} (\ccc , \ccc_{\widehat{i,j}})$ be the downward Brauer diagram given by
\[ \ccc \ni c_k \mapsto \left \{ \begin{array}{ll}
	c_j & k = i\\
	c_k \in \ccc_{\widehat{i,j}} & k \neq i, k \neq j .	\end{array}\right .\] This defines an equivariant contraction $\zeta$ on $(\alg(\ccc))_{\ccc}$ given by 
$\zeta^{i \ddagger j}_{\ccc}\defeq  \alg (\tau^{i \ddagger j}_{\ccc})$. The relations in $\CBD$ imply that $(\alg, \pi, \zeta)$ satisfies (c2) and (c3). (See \cite[Theorem~2.6]{LZ15} or \cite[Proposition~2.15]{Ban16}.) 
For $c \in \CCC$, define $ \epsilon_c \defeq \alg(\cup_c) \colon \alg (\varnothing_{\CCC}) \to \alg(c, \omega c)$. This satisfies (e1) by the triangle identities in $\CBD$.

Conversely, let $\left((A_{\ccc})_{\ccc}, \boxtimes, \eta, \zeta, \epsilon\right)$ satisfy (c1)-(c3) and (e1). Then $(A, \boxtimes, \eta) $ describes a symmetric monoidal functor $\tilde \alg \colon\Sigma^{\CCC} \to \V$.  
By \cite[Theorem~2.6]{LZ15} or \cite[Proposition~2.15]{Ban16}, there is a unique symmetric monoidal functor $ \alg \colon \CBD\to \V$ such that $ \alg    = \tilde \alg $ on $\Sigma^{\CCC}$ and, for all $c \in \CCC$, 
\[\alg (\cap_c) = \zeta_{(\omega c,c)^{1 \ddagger 2}}\colon A_{(\omega c, c)} \to A_{\varnothing_{\CCC}}\] and
\[\alg(\cup_c) \circ \eta  = \epsilon_c \colon I \to A_{(c, \omega c)}.\]

The final statement -- that morphisms of circuit algebras are morphisms of graded symmetric monoids preserving these maps -- is immediate.
\end{proof}

Observe that, in the proof of \cref{prop. product and contraction prop}, the cap morphisms $ \cap_c$ in $\CBD$ induce contractions while the units for the multiplication $\diamond$ are induced by cup morphisms $\cup_c$. In particular, a lax monoidal functor $ \mathcal B \colon \CBDd \to \X$ is, equivalently, a symmetric graded monoid with contraction satisfying (c1)-(c3) but without a unit for the induced multiplication. This motivates the following:
\begin{defn}
\label{def nonunital ca}
A \emph{($\CCC, \omega)$-coloured) nonunital $\V$-circuit algebra} is a symmetric lax monoidal functor $\alg \colon \CBDd \to \V$. 
\end{defn}
 Equivalently, these are algebras over the operad $\mathsf d\CWD$ of \emph{downward $(\CCC, \omega)$-coloured wiring diagrams}. 
\begin{rmk}
	\label{rmk nonunital CA}
	By \cite{SS15}, nonunital monochrome circuit algebras describe algebras in the category of representations of $O_\infty$ and $Sp_{\infty}$. Algebras in the category of representations of $GL_\infty$ are described by nonunital monochrome oriented circuit algebras. See also Sections \ref{ssec. representations of BD}~\&~\ref{ssec nonunital CA}.

\end{rmk}

	\subsection{Circuit algebras and modular operads} \label{subs. MO and Comp}\label{ssec ca mo}

	As usual, let $(\CCC, \omega)$ be an involutive palette and $(\V, \otimes, I)$ a symmetric monoidal category. 
	
	\begin{defn}\label{defn: Modular operad}
		A $(\CCC, \omega)$-coloured \emph{modular operad} with values in $\V$ is a $\listm(\CCC)$-graded symmetric object $S = (S_{\ccc})_{\ccc \in \CCC}$ together with a unital multiplication $(\diamond, \epsilon)$, 
		and a contraction $\zeta$, such that the following axioms are satisfied: 
		\begin{enumerate}[(m1)]
			\item Multiplication is associative: \\
			For all $\bbb = (b_1, \dots, b_{n_1}), \ccc = (c_1, \dots, c_{n_2}), \ddd = (d_1, \dots, d_{n_3})\in \listm(\CCC)$ and all $ 1 \leq i \leq n_1$, $1 \leq j,k \leq n_2$ with $j \neq k$ and $1 \leq m \leq n_3$ such that $b_i = \omega c_j$ and $c_k= \omega d_m$, the following square commutes:

			{\centerfloat
			\begin{minipage}{.4\textwidth}
				\[
			\xymatrix{
				S_{\bbb} \otimes S_{\ccc} \otimes S_{\ddd} 
				\ar[rr]^-{\diamond^{i \ddagger j}_{\bbb,\ccc} \otimes id_{S_{\ddd}}}
				\ar[dd]_{id_{S_{\bbb}} \otimes \diamond_{\ccc,\ddd}^{k \ddagger m}} &&
				S_{\bbb_{\hat i}\ccc_{\hat j}}\otimes S_{\ddd}
				\ar[dd]^{\diamond_ {\bbb_{\hat i}\ccc_{\hat j}, \ddd}^{k' \ddagger m}}\\
				&{}&\\
				S_{\bbb}\otimes S_{\ccc_{\hat k}\ddd_{\hat m}}
				\ar[rr]_-{\diamond_{\bbb, \ccc_{\hat k}\ddd_{\hat m}}^{i \ddagger j'}} &&
				S_{\bbb_{\hat i}\ccc_{\widehat{j,k}\ddd_{\hat m} }}. }
			\]
			\end{minipage}
			\begin{minipage}{.5\textwidth}

	\newcommand{\TL}
	{ \draw[thick]
		(0,0)--(2,2)
		(0,0)--(-2,2)
		(0,0)..controls(-2,-1)..(-3,-3)
		(0,0)..controls (2,-1)..(3,-3)
		(2,2)..controls(4,3)..(5,2)
		(-2,2)..controls(-2.5,3.5)..(-4,4)
		(-2,2)..controls(-3.5,1.5) .. (-5,2);
		\filldraw[white]
		(0,0) circle (20pt)
		(2,2) circle (20pt)
		(-2,2) circle (20pt);
		\draw[red]
		(0,0) circle (20pt)
		(2,2) circle (20pt)
		(-2,2) circle (20pt);
		
		\node at (-2,2){ {$S_{\bbb}$}};
		\node at (0,0){ {$S_{\ccc}$}};
		\node at (2,2){ {$S_{\ddd}$}};
		
	}
	
	\newcommand{\skeletonone}{ \draw[thick]
		(0,0)--(2,2)
		(0,0)--(-2,2)
		(0,0)..controls(-2,-1)..(-3,-3)
		(0,0)..controls (2,-1)..(3,-3)
		(2,2)..controls(4,3)..(5,2)
		(-2,2)..controls(-2.5,3.5)..(-4,4)
		(-2,2)..controls(-3.5,1.5) .. (-5,2);
		\filldraw[red]
		(0,0) circle (5pt)
		(2,2) circle (5pt)
		(-2,2) circle (5pt);}
	
	\newcommand{\skeletontwo}{ \draw[thick]
		(-1,1)--(2,2)
		(-1,1)..controls(-3,0)..(-4,-2)
		(-1,1)..controls (1,0)..(2,-2)
		(2,2)..controls(4,3)..(5,2)
		(-1,1)..controls(-1.5,2.5)..(-3,3)
		(-1,1)..controls(-2.5,.5) .. (-4,1);
		\filldraw[red]
		(-1,1) circle (5pt)
		(2,2) circle (5pt);}

\newcommand{\TR}{ \draw[thick]
	(-1,1)--(2,2)
	(-1,1)..controls(-3,0)..(-4,-2)
	(-1,1)..controls (1,0)..(2,-2)
	(2,2)..controls(4,3)..(5,2)
	(-1,1)..controls(-1.5,2.5)..(-3,3)
	(-1,1)..controls(-2.5,.5) .. (-4,1);
	\filldraw[white]
	(-1,1) ellipse (40pt and 20 pt)
	(2,2) circle (20pt);
	\draw  [ draw=green]
	(-1,1) ellipse (40 pt and 20 pt);
	\node at (-1,1){  $S_{\bbb_{\hat i}\ccc_{\hat j}}$};
	\draw[red]
	(2,2) circle (20pt);
	\node at (2,2){ {$S_{\ddd}$}};}

\newcommand{\skeletonthree}{ \draw[thick]
	(1,1)--(-2,2)
	(1,1)..controls(-1,0)..(-2,-2)
	(1,1)..controls (3,0)..(4,-2)
	(1,1)..controls(3,2)..(4,1)
	(-2,2)..controls(-2.5,3.5)..(-4,4)
	(-2,2)..controls(-3.5,1.5) .. (-5,2);
	\filldraw[red]
	(1,1) circle (5pt)
	(-2,2) circle (5pt);}

\newcommand{\BL}{ \draw[thick]
	(1,1)--(-2,2)
	(1,1)..controls(-1,0)..(-2,-2)
	(1,1)..controls (3,0)..(4,-2)
	(1,1)..controls(3,2)..(4,1)
	(-2,2)..controls(-2.5,3.5)..(-4,4)
	(-2,2)..controls(-3.5,1.5) .. (-5,2);
	\filldraw[white]
	(1,1) ellipse (40 pt and 20 pt)
	(-2,2) circle (20pt);
	\draw[green]
	(1,1) ellipse (40 pt and 20 pt);
	\node at (1,1){ $S_{\ccc_{\hat k}\ddd_{\hat m}}$};
	\draw[red]
	(-2,2) circle (20pt);
	\node at (-2,2){ {$S_{\bbb}$}};
	
}

\newcommand{\BR}{ \draw[thick]
	(0,1)..controls(0,0)..(-3,-2)
	(0,1)..controls (2,0)..(3,-2)
	(0,1)..controls(2,2)..(3,1)
	(0,1)..controls(-.5,2.5)..(-2,3)
	(0,1)..controls(-2.5,.5) .. (-3,1);
	\filldraw[white]
	(0,1) ellipse (55pt and 25pt);
	\draw[blue]
	(0,1) ellipse (55pt and 25pt);
	\node at (0,1){ {$S_{\bbb_{\hat i}\ccc_{\widehat{j,k}}\ddd_{\hat m}}$}};}

\newcommand{\oneA}{ 
	
			\draw  [draw=green,  fill=green, fill opacity=0.2, rotate = -47] 
			(-1.5,-.15) ellipse (60 pt and 30 pt);
			
			\draw[thick]
			(0,0)--(2,2)
			(0,0)--(-2,2)
			(0,0)..controls(-2,-1)..(-3,-3)
			(0,0)..controls (2,-1)..(3,-3)
			(2,2)..controls(4,3)..(5,2)
			(-2,2)..controls(-2.5,3.5)..(-4,4)
			(-2,2)..controls(-3.5,1.5) .. (-5,2);
			\filldraw[red]
			(0,0) circle (5pt)
			(2,2) circle (5pt)
			(-2,2) circle (5pt);
		}

	\newcommand{\oneB}{ 
	
		\draw  [draw=green,  fill=green, fill opacity=0.2, rotate = 47] 
		(1.5,-.15) ellipse (60 pt and 30 pt);
	
		\skeletonone;

	}
	
	\newcommand{\oneC}{ 
	
		\draw  [draw=blue,  fill=blue, fill opacity=0.2, rotate = 20] 
		(1,1.2) ellipse (60 pt and 30 pt);

		\skeletontwo;

	}
	
	\newcommand{\oneD}{ 
		
		\draw  [draw=blue,  fill=blue, fill opacity=0.2, rotate = -20] 
		(-1,1.2) ellipse (60 pt and 30 pt);

		\skeletonthree;

	}

	\newcommand{\oldC}{ \draw  [ draw=green, dashed,  fill=green, fill opacity=0.2]
		plot[smooth, tension=0.7, ] 
		coordinates {(0,-1.5)(-1.5,-1)(-2.5,0)(-3.5,1)(-3.5,2)(-3,3)(-2,3.5)(-1,3.5)(0,3)(1,3.5)(2,3.5)(3,3)(3.5,2)(3.5,1)(2.5,0)(1.5,-1)(0,-1.5)};

		\skeletonone;

		\node at (0,.5){\scriptsize{$\phi_2$}};
		\node at (-1.5,2){\scriptsize{$\phi_1$}};
		\node at (2,2.5){\scriptsize{$\phi_3$}};
		
	}

	\begin{tikzpicture} 
		
		\node at(-4,0){ 
			\begin{tikzpicture}[scale = 0.5]\TL\end{tikzpicture}};
		
		\draw  [->, line width = 2, draw= gray](-1.5,0)--(3.5,0);
		\filldraw[white](1,0) circle (30pt);
		\node at (1,0){
			\begin{tikzpicture}[scale = 0.2]\oneA\end{tikzpicture}
		};
		\draw  [->, line width = 2, draw= gray](-4,-2)--(-4,-6);
		\filldraw[white](-4,-4) circle (30pt);
		\node at (-4,-4){
			\begin{tikzpicture}[scale = 0.2]\oneB\end{tikzpicture}
		};
		
		\node at (6,0){ 
			\begin{tikzpicture}[scale = 0.5]\TR\end{tikzpicture}};
		\draw  [->, line width = 2, draw= gray](6,-2)--(6,-6);
		\filldraw[white](6,-4) circle (30pt);
		\node at (6,-4){
			\begin{tikzpicture}[scale = 0.2]\oneC\end{tikzpicture}
		};
		
		\node at (-4,-8){ 
			\begin{tikzpicture}[scale = 0.5]\BL\end{tikzpicture}};
		\draw  [->, line width = 2, draw= gray](-1.5,-8)--(3.5,-8);
		\filldraw[white](1,-8) circle (30pt);
		\node at (1,-8){
			\begin{tikzpicture}[scale = 0.2]\oneD\end{tikzpicture}
		};
		\node at (6,-8){ 
			\begin{tikzpicture}[scale = 0.5]\BR\end{tikzpicture}};

	\end{tikzpicture}

		\end{minipage}}
	
			\item Contractions commute (see (c2), \cref{prop. product and contraction prop} and \cref{fig. CO axioms})
			\item Multiplication and contraction commute:  \\
			For all $\ccc = (c_1, \dots, c_{n_1}), \ddd = (d_1, \dots, d_{n_1})\in \listm(\CCC)$ and all distinct $ 1 \leq i,j,k \leq n_1$, $1 \leq m \leq n_2$ such that $c_i = \omega c_j$ and $c_k= \omega d_m$, the following square commutes:
			
				{\centerfloat
				\begin{minipage}{.4\textwidth}
					\[
				\xymatrix{
					S_{\ccc } \otimes S_{\ddd}
					\ar[rr]^-{\zeta_{\ccc}^{i \ddagger j} \otimes id_{S_{\ddd}}}
					\ar[dd]_ {\diamond^{k \ddagger m}_{\ccc\ddd}} &&
					S_{\ccc_{\widehat{i,j}}} \otimes S_{\ddd} 
					\ar[dd]^{\diamond_{\ccc_{\widehat{i,j}},\ddd}^{k'\ddagger m}}\\
					&{}&\\
					S_{\ccc_{\hat k}\ddd_{\hat m}}
					\ar[rr]_-{\zeta_{\ccc_{\hat k}\ddd_{\hat m}}^{i' \ddagger j'}} &&
					S_{\ccc_{\widehat{i,j,k}} \ddd_{\hat m} }.}\] 
				\end{minipage}
				\begin{minipage}{.6\textwidth}
					
	\newcommand{\skeletonthree}{ \draw[thick]
		(-1,0) circle (1cm)
		(-1,0)--(2,0)
		(2,0)..controls(3,1)..(4,1)
		(2,0)..controls(3,-1)..(4,-1.5);
		
		\filldraw[red]
		(0,0) circle (5pt)
		(2,0) circle (5pt);}
	
	\newcommand{\oneA}{ \draw[thick]
		(-1.5,0) circle (1.5cm)
		(0,0)..controls (.5, 1.2)..(0, 2.7)
		(0,0)--(2.5,0)
		(2.5,0)..controls(3.5,1)..(5.5,1)
		(2.5,0)..controls(3.5,-1)..(5.5,-1.5);
		\draw  [draw=green,  fill=green, fill opacity=0.2] 
		(-1.5,0) circle (2cm);
		\filldraw[red]
		(0,0) circle (5pt)
		(2.5,0) circle (5pt);}
	
	\newcommand{\oneB}{ \draw[thick]
		(-1.5,0) circle (1.5cm)
		(0,0)..controls (.5, 1.2)..(0, 2.7)
		(0,0)--(2.5,0)
		(2.5,0)..controls(3.5,1)..(5.5,1)
		(2.5,0)..controls(3.5,-1)..(5.5,-1.5);
		\draw  [draw=green,  fill=green, fill opacity=0.2] 
		(1.25,0) ellipse (60 pt and 30 pt);
		\filldraw[red]
		(0,0) circle (5pt)
		(2.5,0) circle (5pt);}
	
	\newcommand{\oneC}{ \draw[thick]
		(0,0)..controls (.5, 1.2)..(0, 2.7)
		(0,0)--(2.5,0)
		(2.5,0)..controls(3.5,1)..(5.5,1)
		(2.5,0)..controls(3.5,-1)..(5.5,-1.5);
		\draw  [draw=blue,  fill=blue, fill opacity=0.2] 
		(1.25,0) ellipse (60 pt and 30 pt);
		\filldraw[red]
		(0,0) circle (5pt)
		(2.5,0) circle (5pt);}
	
	\newcommand{\oneD}{ 
		\draw[thick]
		(1,0) circle (1.5cm)
		(2.5,0)..controls (3, 1.2)..(2.5, 2.7)
		(2.5,0)..controls(3.5,1)..(5.5,1)
		(2.5,0)..controls(3.5,-1)..(5.5,-1.5);
		\draw  [draw=blue,  fill=blue, fill opacity=0.2] 
		(1,0) circle (2cm);
		\filldraw[red]
		(2.5,0) circle (5pt);
		
	}

	\newcommand{\TL}{ \draw[thick]
		(-1.5,0) circle (1.5cm)
		(0,0)..controls (.5, 1.2)..(0, 2.7)
		(0,0)--(2.5,0)
		(2.5,0)..controls(3.5,1)..(5.5,1)
		(2.5,0)..controls(3.5,-1)..(5.5,-1.5);
		\filldraw[white]
		(0,0) circle (20pt)
		(2.5,0) circle (20pt);
		\draw[red]
		(0,0) circle (20pt)
		(2.5,0) circle (20pt);
		\node at (0,0){ {$S_{\ccc}$}};
		\node at (2.5,0){ {$S_{\ddd}$}};
	}
	
	\newcommand{\TR}{ \draw[thick]
		(0,0)..controls (.5, 1.2)..(0, 2.7)
		(0,0)--(2.5,0)
		(2.5,0)..controls(3.5,1)..(5.5,1)
		(2.5,0)..controls(3.5,-1)..(5.5,-1.5);
		\filldraw[white]
		(2.5,0) circle (20pt)
		(0,0) ellipse (40 pt and 20 pt);
		\draw[red]
		(2.5,0) circle (20pt);
		\node at (2.5,0){ {$S_{\ddd}$}};
		\draw[green]
		(0,0) ellipse (40 pt and 20 pt);
		\node at (0,0){ {$S_{\ccc_{\widehat{i,j}}}$}};
	}
	
	\newcommand{\BL}{ \draw[thick]
		(1,0) circle (1.5cm)
		(2.5,0)..controls (3, 1.2)..(2.5, 2.7)
		(2.5,0)..controls(3.5,1)..(5.5,1)
		(2.5,0)..controls(3.5,-1)..(5.5,-1.5);
		\filldraw[white]
		(2.5,0) ellipse (40 pt and 20 pt);
		
		\draw[green]
		(2.5,0) ellipse (40 pt and 20 pt);
		\node at (2.5,0){ $S_{\ccc_{\hat k}\ddd_{\hat m}}$};
	}
	
	\newcommand{\BR}{ \draw[thick]
		(1,0) circle (1.5cm)
		(2.5,0)..controls (3, 1.2)..(2.5, 2.7)
		(2.5,0)..controls(3.5,1)..(5.5,1)
		(2.5,0)..controls(3.5,-1)..(5.5,-1.5);
		\filldraw[white]
		(2.5,0) ellipse (40 pt and 20 pt);
		
		\draw[green]
		(2.5,0) ellipse (40 pt and 20 pt);
		\node at (2.5,0){ $S_{\ccc_{\hat k}\ddd_{\hat m}}$};
	}

\newcommand{\threeA}{ \draw  [ draw=black, dashed,  fill=gray, fill opacity=0.1]
	plot[smooth, tension=0.7, ] 
	coordinates {(-1,2)(2,2)(3.5,1)(3.5,-1)(0,-2)(-1,-2)(-3,-1)(-3,1)(-1,2)};
	
	\draw  [draw=green,  fill=green, fill opacity=0.2] 
	(-1,0) circle (1.5cm)
	(2,0) circle (1 cm);
	\skeletonthree;
	
	\draw[draw = green, dashed]
	(-2.5,0)--(-4.2,0);
	\node at (-5,0){\scriptsize{ $\overline \zeta^{A[\underline c]}_{{w},{x}}(\phi)$}};
	\draw[draw = green, dashed]
	(3,0)--(5,0);
	\node at (5.7,0){\scriptsize{ $\eta A\psi$}};
	\node at (0,-3.5){ {$\overline \diamond ^{TA[\underline c^{\widehat{ {w}, {x}}}], TA[\underline d]}_{{y},z}\left (\overline \zeta^{A[\underline c]}_{{w},{x}}(\phi), \eta A \psi\right) $}};
	
	\node at (-.5,0){\scriptsize{$\phi$}};
	\node at (2.7,0){\scriptsize{$\psi$}};
}

\newcommand{\threeB}{\draw  [ draw=black, dashed,  fill=gray, fill opacity=0.1]
	plot[smooth, tension=0.7, ] 
	coordinates {(-1,2)(2,2)(3.5,1)(3.5,-1)(0,-2)(-1,-2)(-3,-1)(-3,1)(-1,2)};

	\draw  [draw=green,  fill=green, fill opacity=0.2] 
	plot[smooth, tension=0.7, ] 
	coordinates {(2,1)(3,0)(2,-1)(0,-1)(-1,0)(0,1)(2,1)};
	\skeletonthree;
	
	\draw[draw = green, dashed]
	(1,1)--(1,2.5);
	\node at (1,3){\scriptsize{ $\overline \diamond ^{A[\underline c], A[\underline d]}_{{y},z} (\phi  ,\psi )$}};

	\node at (-.5,0){\scriptsize{$\phi$}};
	\node at (2.7,0){\scriptsize{$\psi$}};
	
	\node at (0,-3.5){ {$\overline \zeta^{TA[\underline c^{\hat {y}} \underline d^{\hat z}]}_{{w},{x}}
			\left ( 
			\overline \diamond ^{A[\underline c], A[\underline d]}_{{y},z} (\phi  ,\psi )
			\right )$}};
	
}

\newcommand{\threeC}{\draw  [ draw=green, dashed,  fill=green, fill opacity=0.2]
	plot[smooth, tension=0.7, ] 
	coordinates {(-1,2)(2,2)(3.5,1)(3.5,-1)(0,-2)(-1,-2)(-3,-1)(-3,1)(-1,2)};

	\skeletonthree;

	\node at (-.5,0){\scriptsize{$\phi$}};
	\node at (2.7,0){\scriptsize{$\psi$}};
	
}

\begin{tikzpicture} 
	
	\node at(-4,0){ 
		\begin{tikzpicture}[scale = 0.5]\TL\end{tikzpicture}};
	
	\draw  [->, line width = 2, draw= gray](-1.5,0)--(3.5,0);
	\filldraw[white](1,0) circle (30pt);
	\node at (1,0){
		\begin{tikzpicture}[scale = 0.2]\oneA\end{tikzpicture}
	};
	\draw  [->, line width = 2, draw= gray](-4,-2)--(-4,-6);
	\filldraw[white](-4,-4) circle (30pt);
	\node at (-4,-4){
		\begin{tikzpicture}[scale = 0.2]\oneB\end{tikzpicture}
	};
	
	\node at (6,0){ 
		\begin{tikzpicture}[scale = 0.5]\TR\end{tikzpicture}};
	\draw  [->, line width = 2, draw= gray](6,-2)--(6,-6);
	\filldraw[white](6,-4) circle (30pt);
	\node at (6,-4){
		\begin{tikzpicture}[scale = 0.2]\oneC\end{tikzpicture}
	};
	
	\node at (-4,-8){ 
		\begin{tikzpicture}[scale = 0.5]\BL\end{tikzpicture}};
	\draw  [->, line width = 2, draw= gray](-1.5,-8)--(3.5,-8);
	\filldraw[white](1,-8) circle (30pt);
	\node at (1,-8){
		\begin{tikzpicture}[scale = 0.2]\oneD\end{tikzpicture}
	};
	\node at (6,-8){ 
		\begin{tikzpicture}[scale = 0.5]\BR\end{tikzpicture}};

\end{tikzpicture}

			\end{minipage}}
		 
		 	\vskip 4 ex

			\item  ``Parallel multiplication'' of pairs is well-defined: \\	For all $\ccc = (c_1, \dots, c_{n_1}), \ddd = (d_1, \dots, d_{n_1})\in \listm(\CCC)$ and all distinct $ 1 \leq i,j \leq n_1$ and distinct $1 \leq k,m \leq n_2$ such that $c_i = \omega d_k$ and $c_j= \omega d_m$, the following square commutes:
			
				{\centerfloat
				\begin{minipage}{.4\textwidth}
					\[
				\xymatrix{
					S_{\ccc} \otimes S_{\ddd}
					\ar[rr]^-{ \diamond^{i \ddagger k}_{\ccc\ddd}}
					\ar[dd]_{ \diamond^{j \ddagger m}_{\ccc\ddd}}&&
					S_{\ccc_{\hat i}\ddd_{\hat k}}
					\ar[dd]^{\zeta_{\ccc_{\hat i}\ddd_{\hat k}}^{j' \ddagger m'}}\\&{}&\\
					S_{\ccc_{\hat j}\ddd_{\hat m}}\ar[rr]_-{\zeta_{\ccc_{\hat j}\ddd_{\hat m}}^{i' \ddagger k'}}&&
					S_{\ccc_{\widehat{i,j}}  \ddd _{\widehat{k,m}}}.}\]
				\end{minipage}
				\begin{minipage}{.6\textwidth}
					
	\newcommand{\skeletonfour}{ \draw[thick]
		
		(0,0)..controls (-1.5,2)and(3.5,2)..(2,0)
		(0,0)..controls (-1.5,-2)and(3.5,-2)..(2,0)
		(0,0) ..controls (-1,-.5) ..(-3,0)
		(2,0)..controls(3,1)..(4,1)
		(2,0)..controls(3,-1)..(4,-1.5);
		
		\filldraw[red]
		(0,0) circle (5pt) 
		(2,0) circle (5pt);
		
		\node at (0,.4){\tiny w};
		\node at (2,.4){\tiny x};
		\node at (0,-.4){\tiny y};
		\node at (2,-.4){\tiny z};}
	
	\newcommand{\oneA}{ \draw[thick]
		
		(0,0)..controls (-1.5,2)and(3.5,2)..(2,0)
		(0,0)..controls (-1.5,-2)and(3.5,-2)..(2,0)
		(0,0) ..controls (-1,-.5) ..(-3,0)
		(2,0)..controls(3,1)..(4,1)
		(2,0)..controls(3,-1)..(4,-1.5);
		
		\draw  [draw=green,  fill=green, fill opacity=0.2] 
		(1,.6) ellipse (60 pt and 35 pt);
		\filldraw[red]
		(0,0) circle (5pt) 
		(2,0) circle (5pt);
	}
	
	\newcommand{\oneB}{ \draw[thick]
		
		(0,0)..controls (-1.5,2)and(3.5,2)..(2,0)
		(0,0)..controls (-1.5,-2)and(3.5,-2)..(2,0)
		(0,0) ..controls (-1,-.5) ..(-3,0)
		(2,0)..controls(3,1)..(4,1)
		(2,0)..controls(3,-1)..(4,-1.5);
		
		\draw  [draw=green,  fill=green, fill opacity=0.2] 
		(1,-.6) ellipse (60 pt and 35 pt);
		\filldraw[red]
		(0,0) circle (5pt) 
		(2,0) circle (5pt);
	}
	
	\newcommand{\oneC}{ \draw[thick]
		
		(2,0)..controls (.5,-2)and(3.5,-2)..(2,0)
		(2,0) ..controls (1,-.5) ..(-1,.5)
		(2,0)..controls(3,1)..(4,1)
		(2,0)..controls(3,-1)..(4,-1.5);
		
		\draw  [draw=blue,  fill=blue, fill opacity=0.2] 
		(2,-.8) ellipse (25 pt and 40 pt);
		\filldraw[green]
		(2,0) circle (5pt);}
	
	\newcommand{\oneD}{ \draw[thick]
		
		(2,0)..controls (.5,2)and(3.5,2)..(2,0)
		(2,0) ..controls (1,-.5) ..(-1,.5)
		(2,0)..controls(3,1)..(4,1)
		(2,0)..controls(3,-1)..(4,-1.5);
		
		\draw  [draw=blue,  fill=blue, fill opacity=0.2] 
		(2,.8) ellipse (25 pt and 40 pt);
		\filldraw[green]
		(2,0) circle (5pt);}
	
	\newcommand{\TL}{ \draw[thick]
		
		(0,0)..controls (-1.5,2)and(3.5,2)..(2,0)
		(0,0)..controls (-1.5,-2)and(3.5,-2)..(2,0)
		(0,0) ..controls (-1,-.5) ..(-3,0)
		(2,0)..controls(3,1)..(4,1)
		(2,0)..controls(3,-1)..(4,-1.5);
		
		\filldraw[white]
		(0,0) circle (20pt) 
		(2,0) circle (20pt);
		\draw[red]
		(0,0) circle (20pt) 
		(2,0) circle (20pt);
		\node at (0,0) { {$S_{\ccc}$}};
		\node at (2,0) { {$S_{\ddd}$}};

	}
	
	\newcommand{\TR}{ \draw[thick]
		
		(2,0)..controls (.5,-2)and(3.5,-2)..(2,0)
		(2,0) ..controls (1,-.5) ..(-1,.5)
		(2,0)..controls(3,1)..(4,1)
		(2,0)..controls(3,-1)..(4,-1.5);
		
		\filldraw[white]
		(2,0) ellipse (40 pt and 20pt);
		\draw[green]
		(2,0) ellipse (40 pt and 20pt);
		\node at (2,0) { {$S_{\ccc_{\hat i}\ddd_{\hat k}}$}};

	}
	
	\newcommand{\BL}{ \draw[thick]
		
		(2,0)..controls (.5,2)and(3.5,2)..(2,0)
		(2,0) ..controls (1,-.5) ..(-1,.5)
		(2,0)..controls(3,1)..(4,1)
		(2,0)..controls(3,-1)..(4,-1.5);
		
		\filldraw[white]
		(2,0) ellipse (40 pt and 20pt);
		\draw[green]
		(2,0) ellipse (40 pt and 20pt);
		\node at (2,0) { {$S_{\ccc_{\hat j}\ddd_{\hat m}}$}};

	}

	\newcommand{\BR}{ \draw[thick]
		
		(2,0) ..controls (1,-.5) ..(-1,.5)
		(2,0)..controls(3,1)..(4,1)
		(2,0)..controls(3,-1)..(4,-1.5);
		
		\filldraw[white]
		(2,0) ellipse (50 pt and 30pt);
		\draw[blue]
		(2,0) ellipse (50 pt and 30pt);
		\node at (2,0) { {$S_{\ccc_{\widehat {i,j}}\ddd_{\widehat {k,m}}}$}};

	}
	
	\newcommand{\fourA}{ \draw  [ draw=black, dashed,  fill=gray, fill opacity=0.1]
		plot[smooth, tension=0.7, ] 
		coordinates {(3,0)(3,-2)(1,-3)(-1,-2)(-1,2)(1,3)(3,2)(3,0)};
		
		\draw  [draw=green,  fill=green, fill opacity=0.2] 
		(1,0.7) circle (1.5cm)
		;
		\skeletonfour;

		\draw[draw = green, dashed]
		(-.5, 0.7)--(-2,.7);
		\node at (-3.5,.7){\scriptsize{ $\overline \diamond^{A[\underline c], A[\underline d]}_{{x},{z}} ( \phi ,\psi)$}};
		\node at (0.3,-3.7){ {$
				\overline \zeta^{TA[\underline c^{\hat {w}}\underline d^{\hat {y} }]}_{{x},{z}}\left(
				\overline \diamond^{A[\underline c], A[\underline d]}_{{x},{z}} ( \phi ,\psi)
				\right )
				$}};
		
		\node at (-.7,0){\scriptsize{$\phi$}};
		\node at (2.7,0){\scriptsize{$\psi$}};
		
	}
	
	\newcommand{\fourB}{\draw  [ draw=black, dashed,  fill=gray, fill opacity=0.1]
		plot[smooth, tension=0.7, ] 
		coordinates {(3,0)(3,-2)(1,-3)(-1,-2)(-1,2)(1,3)(3,2)(3,0)};

		\draw  [draw=green,  fill=green, fill opacity=0.2] 
		(1,-0.7) circle (1.5cm)
		;
		
		\skeletonfour;
		
		\draw[draw = green, dashed]
		(-.5,- 0.8)--(-2,-.8);
		\node at (-3.7,-.8){\scriptsize{ $\overline \diamond^{A[\underline c], A[\underline d]}_{{x},{z}} (\phi, \psi)$}};
		
		\node at (-.7,0){\scriptsize{$\phi$}};
		\node at (2.7,0){\scriptsize{$\psi$}};  
		\node at (.3,-3.7){ {$
				\overline \zeta^{TA[\underline c^{\hat {x}}\underline d^{\hat {z} }]}_{{w},{y}}
				\left (
				\overline \diamond^{A[\underline c], A[\underline d]}_{{w},{y}}(\phi , \psi)\right)$}};
		
	}
	
	\newcommand{\fourC}{\draw  [ draw=green, dashed,  fill=green, fill opacity=0.2]
		plot[smooth, tension=0.7, ] 
		coordinates {(3,0)(3,-2)(1,-3)(-1,-2)(-1,2)(1,3)(3,2)(3,0)};

		\skeletonfour;

		\node at (-.7,0){\scriptsize{$\phi$}};
		\node at (2.7,0){\scriptsize{$\psi$}};
	}

\begin{tikzpicture} 
	
	\node at(-4,0){ 
		\begin{tikzpicture}[scale = 0.5]\TL\end{tikzpicture}};
	
	\draw  [->, line width = 2, draw= gray](-1.5,0)--(3.5,0);
	\filldraw[white](1,0) circle (30pt);
	\node at (1,0){
		\begin{tikzpicture}[scale = 0.2]\oneA\end{tikzpicture}
	};
	\draw  [->, line width = 2, draw= gray](-4,-2)--(-4,-6);
	\filldraw[white](-4,-4) circle (30pt);
	\node at (-4,-4){
		\begin{tikzpicture}[scale = 0.2]\oneB\end{tikzpicture}
	};
	
	\node at (6,0){ 
		\begin{tikzpicture}[scale = 0.5]\TR\end{tikzpicture}};
	\draw  [->, line width = 2, draw= gray](6,-2)--(6,-6);
	\filldraw[white](6,-4) circle (30pt);
	\node at (6,-4){
		\begin{tikzpicture}[scale = 0.2]\oneC\end{tikzpicture}
	};
	
	\node at (-4,-8){ 
		\begin{tikzpicture}[scale = 0.5]\BL\end{tikzpicture}};
	\draw  [->, line width = 2, draw= gray](-1.5,-8)--(3.5,-8);
	\filldraw[white](1,-8) circle (30pt);
	\node at (1,-8){
		\begin{tikzpicture}[scale = 0.2]\oneD\end{tikzpicture}
	};
	\node at (6,-8){ 
		\begin{tikzpicture}[scale = 0.5]\BR\end{tikzpicture}};

\end{tikzpicture}

			\end{minipage}}

		\end{enumerate}
	\vskip 4 ex
	
Morphisms in the category $\VCMO$ of \emph{$(\CCC, \omega)$-coloured modular operads with values in $\V$} are morphisms of the underlying symmetric graded objects that preserve multiplication, contraction and units.

	Symmetric graded objects with multiplication and contraction satisfying (m1)-(m4) but without a unit for the multiplication are called \emph{nonunital modular operads}. The category of $(\CCC, \omega)$-coloured nonunital modular operads and levelwise maps that preserve multiplication and contraction is denoted $\nuVCMO$.

	\end{defn}
	
		\begin{rmk}\label{rmk different definitions of modular operads}
		This paper considers (coloured) modular operads and circuit algebras, \emph{enriched} in a symmetric monoidal category $\V$, in the sense of \cite{HRY19a, HRY19b}. In particular, their definition is relative to a fixed palette $(\CCC, \omega)$, which can be thought of as the set of objects.

		In \cite[Section~3]{RayCA2}, modular operads are defined \textit{internal} to a category $\E$ with sufficient (co)limits. Under this definition, which is based on \cite{JK11} and follows the construction of \cite{Ray20}, the object set is replaced with an involutive object \emph{object} in $\E$. 
	The two versions coincide (up to equivalence) in $\Set$. 
	\end{rmk}

	The assignment $(\CCC, \omega) \mapsto \V\mathdash\CCA$ defines a $\Cat$-valued presheaf $ca_\V$ on the palette category ${\pal}$: a morphism $\phi \colon (\CCC, \omega) \to (\CCC', \omega')$ in $\pal$ induces a strict symmetric monoidal functor $\CBD \to \CBD[(\CCC', \omega')]$, and hence $ \alg' \in \CCA[(\CCC', \omega')]_\V $ may be pulled back to 
	a $(\CCC, \omega )$-coloured circuit algebra $\phi^* \mathcal A' \in \V\mathdash\CCA$.
	
	For a symmetric monoidal category $(\V, \otimes, I)$, let $\V\mathdash \bigCA$ be the category of all $\V$-circuit algebras: objects are pairs $((\CCC, \omega), \mathcal A)$ of a palette $(\CCC, \omega)$ and a $(\CCC, \omega)$-coloured $\V$-circuit algebra $\mathcal A$, and morphisms $((\CCC, \omega), \mathcal A) \to ((\CCC', \omega'), \mathcal A') $ are pairs $(\phi, \gamma )$ where $\phi \colon \CCC\to \CCC'$ satisfies $\phi \omega = \omega' \phi$ and $ \gamma \colon \phi^* \mathcal A' \to \mathcal A$. 
When $\V = \Set$, 
	 write $\bigCA\defeq \V\mathdash \bigCA$. 
	
	The categories $\V\mathdash \nubigCA$ of all nonunital $\V$-circuit algebras, and $\VMO$ (and $\VnuMO$) of all (nonunital) $\V$-modular operads are defined similarly.

	\begin{rmk} \label{rmk not operad alg}
		Note that $\V\mathdash \bigCA$ is not a category of algebras for some single operad since the operad composition in each $\CWD$ is dependent on $(\CCC, \omega)$.
		However, when $\V  = \Set$, $ \bigCA$ can be obtained as a category of algebras for a monad. In this construction -- based on \cite{JK11, Ray20} -- the palette $(\CCC, \omega)$ is just part of the data of any given object. More generally, if $\mathsf E$ is a symmetric monoidal category with all finite limits, then it is possible to construct a monad whose algebras are circuit algebras \emph{internal} to $\mathsf E$ with palettes replaced by involutive objects in $\mathsf E$. This is described in detail in \cite{RayCA2}.
	\end{rmk}

	\begin{ex}

		\label{ex oriented ca} A morphism $\alg \to \oWD$ 
		in $\bigCA$ pulls back to an orientation on $\alg$. Hence, by \cref{prop: wheeled prop oriented brauer}, the category $\bigWP$ of wheeled props (of all colours) in $\V$  is equivalent to the slice category $ \bigOCA\simeq \bigCA\ov \oWD$. 
		
		More generally, a morphism of palettes $(\CCC, \omega) \to \dipal$ induces an orientation on $ (\CCC, \omega)$. Objects of the category $\V\mathdash \bigWP$ of wheeled props (of any colour) with values in $ \V$ are equivalent to pairs  $( \theta, \alg)$ with $\alg$ a $(\CCC, \omega)$-coloured $\V$-circuit algebra and $ \theta \colon (\CCC, \omega) \to \dipal$ a morphism of palettes. Morphisms in $\V\mathdash \bigWP$ are described by morphisms on the underlying circuit algebras that preserve the orientation on palettes. 
	
	\end{ex}

\begin{prop}\label{prop CA MO} There are canonical inclusions of categories 
	\[ \xymatrix{
\V\mathdash \bigCA	 \ar[rr]\ar[d]&	&\VMO \ar[d]\\
	\V\mathdash \nubigCA	\ar[rr]&& \VnuMO.
	}\]
	
\end{prop}
\begin{proof}
	Since multiplicative units are unique, the vertical inclusions are full and induced by simply forgetting units. 
	
	Let $(\alg, \pi, \eta)\colon \CBDd \to \V$ define a nonunital $(\CCC, \omega)$-coloured circuit algebra. By \cref{prop. product and contraction prop}, $\alg$ admits a contraction $\zeta$ such that $(\alg, \pi, \eta, \zeta)$ satisfies (c1)-(c3). 
	
Since (m2) coincides with (c2), and $\alg$ satisfies (e1), it is only necessary to check that $(\alg, \zeta, \diamond )$ satisfies (m1), (m3), (m4).

Let $\bbb = (b_1, \dots, b_{n_1}), \ccc = (c_1, \dots, c_{n_2}), \ddd = (d_1, \dots, d_{n_3})\in \listm(\CCC)$ and all $ 1 \leq i \leq n_1$, $1 \leq j,k \leq n_2$ with $j \neq k$ and $1 \leq m \leq n_3$ such that $b_i = \omega c_j$ and $c_k= \omega d_m$. 
The composition 
\[\xymatrix{ \alg (\bbb)\otimes \alg (\ccc)\otimes \alg (\ddd) 
	\ar[rrr]^-{\diamond_{\bbb,\ccc}^{i \ddagger j} \otimes id_{\ddd} }&&&
	 \alg(\bbb_{\hat i}\ccc_{\hat j}) \otimes \alg (\ddd)
	 \ar[rrr]^-{\diamond_{\bbb_{\hat i} \ccc_{\hat j}, \ddd}^{k'\ddagger m}}&&&
	  \alg (\bbb_{\hat i}\ccc_{\widehat{j,k}}\ddd_{\hat m}) }\]
	  is given by 
	  \[   \zeta^{k'\ddagger m'}_{\bbb_{\hat i}\ccc_{\hat j}\ddd}\circ \pi_{\bbb_{\hat i}\ccc_{\hat j}, \ddd}     \circ \left(  \zeta^{i \ddagger (n_1+j)}_{\bbb\ccc}\circ \pi_{\bbb,\ccc} \otimes id_{\alg (\ddd)}\right)\]
	  where $ k'   = (n_1 -1) + k, m'  = (n_1 + n_2-2)+m$.
	By (c1)-(c3) this is 
	\begin{multline*}
		\zeta^{k'\ddagger m'}_{\bbb_{\hat i}\ccc_{\hat j}\ddd}\circ \pi_{\bbb_{\hat i}\ccc_{\hat j}, \ddd}     \circ \left(  \zeta^{i \ddagger (n_1+j)}_{\bbb\ccc}\circ \pi_{\bbb,\ccc} \otimes id_{\alg (\ddd)}\right )
		\stackrel{(c3)}{=}	
		\zeta^{k'\ddagger m'}_{\bbb_{\hat i}\ccc_{\hat j}\ddd}\circ \zeta^{i \ddagger (n_1+j)}_{\bbb\ccc\ddd}\circ \pi_{\bbb\ccc, \ddd}     \circ \left(  \pi_{\bbb,\ccc} \otimes id_{\alg (\ddd)}\right )\\[5pt]
			\stackrel{(c1)}{=}	
		\zeta^{k'\ddagger m'}_{\bbb_{\hat i}\ccc_{\hat j}\ddd}\circ \zeta^{i \ddagger (n_1+j)}_{\bbb\ccc\ddd}\circ \pi_{\bbb,\ccc \ddd}     \circ \left( id_{\alg (\bbb)} \otimes  \pi_{\ccc, \ddd} \right )
			\stackrel{(c2)}{=}	
		\zeta^{i'\ddagger j'}_{\bbb\ccc_{\hat k}\ddd_{\hat m }}\circ \zeta^{(n_1 + k) \ddagger (n_1+ n_2 + m)}_{\bbb\ccc\ddd}\circ \pi_{\bbb,\ccc \ddd}     \circ \left( id_{\alg (\bbb)} \otimes  \pi_{\ccc, \ddd} \right )
		\\[5pt] 	\stackrel{(c3)}{=}	
			\zeta^{i'\ddagger j'}_{\bbb\ccc_{\hat k}\ddd_{\hat m }}\circ  \pi_{\bbb,\ccc_{\hat k} \ddd_{\hat m}} \circ  \left( id_{\alg (\bbb)} \otimes \zeta^{k \ddagger m}_{\ccc\ddd} \circ \pi_{\ccc, \ddd} \right ). \hskip 3cm
		 \end{multline*}
	  And this is precisely the composition 
	  \[\xymatrix{ \alg (\bbb)\otimes \alg (\ccc)\otimes \alg (\ddd) 
	  	\ar[rrr]^-{id_{\bbb} \otimes \diamond_{\ccc, \ddd}^{k\ddagger m} }&&&
	  \alg (\bbb)\otimes	\alg(\ccc_{\hat k}\ddd_{\hat m}) 
	  	\ar[rrr]^-{\diamond_{\bbb, \ccc_{\hat k} \ddd_{\hat m}}^{i\ddagger j'}}&&&
	  	\alg (\bbb_{\hat i}\ccc_{\widehat{j,k}}\ddd_{\hat m}) .}\]
	  	Hence $(\alg, \diamond, \zeta)$ satisfies (m1). Axioms (m3) and (m4) follow similarly, whereby $(\alg, \pi, \eta)$ defines a modular operad. Hence, since $\diamond$ is obtained as a composition of $\zeta$ and $\pi$, this defines a functorial inclusion inclusion of categories $ \V \mathdash\nubigCA \hookrightarrow \VnuMO$.  
	  	
	  Finally, if $\alg $ extends to a functor from $ \CBD$, it admits a unital multiplication $(\diamond, \epsilon)$ with $\diamond$ defined as in (\ref{eq CA mult}) and $\epsilon$ induced by $\alg (\cup_c)$. 
\end{proof}

	\begin{rmk}
	\label{rmk properads}
	The relationship between circuit algebras and modular operads observed in \cref{prop CA MO} generalises that between wheeled props and wheeled properads (c.f.~ \cite{Val07}). 
	
	The image of a wheeled prop (viewed as a circuit algebra with oriented palette) under the forgetful functor $\V \mathdash \bigCA \to \VMO $ is its underlying  \textit{wheeled properad} (see \cite{HRY15,JY15}).

\end{rmk}

By \cref{thm. lax functor ca}, circuit algebras may be characterised categorically, as lax monoidal functors from a category of Brauer diagrams, or operadically, as algebras over an operad of wiring diagrams. By contrast, the modular operad structure is inherently operadic: modular operads cannot be described by functors from some subcategory of $\CBD$. They do, however, admit a straightforward description in terms of wiring diagrams: 

A \emph{connected wiring diagram} in $\overline f$ in $\CWD$ is one that cannot be obtained as a disjoint sum $ \overline f  = \overline{f_1} \oplus \overline{f_2}$ of non-trivial wiring diagrams as in (\ref{eq. op monoid}). Note that this notion of connectedness only makes sense in the operad $\WD$ and not in the category $\BD$. Connected wiring diagrams form a suboperad of $\WD$ (or $\CWD$) and modular operads are algebras over this suboperad of connected wiring diagrams. See \cite[Section~6]{RayCA2} for more details.

 In fact, the inclusions in \cref{prop CA MO} are the right adjoints in a square of monadic adjunctions. The left adjoints for the vertical pairs are obtained by formally adjoining units, and the left adjoints for the horizontal pairs are induced by the free graded monoid monad on the underlying symmetric graded objects. This is discussed in detail in \cite{RayCA2}.

\section{Circuit algebras and invariant theory}\label{sec. invariants}

Henceforth, unless otherwise stated, all circuit algebras will take values in the category $\vect$ of vector spaces over a field $\Bbbk$ of characteristic $0$. 

Derksen and Makam \cite{DM23} have described algebras over the finite dimensional general linear groups $GL_d$ in terms of wheeled props. The aim of this  section is to adapt their methods to provide a circuit algebra characterisation of the categories of algebras for the orthogonal and symplectic groups.

\subsection{Unital circuit algebras and finite dimensional classical groups}\label{ssec unital inv}

An action of an algebraic group $G$ on a (possibly infinite dimensional) $\Bbbk$-vector space $W$ is \emph{rational} if for all $w \in W$, there is a finite dimensional $G$-stable subspace $W_w \subset   W$ containing $w$. In other words, there is a $\Bbbk$-linear morphism $\gamma \colon W \to \Bbbk[G]\otimes W$ such that, if $ \gamma (w) = \sum _{i =1}^k f_i \otimes w_i$, then $G$ acts by $g \cdot w  =  \sum _{i =1}^k f_i(g)  w_i$.

\begin{defn}
	\label{defn G algebra}
	A \emph{$G$-algebra} is a commutative $\Bbbk$-algebra $R$ equipped with a rational action of $G$ by  $\Bbbk$-algebra automorphisms. The category of $G$-algebras and $G$-equivariant ring homomorphisms is denoted by $Alg (G)$.
\end{defn}

As in \cref{ex initial CA}, let $U$ be the initial $\vect$-valued wheeled prop and, 
for any wheeled prop $P$, let $a_P \colon U \to P$ denote the unique wheeled prop map. Note that, for all $k \geq 0$, there are distinguished morphisms $e(k) = \sum_{\sigma \in \Sigma _k} \mathrm{sgn}(\sigma)\sigma\in U(k,k)$ and $ (\bigcirc - k)\in U(0,0)$.


 \begin{thm}
 	[Derksen-Makam~`23 \cite{DM23}, Theorems~5.2~\&~7.3]\label{thm. MD wheeled props}
 
 	There is an equivalence of categories between $Alg(GL_d)$ 
 	and the category of wheeled props $P$ such that $e(d+1)$ and $ \bigcirc - d$ are in the kernel of $a_P \colon U \to P$. 
 \end{thm}

By \cref{prop: wheeled prop oriented brauer}, \cref{thm. MD wheeled props} may be restated in terms of oriented circuit algebras. Therefore, given the relationship between Theorems \ref{thm lz15} and \ref{thm Weyl FFT SFT}, it is natural to ask whether there is an undirected circuit algebra version of \cref{thm. MD wheeled props} that characterises algebras over the (finite dimensional) orthogonal and symplectic groups. 

To this end, let $V$ be a (finite) $d$-dimensional vector space equipped with a nondegenerate bilinear form $\theta \colon V \otimes V \to \Bbbk$ that is either symmetric or skew-symmetric (so $d = 2k$). As in \cref{ssec. representations of BD}, let $\delta  = d$ if $\theta$ is symmetric, and $\delta  = - k  $ if $\theta$ is skew-symmetric. If $G \subset   GL (V)$ is the isometry group of $\theta$, then for $\theta$ symmetric, $G \cong O_\delta$ and, for $\theta$ skew-symmetric $G \cong Sp_{(-\delta)}$. 
	
As in \cref{ex lz ca}, let $\mathcal I_{\theta} \subset \kinica$ be the ideal generated by $ \coev{e(|\delta| + 1)} \in \kinica (2(|\delta| +1))$ (defined in (\ref{eq coev e(k)})) and $(\bigcirc - \delta )\in \kinica (0)$. Let $\CA_{\theta}\subset \vect\mathdash\CA $ be the subcategory of (monochrome $\vect$-) circuit algebras $\alg $ such that $\mathcal I_{\theta}$ is in the kernel of the unique circuit algebra morphism  $ \uniA \colon \kinica \to \alg $.

The remainder of this section is devoted to the proof of the following theorem:	 

 
	    \begin{thm}
	  	\label{thm. CA inv}
	  	The categories $Alg(G)$ and $\CA_{\theta}$ are equivalent. 
	  \end{thm}

	 By \cref{lem universal property}, a circuit algebra $\alg$ such that $ \alg (\bigcirc ) = \delta$ factors through $\Br$. Hence, \cref{thm. CA inv} may be reformulated as the statement that $Alg(G)$ is equivalent to the category of symmetric monoidal $\vect$-functors $\Br \to \vect$ for which $e(|\delta|+1) \in \Br(|\delta|+1, |\delta|+1)$ vanishes.

	 The proof of \cref{thm. CA inv} is closely based on the proof method of 
	 \cite[Sections~5-7]{DM23} and involves showing that $Alg(G)$ and $\CA_{\theta}$ are each equivalent to a third category $\mathsf{K}_{\theta}$ that will now be described.

Recall from \cref{ex lz ca} that $\vGca \colon \BD \to \vect $ is the circuit algebra described by the unique symmetric strict monoidal functor such that $1 \mapsto V$ and $ \cap \mapsto \theta$. For any $\Bbbk$-algebra $R$, we may construct a circuit algebra $ R \otimes \vGca $ with $ ( R \otimes \vGca )(n)=  R \otimes V^{\otimes n}$ in the obvious way: 
	contraction in $\vGca  $ extends to $R \otimes \vGca$, and the monoidal product on $R \otimes \vGca$  is induced by 
	 \[ \left((\sum_i r_i \otimes v_i),(\sum_j r_j \otimes v_j)\right) \mapsto \sum_{i,j} r_ir_j \otimes ( v_i \otimes v_j). \]

	Let $\mathsf{K}_{\theta}$ be the full subcategory of $\vect\mathdash\CA $ whose objects are circuit algebras $\alg$ such that there exists a $\Bbbk$-algebra $R$ and an injective morphism of circuit algebras $\alg \hookrightarrow R \otimes \vGca$.

For $R \in Alg(G)$, 
the subspace $ \bigoplus_n (R \otimes V^{\otimes n})^G$ of $G$-invariant elements in the image of $R \otimes \vGca$ is closed under the image of $\BD$ morphisms by \cref{thm lz15}, and hence describes a circuit algebra $(R \otimes \vGca)^{G}$. And, if $\phi \colon R \to S$ is a $G$-algebra homomorphism, then the induced morphism of circuit algebras $\phi \otimes id \colon R \otimes \vGca \to S \otimes \vGca$ is $G$-equivariant. Hence, the assignment $ R\mapsto R \otimes \vGca ^G$ extends to a functor $ \Phi \colon Alg(G) \to \mathsf{K}_{\theta}$. 	
	
The construction of the converse functor $\Psi\colon \mathsf{K}_{\theta} \to Alg(G)$ is more involved.
	
	Let $R$ be a $\Bbbk$-algebra. For each $n \in \N$, the pairing $\theta$ on $V$ extends to a pairing $V^{\otimes n}\otimes V^{\otimes n} \to \Bbbk$ by 
	 \[v_1 \otimes  \dots \otimes v_n\otimes w_1 \otimes  \dots \otimes w_n \mapsto \prod_{i = 1}^{n} \theta(v_i, w_i),\]  and hence to a $\Bbbk$-algebra map $ R \otimes V^{\otimes n} \otimes V^{\otimes n}  \to R$ that 
	 will also be denoted by $\theta$.

	For any morphism $\rho \colon \alg \to R \otimes \vGca$ of circuit algebras, we may consider the subspace $T_{\rho }\subset R$ spanned by elements of the form $ \theta (\rho (a), v)$, where $ a \in \alg (n)$ if $v \in V^{\otimes n}$. This is a  $\Bbbk$-algebra since \[ \theta (\rho (a), v) \theta (\rho (b), w) = \theta (\rho(a) \otimes \rho (b), v \otimes w)\] for all $a, b \in \alg$ and $v, w \in \vGca$ such that $\theta (\rho (a), v) ,\theta (\rho (b), w)$ are defined.

	 Observe that $\rho\colon \alg \to R \otimes \vGca$ factors through the inclusion $T_\rho \otimes \vGca \hookrightarrow R \otimes \vGca$ induced by $ T_\rho \subset   R$: Namely, for any non-zero $w \in V^{\otimes n}$, let $w^* \in V^{\otimes n}$ be the element defined by $\theta (w, w^*) = 1$. Let $ a \in \alg (n)$. Since $(V, \theta)$ is an orthogonal (or symplectic) space then, for all $n$ there exists a basis $(w_i)_i$ for $V^{\otimes n}$ such that for all $a \in \alg (n)$,
	\[ \rho (a) = \sum_i \theta (\rho(a), w^*_i) \otimes w_i  \in T_{\rho} \otimes V^{\otimes n}. \]

	Assume now that $\rho \colon \alg \to R \otimes \vGca$ is, moreover, an \textit{injective} morphism of circuit algebras. Then, the subspace $T_\rho\subset   R$ has the following universal property (c.f.,~\cite[Lemma~5.1]{DM23}): 
	 \begin{lem}
	 	\label{lem. universal}
	 If $\rho \colon \alg \to R \otimes \vGca$ is injective, then, for any $\Bbbk$-algebra $S$ and morphism $\lambda \colon \alg \to S \otimes \vGca$ in $\CA$, there is a unique $\Bbbk$-algebra homomorphism $\phi\colon T_\rho \to S$ such that the following diagram commutes
	 \begin{equation}\label{eq. T universal}
	 		\xymatrix{ \alg \ar[rr]^{\rho} \ar[rrd]_{\lambda} && T_\rho \otimes \vGca \ar@{..>}[d]^{\phi \otimes id} \\ 
	 &&  S \otimes \vGca. }
	 \end{equation}

	Moreover, $T_\rho \in Alg(G)$ and the assignment $\alg \mapsto T_{\rho}$ extends to a functor $\Psi \colon \mathsf{K}_{\theta} \to Alg (G)$.
	
	 \end{lem}
	 
	 \begin{proof}

	 	Let $S$ be a $\Bbbk$-algebra and let $ \lambda \colon \alg \to S \otimes \vGca$ be a circuit algebra morphism.
	 	
	 	For all $n \in \N$, there exist $w_i \in V^{\otimes n}$ such that, for all $ a \in \alg(n)$, 
	 	\[ \rho (a) = \sum_i \theta (\rho (a), w^*_i) \otimes w_i \ \text{ and }  \lambda (a) = \sum_i \theta (\lambda (a), w^*_i) \otimes w_i. \] 
	 Since $\rho$ is injective, the elements $\phi (\theta (\rho(a), w^*_i)) \defeq\theta (\lambda(a), w^*_i) \in S$ are well-defined. This assignment extends linearly to a unique $\Bbbk$-algebra homomorphism $\phi \colon T_\rho \to S$ such that Diagram (\ref{eq. T universal}) commutes. 
	
	 	Following \cite[Proof of Lemma~5.1]{DM23}, to obtain a $G$-algebra structure on $ T_{\rho}$, let $ \gamma_{\vGca} \colon \vGca \to \Bbbk[G] \otimes \vGca$ describe the rational $G$-action on $\vGca$. By the universal property of $T_{\rho}$, there is a unique $\Bbbk$-algebra homomorphism $\mu \colon T_{\rho} \to T_{\rho} \otimes \Bbbk[G]$ such that the following diagram commutes: 
	 	\begin{equation}
	 		\label{eq. T is G alg}
	 		\xymatrix{\alg \ar[rr]^{\rho} \ar[d]_{\rho}&& T_{\rho} \otimes \vGca \ar[d]^{id \otimes \gamma_{\vGca}}\\
	 			T_{\rho} \otimes \vGca \ar@{..>}[rr]_{\mu \otimes id}&& 	T_{\rho} \otimes \Bbbk[G] \otimes \vGca.}
	 	\end{equation}
	 	
	 	In particular, $\mu$ defines a rational right action of $G$ on $T_{\rho}$ such that, if $\mu (r) = \sum_i r_i \otimes f'_i$, then $ r \cdot g = \sum_i  r_i f_i (g)$, and hence a rational left action of $G$ on $R$ by $g \cdot r = r \cdot g^{-1}$. 
	 	
	 	Finally, observe that, if $\rho \colon \alg \to R \otimes \vGca$ and $ \lambda \colon \alg \to S \otimes \vGca$ are both injective morphisms of circuit algebras, then it follows from the universal property that $\rho (a) \mapsto \lambda (a)$ induces an isomorphism $T_{\rho} \cong T_{\lambda}$. Hence, we may define $ T_{\alg} \cong T_{\rho}$ to be the limit of $T_{\rho}$ where $\rho$ varies over all injective circuit algebra morphisms of the form $\alg \hookrightarrow R \otimes \vGca$ (with $R $ a $\Bbbk$-algebra). 
	 	
	 	By (\ref{eq. T is G alg}), if $ \alg, \mathcal B \in \mathsf{K}_{\theta}$, then  $T_{\alg} , T_{\mathcal B} \in Alg(G)$ and, if $\gamma \colon \alg \to \mathcal B$ is a morphism of circuit algebras, then, by the universal property (\ref{eq. T universal}), there is a $\Bbbk$-algebra morphism $T_{\alg} \to T_{\mathcal B}$ that commutes with the $G$-algebra structure by construction. Hence, $\alg \mapsto T_{\alg}$ extends to a functor $\Psi \colon \mathsf{K}_{\theta} \to Alg (G)$. 
	 	\end{proof}

\begin{prop}
	\label{prop. AlgG equiv WG}
	
The functors $ \Phi \colon Alg (G) \leftrightarrows \mathsf{K}_{\theta} \colon \Psi$ define an equivalence of categories.
\end{prop}
\begin{proof}
	The proof follows that of \cite[Theorem~5.2]{DM23}.

To see that $ \Phi \circ \Psi $ is equivalent to the identity functor on 
${\mathsf{K}_{\theta}}$, observe first that, if $ \rho\colon \alg \to R \otimes \vGca$ is an injective morphism of circuit algebras with $R$ a $\Bbbk$-algebra, then its image $\rho (\alg)$ is invariant under the $G$-action on $T_{\rho} \otimes \vGca$: 
Namely, for $g \in G$, let $L_g$ and $R_g$ respectively define left and right multiplication by $g$ in $T_\rho$ and $\vGca$. So $G$ acts on $T_\rho \otimes \vGca$ by $g \mapsto L_g \otimes L_g = R_{g^{-1}} \otimes L_g$. By (\ref{eq. T is G alg}) above,
\[ (L_g \otimes L_g ) \circ \rho =(R_{g^-1} \otimes L_g ) \circ \rho =(id\otimes L_g ) \circ (R_{g^{-1}}\otimes id ) \circ\rho  = (id\otimes L_g ) \circ (id\otimes L_{g^{-1}}) \circ\rho = \rho. \]
So, $\rho (\alg) \subset   (T _\rho \otimes \vGca)^G$ is $G$-invariant.

To prove that $ (T _\rho \otimes \vGca)^G\subset \rho (\alg) $ and therefore $\rho (\alg)= (T _\rho \otimes \vGca)^G$, let 
 $u \in T _\rho \otimes V^{\otimes n}$. So,
 \[
 u = \sum_i \theta (\rho(a_i), v_i) \otimes w_i  \  \text{ with }a_i \in \alg(n_i), v_i \in V^{\otimes n_i} \text{ and } w_i \in V^{\otimes n}. \]
Writing $ f_i \defeq \theta(-, v_i)\otimes w_i \colon V^{\otimes n_i } \to  V^{\otimes n}$, gives $u = \sum_i f_i (\rho (a_i))$. 

The elements $\rho (a_i)$ are $G$-invariant since $\rho(\alg)\subset  (T _\rho \otimes \vGca)^G$.
So, if $u (T _\rho \otimes \vGca)^G$ is also $G$-invariant, then, by applying the Reynolds operator to $u$ and $\sum f_i (\rho (a_i))$, each $f_i$ may also be assumed to be $G$-invariant. Hence, by \cref{thm lz15}, $f_i$ is a linear combination of morphisms in the image of $\vGca \colon  {\BD} \to \vect$ whereby $ u \in \rho (\alg)$ and  $\rho (\alg) = (T _\rho \otimes \vGca)^G$.

In particular, since $\rho$ is injective, $ \alg \cong (T_{\rho} \otimes \vGca)^G$. 
It follows from the definitions of $\Phi$ and $ \Psi$ that this extends to an equivalence of functors $ \Phi \circ \Psi \simeq id_{\mathsf{K}_{\theta}}$.

 For the converse, let $R \in Alg(G)$. Let $\rho \colon (R \otimes \vGca)^G \hookrightarrow R \otimes \vGca$ denote the inclusion. This factors through $T_{\rho} \otimes \vGca$, where $T_{\rho}  \subset  R$ is a $G$-subalgebra. In particular, restricting $\iota \otimes id_{\vGca}$ to $G$-invariant subspaces gives $(T_{\rho}\otimes \vGca)^G  = (R \otimes \vGca)^G  = \Phi (R)$. 
 
 Since $T_{\rho} \cong (\Psi \circ \Phi)(R)$, we want to show that $T_{\rho}  = R$. Let $\iota \colon T_{\rho} \to R$ denote the inclusion. This is a morphism of $G$-algebras by \cref{lem. universal}. In particular, $ T_{\rho} \cong \bigoplus_W T_W$ and $ R\cong \bigoplus_W R_W$, where the sum is over all irreducible $G$-representations $W$, and $T_W \subset  T_{\rho}$ and $R_W \subset  R$ are the corresponding $W$-isotypic components of $T_{\rho}$ and $R$.  Since $\iota$ preserves $G$-subrepresentations, it follows that $\iota = \bigoplus_{W} \iota_W$ where $ \iota_W \colon T_W \to R_W$ is the restriction. 
 
 Hence, to show that $T_{\rho} = R$, it suffices to show that $\iota_W$ is an isomorphism for all irreducible representations $W$ of $G$. 
 
 Let $W \subset   V^{\otimes n}$ be an irreducible representation. Then 
$ \theta \colon (R_W \otimes W)^G \otimes W \to R$ induces isomorphisms $ (R_W \otimes W)^G \otimes W \cong R_W$ and $ (T_W \otimes W)^G \otimes W \cong T_W$, and hence there is an equivariant map $ \psi_W \colon (T_{\rho} \otimes W)^G \otimes W \to (R\otimes W)^G \otimes W$ -- of the form  $\psi_W= \tilde \psi_W \otimes id_W $ for some $\tilde \psi_W \colon  (T_{\rho} \otimes W)^G \to (R\otimes W)^G$ -- such that the following diagram commutes
\begin{equation}\label{eq isotypic}
	\xymatrix
	{ T_W \ar[rr]^{\iota_W} \ar[d]_{\cong} && T_W \ar[d]^{\cong}  \\
	(T_{\rho} \otimes W)^G \otimes W \ar[rr]_{\psi_W} &&(R\otimes W)^G \otimes W.}
\end{equation}
 
Since $\iota_W$ is injective, so is $ \tilde \psi_W$. Hence, by the universal property of $T_{\rho}$, $\tilde \psi_W$ is the restriction to $	(T_{\rho} \otimes W)^G$ of $\iota \otimes id_W$ and therefore an isomorphism. Therefore $T_W = R_W$ for all irreducible representations $W$ of $G$ whereby $R = T_{\rho} \cong( \Psi \circ \Phi) (R) $ in $Alg(G)$.

This extends, by $G$-equivariance of morphisms in $Alg(G)$ and $\mathsf{K}_{\theta}$, to an equivalence of functors $\Psi \circ \Phi \simeq id_{Alg(G)}$, and therefore the categories $Alg (G)$ and $\mathsf{K}_{\theta}$ are equivalent.
\end{proof}

It remains to prove that $ \mathsf{K}_{\theta}$ is also equivalent to $\CA_{\theta}$. As in \cite[Proposition~5.3~\&~Remark~5.4]{DM23}, this rests on the following:

\begin{lem}
	\label{lem WG closed}
 If $\alg \cong (R \otimes \vGca)^G$ for some $R \in Alg(G)$ and $\mathcal J \subset  \alg$ is a circuit algebra ideal, then there exists an ideal $J \subset  R$ such that $ \mathcal J  = (J \otimes \vGca)^G$.

	Furthermore, if $ \alg \in \mathsf{K}_{\theta}$ and $ \phi \colon  \alg \to \mathcal B$ is a morphism of circuit algebras, then $\phi(\alg) \in \mathsf{K}_{\theta}$. 
\end{lem}

\begin{proof}

Let $R $ be a $G$-algebra and $\alg  = (R \otimes \vGca)^G$ its image under $\Phi$. Let $ \rho \colon \mathcal J  \subset   \alg$ be the inclusion of a circuit algebra ideal and $J \defeq T_{\rho} \subset R$. By the proof of \cref{prop. AlgG equiv WG}, $\mathcal J \cong (J \otimes \vGca)^G$ as circuit algebras. To show that $J$ is an ideal of $R$, let $r \in R$ and let $u = \sum_j \theta (\rho (\beta_j ), w_j)$ -- with $ \beta_j \in \mathcal J(m_i)$ and $w_j\in V^{\otimes m_j}$ -- be an element of $J$. By the proof of \cref{prop. AlgG equiv WG}, $r= \sum_i \theta (\alpha_i, v_i)$ for some $ \alpha _i \in \alg (n_i)$, $v_i \in V^{\otimes n_i}$. Hence,
\[ru = \sum_{i,j} \theta (\alpha_i, v_i) \theta (\rho (\beta_j ), w_j) = \sum_{i,j} \theta (\alpha_i \otimes \rho (\beta_j), v_i \otimes w_j).\]
Since $ \mathcal J \subset  \alg $ is a circuit algebra ideal, $\alpha_i \otimes \rho (\beta_j) \in \mathcal J$ for all $i,j$, and therefore $ru \in J$, whereby $J$ is an ideal of $R$.

For the second statement, let $\phi \colon \alg = (R \otimes \vGca)^G\to \mathcal B$ be a morphism of circuit algebras with kernel $ \mathcal J  \subset   \alg$. So, there is an isomorphism $\phi (\alg) \cong \alg / \mathcal J$ of circuit algebras.

Let $\iota \colon J \hookrightarrow R$ denote the inclusion of the ideal $J$ such that $ \mathcal J \cong (J \otimes \vGca)^G$, and let $q \colon R \to R/J$ be the quotient. The inclusion $\mathcal J \hookrightarrow \alg$ is given by the restriction to $(J \otimes \vGca)^G$ of $ \iota \otimes id_{\vGca}$ .

Then the following diagram -- where the vertical arrows are inclusions -- commutes:
\begin{equation}
	\xymatrix{ 0 \ar[rr] && \mathcal J \ar[rr]^{ \iota \otimes id_{\vGca}|_{\mathcal J}} \ar[d]&& \alg  \ar[rr]\ar[d]&& \alg / \mathcal J \ar[rr]\ar@{..>}[d]&&0 \\
		0 \ar[rr] &&J \otimes \vGca \ar[rr]_{\iota \otimes id_{\vGca}}&& R \otimes \vGca \ar[rr]_{q \otimes id_{\vGca}}&& R/J \otimes \vGca \ar[rr] &&0.
	}
\end{equation}
 It follows that $\alg/ \mathcal J$ is isomorphic to the image of the restriction to $\alg =(R\otimes \vGca)^G$  of the quotient $q \otimes id_{\vGca} \colon R \otimes \vGca \to R/J \otimes \vGca$. Since $\iota\colon J \to R$, and hence also $q \colon R \to R/J$, is $G$-equivariant, so is $q \otimes id_{\vGca} $. Hence, the image of its restriction to $ \alg$ is $G$-invariant, and therefore $  \mathcal J \subset  (R/J \otimes \vGca)^G$ is in $\mathsf{K}_{\theta}$ and the lemma is proved.
\end{proof}

\begin{prop}\label{prop WG CAtheta equiv}
The categories $\mathsf{K}_{\theta}$ and $\CA_{\theta}$ are equivalent. 
\end{prop}
\begin{proof}

 For all $G$-algebras $R$, since $\coev{e(|\delta| + 1)}$ and $ \bigcirc - \delta$ are in the kernel of $\uniG \colon \kinica \to \vGca$, they are in the kernel of the unique circuit algebra morphism $ z\colon \kinica \to (R \otimes \vGca )^G$. It follows, from \cref{prop. AlgG equiv WG}, that $\mathsf{K}_{\theta}$ is a full subcategory of $ \CA_{\theta}$. It therefore suffices to show that each $\alg \in \CA_{\theta}$ is equivalent to some object of $\mathsf{K}_{\theta}$. The proof follows that of \cite[Theorem~7.3]{DM23}.

Let $\alg \in \CA_{\theta}$  with underlying graded set $ A = (A_n)_n$ and let $\kinica\langle A \rangle $ be the free $\vect$-circuit algebra generated by $A$ (\cref{ex. free ca}). Then, there is a circuit algebra ideal $\mathcal I \subset  \kinica \langle A \rangle $ such that $\alg \cong \kinica \langle A \rangle  / \mathcal I$. If $\mathcal I_{\theta} \subset  \kinica \langle A \rangle $ is the ideal generated by $\coev{e(d+1)}$ and $\bigcirc - \delta$, and $ \mathcal B \defeq \kinica \langle A \rangle / \mathcal I_{\theta}$, then $\mathcal I ^\theta \subset  \mathcal I$ since $\alg \in \CA_{\theta}$. So, there exists a circuit algebra ideal $\mathcal J \subset  \mathcal B$ such that $\alg \cong \mathcal B / \mathcal J$.

To prove the proposition, it therefore suffices (by \cref{lem WG closed}) to show that there is a $\Bbbk$-algebra $R$ and an inclusion of circuit algebras $\mathcal B \subset  R \otimes \vGca $.

So, let $\{e_1, \dots, e_d\}$ be a basis for $V$ 
and, for each $n$, let $\{e_{j_1, \dots, j_n}\}_{1 \leq j_i \leq d}$ denote the induced basis for $V^{\otimes n}$. 
For each $\alpha \in A_n \subset  A$, introduce formal variables $\{a^{\alpha}_{j_1, \dots, j_n} \}_{ (j_1, \dots, j_n) \in \{ 1, \dots, d\}^n}$ and define
\[ R \defeq \Bbbk [ a^{\alpha}_{j_1, \dots, j_n}| (j_1, \dots, j_n) \in \{ 1, \dots, d\}^n, \alpha \in A_n, n \in \N].\] 

Then, the circuit algebra morphism $\rho \colon \mathcal B \to R \otimes \vGca$ given by \[ \alpha \mapsto \sum_{(j_1, \dots, j_n)} a^{\alpha}_{j_1, \dots, j_n}\otimes e_{j_1, \dots, j_n}, \ \alpha \in A_n\] is well defined since $I^{\theta}$ vanishes in $\vGca$ and therefore also in $R \otimes \vGca$. Moreover, by \cref{thm lz15}, $\rho$ is injective. Hence $\mathcal B \in \mathsf{K}_{\theta}$ and therefore, by \cref{lem WG closed}, so is $\kinica \langle A \rangle  / \mathcal I \cong \alg$.

It follows that $\CA_{\theta} \simeq \mathsf{K}_{\theta}$ as required. 
\end{proof}

\cref{thm. CA inv} follows immediately from Propositions \ref{prop. AlgG equiv WG} and \ref{prop WG CAtheta equiv}.

\begin{rmk}
	\label{rmk. ideals}
The ideals of the initial $\vect$-wheeled prop $U$ are classified in \cite[Section~4]{DM23}, and the ideals of the initial circuit algebra $\kinica$ may be similarly described. 
It is therefore natural to whether there are interesting statements, analogous to \cref{prop. AlgG equiv WG}, that consider quotients of $\kinica$ by different ideals, and whether this leads to a (partial) classification of monochrome (oriented) circuit algebras via duality results like \cref{thm. CA inv} and \cref{thm. MD wheeled props}.

\end{rmk}

\subsection{Nonunital circuit algebras and representations of infinite dimensional groups } \label{ssec nonunital CA}

As in \cref{rmk stable rep}, let $G_{\infty}= \bigcup_d G_d$ be the infinite dimensional orthogonal or symplectic group with standard representation $\bm V = \bigcup_d V_d$ and induced symmetric (or skew-symmetric) form ${\bm {\theta}} \defeq \bigcup_d \theta_d$. For all $d \geq 0$, $\theta_d$ is a nondegenerate (orthogonal or symplectic) form on the finite dimensional space $V_d$ and $\theta_{d+1}$ restricts to $\theta_d$ on $V_d \subset V_{d+1}$. An algebra $W$ over $G_{\infty}$ is, in particular, an algebra over $G_d$ for all $d \geq 1$. Hence by \cref{prop. AlgG equiv WG}, there is a compatible sequence of circuit algebras $(\alg_d)_d$ with $ \alg_d \defeq (W \otimes\vGca[{\theta_d}] )_d^{G_d}$. 

For $d \geq 0$, let $\dvGca[{\theta_d}]$ be the nonunital circuit algebra given by the restriction of $\vGca[{\theta_d}] $ to $\BDd$. Since $\dvGca[{\theta_d}] \colon \BDd \to \vect$ is a strict monoidal functor, $\vGca[{\theta_d}]$ is the unique extension of $\dvGca[{\theta_d}] $ to $\BD$.

As in \cref{rmk stable rep}, let $\bm F\colon \BDd \to \vect$ be the strict symmetric monoidal functor $1 \mapsto \bm V$, $ \cap \mapsto \bm \theta$. Then $ {\bm F} = \mathrm{colim}_d\dvGca[{\theta_d}] $ and is $G_{\infty}$-equivariant. In particular, for each $d \geq 0$, there is a morphism of nonunital circuit algebras $ p_d \colon {\bm F} \to \dvGca[{\theta_d}]$ that commutes with the actions of  $G_{\infty}$ and $G_d$ on either side.

Let $\tilde{ \mathsf{K}}_{\bm \theta}$ be the category of nonunital circuit algebras $ \tilde \alg$ for which there exists a $\Bbbk$-algebra $R$ and an inclusion of nonunital circuit algebras $\rho \colon\tilde \alg \hookrightarrow R \otimes \bm F$.

\begin{thm}
	\label{thm. Ginfty}
	There is an equivalence of categories $\widetilde {\mathsf{K}}_{\bm \theta} \simeq Alg(G_{\infty})$.
\end{thm}
\begin{proof}
Given a $G_{\infty}$-algebra $R$, we may construct the nonunital circuit algebra $(R \otimes {\bm F})^{G_{\infty}} \in \widetilde {\mathsf{K}}_{\bm \theta}$. The assignment $R \mapsto (R \otimes {\bm F})^{G_{\infty}} $ clearly extends to a functor $\tilde \Phi\colon Alg(G_{\infty}) \to \widetilde {\mathsf{K}}_{\bm \theta} $.

Conversely, let $R$ be a $\Bbbk$-algebra and let $ \rho \colon \tilde \alg \to R \otimes \bm F$ be an inclusion of nonunital circuit algebras.

 Let $\tilde T$ be the space generated by ${\bm \theta}(\rho (a), v)$ for all $a \in \tilde \alg (n), v \in {\bm F}(n)$ and all $n \in \N$. This is a $\Bbbk$-algebra as ${\bm \theta} (\rho (a), v)  {\bm \theta} (\rho (b), w) =  {\bm \theta} (\rho(a) \otimes \rho (b), v \otimes w)$. 
 
 To show that $\tilde T$ is a $G^{\infty}$-algebra, observe that, since $\dvGca[{\theta_d}]$ admits a unique extension to a circuit algebra (namely $\vGca[\theta_d]$) for all $d$, there is an increasing  sequence of circuit algebras $(R \otimes\vGca[{\theta_d}] )_{d}$.  
	 Moreover, for all $d \geq 0$, there is an injection $ \tilde \alg / ker_d \to R \otimes \dvGca[{\theta_d}]$ of nonunital circuit algebras, where $ker_d \subset \tilde \alg$ is the kernel of the nonunital circuit algebra morphism $\rho_d = p_d \circ \rho \colon \tilde \alg \to R \otimes \dvGca[{\theta_d}]$. Since $\rho$ is an injection, $ker_d$ does not depend on $R$. 
	 
	 As $\dvGca[{\theta_d}]$ admits a unique extension to a circuit algebra (namely $\vGca[\theta_d]$), so does $\tilde \alg /ker_d$. Let $\alg_d$ be the circuit algebra so defined. Then there is an inclusion $\alg_d \hookrightarrow R \otimes\vGca[{\theta_d}] $ and hence, by \cref{prop. AlgG equiv WG}, there is a $G_d$ algebra $T_d\subset R$  such that $\alg_d \cong (T_d \otimes\vGca[{\theta_d}] )^{G_d} .$
	 
	  Moreover, $T_d$ is generated by elements of the form ${\bm \theta}(\rho (a), v)$ for all $a \in \alg_d (n), v \in\vGca[{\theta_d}] (n) =  \dvGca[{\theta_d}](n) $ and all $n \in \N$ and is, up to isomorphism, independent of $R$.

	It follows, in particular, that $\tilde T = \bigcup_d T_d$ describes a filtration and hence $\tilde T $ is independent of $R$ and a $G^{\infty}$ algebra. The assignment $\tilde \alg \mapsto \tilde T$ clearly extends to a functor $\tilde \Psi \colon\widetilde {\mathsf{K}}_{\bm \theta}\to Alg (G_\infty)$ and $\tilde \Phi \colon Alg (G_\infty) \leftrightarrows\widetilde { \mathsf{K}}_{\bm \theta} \colon \tilde \Psi$ describes an equivalence of categories by \cref{prop. AlgG equiv WG} and the constructions of $\tilde \Phi, \tilde \Psi$.
\end{proof}

A directed version of \cref{thm. Ginfty}, relating nonunital wheeled props and $GL_\infty$-algebras may be obtained by similarly modifying the results of \cite{DM23}.

\bibliography{Compactbib}{}
\bibliographystyle{plain}

\end{document}